\input ./style/arxiv-general.cfg
\documentclass[aop,MSNbibl,seceqn,dvips]{arximspdf}
\makeatletter
   \@ifpackageloaded{graphicx}{}{\usepackage{graphicx}}
\makeatother

\doi{10.1214/13-AOP887} %kopijuoti is PTS
\volume{43}
\issue{3}
\pubyear{2015}
\firstpage{1202}
\lastpage{1273}

\makeatletter
\newtheorem{Lemma}{Lemma}[section]
\newtheorem{Hypothesis}{Hypothesis}
\newproclaim{Definition}[Lemma]{Definition}
\newproclaim{Remark}[Lemma]{Remark}
\newproclaim{remark}[Lemma]{Remark}
\newtheorem{Theorem}[Lemma]{Theorem}
\newtheorem{Proposition}[Lemma]{Proposition}
\newtheorem{theorem}{Theorem}[section]
\newtheorem{lemma}{Lemma}[section]
\newproclaim{definition}{Definition}[section]
\renewcommand{\emptyset}{\varnothing}
\makeatother

\begin{document}
\begin{frontmatter}

\title{Sub and supercritical stochastic quasi-geostrophic~equation\thanksref{T1}}
\runtitle{Stochastic quasi-geostrophic equation}
\thankstext{T1}{Research supported in part by NSFC (11301026) and China
Postdoctoral
Science Foundation funded project (2012M520153) and DFG through IRTG
1132 and CRC 701.}

\begin{aug}
\author[A]{\fnms{Michael} \snm{R\"{o}ckner}\ead[label=e1]{roeckner@math.uni-bielefeld.de}},
\author[B]{\fnms{Rongchan} \snm{Zhu}\corref{}\ead[label=e2]{zhurongchan@126.com}\thanksref{T2}}
\and
\author[C]{\fnms{Xiangchan} \snm{Zhu}\ead[label=e3]{zhuxiangchan@126.com}}

\runauthor{M. R\"{o}ckner, R. Zhu and X. Zhu}
\thankstext{T2}{Corresponding author.}
\affiliation{University of Bielefeld,
Beijing Institute of Technology  and Beijing Jiaotong University}
\address[A]{M. R\"{o}ckner\\
Department of Mathematics\\
University of Bielefeld\\
D-33615 Bielefeld\\
Germany\\
\printead{e1}} %adresu isvedimo komanda gale!
\address[B]{R. Zhu\\
Department of Mathematics\\
Beijing Institute of Technology\\
Beijing 100081\\
China\\
\printead{e2}}
\address[C]{X. Zhu\\
School of Science\\
Beijing Jiaotong University\\
Beijing 100044\\
China\\
\printead{e3}}
\end{aug}

% HISTORY:
\received{\smonth{2} \syear{2013}}

% ABSTRACT
%
\begin{abstract}
In this paper, we study the 2D stochastic quasi-geostrophic equation on
$\mathbb{T}^2$ for general parameter $\alpha\in(0,1)$ and
multiplicative noise. We prove the existence of weak solutions and
Markov selections for multiplicative noise for all $\alpha\in(0,1)$.
In the subcritical case $\alpha>1/2$, we prove existence and uniqueness
of (probabilistically) strong solutions. Moreover, we prove ergodicity
for the solution of the stochastic quasi-geostrophic equations in the
subcritical case
driven by possibly degenerate noise. The law of large numbers for the
solution of the stochastic quasi-geostrophic equations in the
subcritical case is also established. In the case of nondegenerate
noise and $\alpha>2/3$ in addition exponential ergodicity is proved.
\end{abstract}

% KEYWORDS

\end{frontmatter}

%s1 #&#
\section{Introduction}\label{intro}

Consider the following two-dimensional (2D) stochastic
quasi-geostrophic equation in the periodic domain $\mathbb
{T}^2=\mathbb
{R}^2/(2\pi\mathbb{Z})^2$:
%
%e1.1 #&#
%
\begin{equation}
\label{eq1.1} \frac{\partial\theta(t,
\xi)}{\partial t}=-u(t,\xi)\cdot\nabla\theta(t,\xi)-\kappa (-
\triangle )^\alpha\theta(t,\xi)+\bigl(G(\theta)\eta\bigr) (t,\xi),
\end{equation}
with initial condition
%
%e1.2 #&#
%
\begin{equation}
\label{eq1.2} \theta(0,\xi)=\theta_0(\xi),
\end{equation}
where $\theta(t,\xi)$ is a real-valued function of $\xi\in\mathbb
{T}^2$ and $t\geq0$, $0<\alpha<1, \kappa>0$ are real numbers. $u$ is
determined by $\theta$ via the following relation:
%
%e1.3 #&#
%
\begin{equation}
\label{eq1.3} u=(u_1,u_2)=(-R_2
\theta,R_1\theta)=R^\bot\theta.
\end{equation}

Here, $R_j$ is the $j$th periodic Riesz transform and $\eta(t,\xi)$ is
a Gaussian random field, white noise in time, subject to the
restrictions imposed below.
The case $\alpha=\frac{1}{2}$ is called the critical case, the case
$\alpha>\frac{1}{2}$ subcritical and the case $\alpha<\frac{1}{2}$
supercritical.

In the deterministic case ($G\equiv0$), such equations are important
models in geophysical fluid dynamics. Indeed,
they are special cases of general quasi-geostrophic approximations for
atmospheric
and oceanic fluid flows with small Rossby and Ekman numbers. These
models arise under
the assumptions of fast rotation, uniform stratification and uniform
potential vorticity. The case $\alpha=1/2$ exhibits similar features
(singularities) as the 3D Navier--Stokes equations and can therefore
serve as a
model case for the latter. For more details about the geophysical
background, see, for instance, \cite{CMT94,P87}. In the deterministic case,
this equation has been intensively investigated because of both its
mathematical importance and its background in geophysical fluid
dynamics (see, e.g., \cite{CV10,Re95,CW99,Ju04,Ju05,KNV07,CC04,KN09}
and the references therein). In the deterministic case, the
global existence of weak solutions has been obtained in \cite{Re95}
and one
most remarkable result in \cite{CV10} gives the existence of a classical
solution for $\alpha=1/2$. In \cite{KNV07}, another very important
result is
proved, namely that solutions for $\alpha=1/2$ with periodic $C^\infty$
data remain $C^\infty$ for all times.

There is another model considering a simplified geophysical fluid model
at asymptotically high rotation rate or with small Rossby number. This
geophysical model with random perturbation has been studied in
\cite{BDW98,HGH08} and the references therein. The equation is of a
different type
compared with our equation.

In this paper, we study the 2D stochastic quasi-geostrophic equation on
the torus
$\mathbb{T}^2$ for general parameter $\alpha\in(0,1)$ and for both
additive as well as multiplicative noise. Here, since the dissipation
term is not strong enough to control the nonlinear term, we have to
work in $L^p$ and to prove appropriate $L^p$-norm estimates. This leads
to considerable complications in comparison to the stochastic
Navier--Stokes equation, for example, when one wants to prove
$L^p$-norm estimates for the weak solutions (see Theorem~\ref
{the3.3}), which
are essential to obtain pathwise uniqueness, and the improved
positivity lemma to obtain uniform $L^p$-norm estimates (see Lemma~\ref
{lem5.1.4}
and Proposition~\ref{pro5.1.5}) which will be used to prove ergodicity.

\textit{Main results for general $\alpha\in(0,1)$}: We prove the
existence of weak solutions for multiplicative noise (Theorem~\ref{the3.3}).
In order to prove the existence of (probabilistically strong) solutions
and ergodicity in subsequent sections, we need $L^p$ norm estimates for
the solutions, which are
obtained using the $L^p$-It\^{o} formula proved in \cite{Kr10}. But these
$L^p$-norm estimates we cannot prove by Galerkin approximation;
instead, we use another approximation which can be seen as a piecewise
linear equation on small subintervals [see (\ref{eq3.4})]. To piece together
martingale solutions on each subinterval and to get the existence of a
martingale solution for the approximation, we first use the measurable
selection theorem to find a martingale solution measurable with respect
to the initial condition and apply a classical theorem from \cite
{SV79} (see
Theorem~\ref{the3.2}). Using an abstract result for obtaining Markov selections
from \cite{GRZ09}, we prove the existence of an a.s. Markov family in
Appendix~\hyperref[appC]{C} (Theorem~\ref{theC.5}).

\textit{Main results for the subcritical case $\alpha>1/2$}:
We obtain
pathwise uniqueness in a larger space by using $L^p$-norm estimates
(Theorem~\ref{the4.2}) and, therefore, get a (probabilistically
strong) solution
(Theorem~\ref{the4.3}) by the Yamada--Watanabe theorem. In particular, it
follows that the laws of the solutions form a Markov process.
Subsequently, in Section~\ref{sec5} we use a coupling method to study
the long
time behavior of the solution for the 2D stochastic quasi-geostrophic
equation and we obtain ergodicity, that is, the existence (Theorem~\ref
{the5.3.2})
and uniqueness (Theorem~\ref{the5.2.3}) of an invariant measure, for the
solution to the 2D stochastic quasi-geostrophic equation (in case
$\alpha>1/2$) driven by possibly degenerate noise. Furthermore, the
Markov semigroup $P_t$ converges to the unique invariant measure
polynomially fast (Theorem~\ref{the5.3.3}). Finally, we prove that a
law of
large numbers holds in our case, that is, the times averages $\frac
{1}{T}\int_0^T\psi(\theta_t)\,dt$ converge to a constant in probability
if $\psi\dvtx H^1\mapsto\mathbb{R}$ is smooth (Theorem~\ref{the5.4.1}).

We add a detailed discussion on our approach to ergodicity via
coupling, in particular, on its justification and on its relation to
other approaches in Remark~\ref{rem5.2.4} below. In this paper, we are inspired
by \cite{O08} to construct an intermediate process $\tilde{\theta}$ such
that $\theta-\tilde{\theta}$ has a strong dissipation term and $\|
\theta
(t)-\tilde{\theta}(t)\|_{H^{-1/2}}\rightarrow0$ as $t\rightarrow
\infty
$. Using this intermediate process, we can prove $E\|\theta_1(t,\theta
_0^1,\theta_0^2)-\theta_2(t,\theta_0^1,\theta_0^2)\|_{H^{-1/2}}$
converges to zero polynomially fast when time goes to infinity, where
$(\theta_1(t,\theta_0^1,\theta_0^2),\theta_2(t,\theta_0^1,\theta_0^2))$
denotes a coupling of two solutions to (\ref{eq3.1}) starting from two
different initial values $\theta_0^i\in H^1, i=1,2$. Then we can deduce
the uniqueness of invariant measures (Theorem~\ref{the5.2.3}). Also by a
suitable choice of the metrics the asymptotically strong Feller
property of the semigroup associated with the solution to the 2D
stochastic quasi-geostrophic equation is also established (Remark~\ref
{rem5.2.4}).
Here, we want to emphasize that although we consider the
semigroup in $H^1$, the convergence is in $H^{-1/2}$ norm. Moreover, we
obtain the existence of the invariant measure, which lives on $H^1$,
by using the uniform $L^p$-estimates (Theorem~\ref{the5.3.2}), which
require the
improved positivity lemma (Lemma~\ref{lem5.1.4}). Thus, we obtain
ergodicity for
the solution of the quasi-geostrophic equation in the subcritical case
(Theorem~\ref{the5.3.3}).

\textit{Additional results in the subcritical case $\alpha>2/3$}:
In
Section~\ref{sec6}, we prove the exponential convergence of the
solution under a
stronger condition on the noise and on $\alpha$.
In order to prove the exponential convergence (Theorem~\ref
{the6.3.3}), we first
show the strong Feller property of the associated semigroup
(Theorem~\ref{the6.1.2}), which follows from employing the weak-strong
uniqueness
principle in \cite{FR08} (Theorem~\ref{the6.1.3}) and the
Bismut--Elworthy--Li formula.
As the dynamics only exist in the (analytically) weak sense and
standard tools of stochastic analysis are not available, the
computations are made for an approximating cutoff dynamics,
which are equal to the original dynamics on a small random time interval.
Since in our case $\alpha<1$, it is more difficult to use the
$H^\alpha
$-norm to control the nonlinear term even though the equation is on
$\mathbb{T}^2$.
To prove the weak-strong uniqueness principle, we need some regularity
for the trajectories of the noise. Therefore, we need conditions on $G$
so that it is enough regularizing. However, in order to apply the
Bismut--Elworthy--Li formula,
we also need $G^{-1}$ to be regularizing enough. As a result, $\alpha
>2/3$ is required (see Remark~\ref{rem6.1.1} below for details).
It seems difficult to use the Kolmogorov equation method as in \cite
{DD03,DO06} (see Remark~\ref{rem6.1.1} below).

This paper is organized as follows. In Section~\ref{sec2}, we
introduce some
notation as preparation.
In Section~\ref{sec3}, we prove the existence of weak solutions for general
parameter $\alpha\in(0,1)$ and multiplicative noise.
In Section~\ref{sec4}, we prove pathwise uniqueness for all $\alpha
\in(\frac
{1}{2},1)$. Furthermore, we get the existence and uniqueness of
(probabilistically strong) solutions for multiplicative noise in the
subcritical case.
Moreover, we prove the Markov property for this unique solution. In
Section~\ref{sec5}, we use the coupling method to prove the uniqueness
of an
invariant measure in the subcritical case. Moreover, we obtain that the
semigroup $P_t$ converges to the invariant measure polynomially fast.
The law of large numbers for the solution to the 2D stochastic
quasi-geostrophic equation is also established in this section.
In Section~\ref{sec6}, for $\alpha>2/3$, and provided the noise is
nondegenerate, we prove the exponential convergence to the (unique)
invariant measure.
Appendix~\hyperref[appA]{A} is devoted to a measurability problem (see
Theorem~\ref{theA.4})
which arises in implementing the coupling method in Section~\ref
{sec5}. In
Appendix~\hyperref[appB]{B}, we prove existence of measurable selections for the
solutions to the martingale problem in Section~\ref{sec3}, and finally
Appendix~\hyperref[appC]{C}
is devoted to the existence of the corresponding Markov selection.

%s2 #&#
\section{Notations and preliminaries}\label{sec2}

In the following, we will restrict ourselves to flows which have zero
average on the torus ${\mathbb{T}^2}$, that is,
\[
\int_{\mathbb{T}^2}\theta\,d\xi=0,
\]
where $d\xi$ denotes the volume measure on ${\mathbb{T}^2}$.
Thus, (\ref{eq1.3}) can be restated as
\[
u=\biggl(-\frac{\partial\psi}{\partial\xi_2},\frac{\partial\psi
}{\partial
\xi_1}\biggr) \quad\mbox{and}\quad(-
\Delta)^{1/2}\psi=-\theta.
\]
Set $H=\{f\in L^2(\mathbb{T}^2)\dvtx\int_{\mathbb{T}^2}f \,d\xi=0\}
$ and let
$|\cdot|$ and $\langle\cdot,\cdot\rangle$ denote the norm and inner
product in $H$, respectively. $L^p(\mathbb{T}^2), p\in(0,\infty]$
denote the standard $L^p$ spaces on ${\mathbb{T}^2}$ with norm $\|
\cdot
\|_{L^p}$. On the periodic domain $\mathbb{T}^2$, $\{\sin\langle
k,\cdot\rangle_{\mathbb{R}^2}|k\in\mathbb{Z}^2_+\}\cup\{\cos
\langle
k,\cdot\rangle_{\mathbb{R}^2}|k\in\mathbb{Z}^2_-\}$ form\vspace*{1pt} an eigenbasis
of $-\Delta$ (we denote it by $\{e_k\}$). Here, $ \mathbb{Z}^2_+=\{
(k_1,k_2)\in\mathbb{Z}^2|k_2>0\}\cup\{(k_1,0)\in\mathbb
{Z}^2|k_1>0\}$,\vspace*{1pt}
$\mathbb{Z}^2_-=\{(k_1,k_2)\in\mathbb{Z}^2|{-}k\in\mathbb{Z}^2_+\}$,
$\xi
\in\mathbb{T}^2$, and the corresponding eigenvalues are $|k|^2$. For
$s>0$, define
\[
\|f\|_{H^s}^2=\sum_k
|k|^{2s}\langle f,e_k\rangle^2
\]
and let $H^s$ denote the Sobolev space of all $f\in H$ for which $\|f\|
_{H^s}$ is finite.
For $s<0$, define $H^s$ to be the dual of $H^{-s}$.
Set $\Lambda=(-\Delta)^{1/2}$. Then
\[
\|f\|_{H^s}=\bigl |\Lambda^s f\bigr |.
\]
For $s\geq0, p\in[1,+\infty]$ we use $H^{s,p}$ to denote a subspace of
$L^p(\mathbb{T}^2)$, consisting of all $f$ which can be written in the
form $f=\Lambda^{-s}g, g\in L^p(\mathbb{T}^2)$ and the $H^{s,p}$ norm
of $f$ is defined to be the $L^p$ norm of $g$, that is, $\|f\|
_{H^{s,p}}:=\|\Lambda^s f\|_{L^p}$.

By the singular integral theory of Calder\'{o}n and Zygmund (cf. \cite{St70},
Chapter~3), for any $s\geq0, p\in(1,\infty)$, there is a constant
$C_R=C_R(s,p)$, such that
%
%e2.1 #&#
%
\begin{equation}
\label{eq2.1} \bigl \|\Lambda^su\bigr \|_{L^p}\leq C_R(s,p)
\bigl \|\Lambda^s\theta\bigr \|_{L^p}.
\end{equation}

Fix $\alpha\in(0,1)$ and define the linear operator $A_\alpha\dvtx
D(A_\alpha)=H^{2\alpha}(\mathbb{T}^2)\subset H\rightarrow H$ as
$A_\alpha u:=\kappa(-\Delta)^\alpha u$. The operator $A_\alpha$ is
positive definite and self-adjoint with the same eigenbasis as that of
$-\Delta$ mentioned above. Denote the eigenvalues of $A_\alpha$ by
$0<\lambda_1\leq\lambda_2\leq\cdots$\,, and renumber the above
eigenbasis correspondingly as $e_1,e_2,\ldots$\,.

First, we recall the following important product estimates (cf. \cite{Re95},
Lem\-ma~A.4):

%le2.1 #&#
%
\begin{Lemma}\label{lem2.1}
Suppose that $s>0$ and $p\in(1,\infty)$. If $f,g\in
C^{\infty}(\mathbb{T}^2)$ then
%
%e2.2 #&#
%
\begin{equation}
\label{eq2.2} \bigl \|\Lambda^s(fg)\bigr \|_{L^p}\leq C\bigl(\|f
\|_{L^{p_1}}\bigl \|\Lambda ^sg\bigr \|_{L^{p_2}}+\|g
\|_{L^{p_3}}\bigl \|\Lambda^sf\bigr \|_{L^{p_4}}\bigr),
\end{equation}
with $p_i\in(1,\infty]$, $i=1,\ldots,4$ such that
\[
\frac{1}{p}=\frac{1}{p_1}+\frac{1}{p_2}=\frac{1}{p_3}+
\frac{1}{p_4}.
\]
\end{Lemma}

We shall use as well the following standard Sobolev inequality (cf.
\cite{St70}, Chapter~V):

%le2.2 #&#
%
\begin{Lemma}\label{lem2.2}
Suppose that $q>1$, $p\in[q,\infty)$ and
\[
\frac{1}{p}+\frac{\sigma}{2}=\frac{1}{q}.
\]
Suppose that $\Lambda^\sigma f\in L^q$, then $f\in L^p$ and there is a
constant $C_S\geq0$ independent of $f$ such that
\[
\|f\|_{L^p}\leq C_S\bigl \|\Lambda^\sigma f
\bigr \|_{L^q}.
\]
\end{Lemma}

The following commutator estimate from \cite{Ju04}, Lemma~3.1, is very
important for later use.

%le2.3 #&#
%
\begin{Lemma}[(Commutator estimates)]\label{lem2.3}
Suppose that $s>0$ and $p\in
(1,\infty)$. If $f,g\in C^\infty(\mathbb{T}^2)$, then
\[
\bigl \|\Lambda^s(fg)-f\Lambda^sg\bigr \|_{L^p}\leq C
\bigl(\|\nabla f\|_{L^{p_1}}\bigl \| \Lambda^{s-1}g\bigr \|_{L^{p_2}}+\|g
\|_{L^{p_3}}\bigl \|\Lambda^sf\bigr \|_{L^{p_4}}\bigr),
\]
with $p_i\in(1,\infty)$, $i=1,\ldots,4$ such that
\[
\frac{1}{p}=\frac{1}{p_1}+\frac{1}{p_2}=\frac{1}{p_3}+
\frac{1}{p_4}.
\]
\end{Lemma}

We will also use the following classical interpolation inequality (see,
e.g., \cite{D04}, (5.5)).

%le2.4 #&#
%
\begin{Lemma}\label{lem2.4}
For $f\in C^\infty(\mathbb{T}^2)$, we have
%
%e2.3 #&#
%
\begin{equation}
\label{eq2.3} \|f\|_{H^s}\leq C\|f\|_{H^{s_1}}^{
({s_2-s})/({s_2-s_1})}\|f
\|_{H^{s_2}}^{({s-s_1})/({s_2-s_1})},\qquad s_1<s<s_2.
\end{equation}
\end{Lemma}

%s3 #&#
\section{Weak solutions in the general case}\label{sec3}

In this section, we consider the following abstract stochastic
evolution equation in place of equations (\ref{eq1.1})--(\ref{eq1.3}):
%
%e3.1 #&#
%
\begin{equation}
\label{eq3.1} \cases{ %
d\theta(t)+A_\alpha
\theta(t)\,dt+u(t)\cdot\nabla\theta (t)\,dt= G\bigl(\theta(t)\bigr)\,dW(t),
\cr
\theta(0)=\theta_0\in H, %\end{array}
}
\end{equation}
where $u$ satisfies (\ref{eq1.3}) and $W(t)$, $t\in[0,T]$, is a
cylindrical Wiener
process in a separable Hilbert space $U$ defined on a filtered
probability space $(\Omega,\mathcal{F},\{\mathcal{F}_t\}_{t\in
[0,T]},P)$. Here, $G$ is a measurable mapping from $H^\alpha$ to
$L_2(U,H)$ ($=$ all Hilbert--Schimit operators from $U$ to $H$). Let
$f_n$, $n\in\mathbb{N}$, be an ONB of $U$.

In the following, we assume the following conditions on $G$:

\renewcommand{\theHypothesis}{G.1}
%hy1 #&#
%
\begin{Hypothesis}\label{hypG.1}
\textup{(i)} $\|G(\theta)\|^2_{L_2(U,H)}\leq\lambda_0|\theta|^2+\rho
_1|\Lambda
^\alpha\theta|^2+\rho_2$, $\theta\in H^\alpha$, for some positive real
numbers $\lambda_0$, $\rho_2$ and $\rho_1<2\kappa$. Moreover, for some
$\beta>3$, $\|G(\theta)\|^2_{L_2(U,H^{-\beta})}\leq\rho_3(|\theta
|^2+1),\theta\in H^\alpha$, for some positive real numbers $\rho_3$.

\textup{(ii)} If $\theta,\theta_n\in H^\alpha$ such that $\theta
_n\rightarrow
\theta$ in $H$, then $\lim_{n\rightarrow\infty}\|G(\theta
_n)^*(v)-G(\theta)^*(v)\|_U=0$ for all $v\in C^\infty(\mathbb{T}^2)$,
where the asterisk denotes the adjoint operator of $G(\theta)$.
\end{Hypothesis}

First, we introduce the following definition of a weak solution.

%de3.1 #&#
%
\begin{Definition}\label{def3.1}
We say that there exists a weak solution of
equation (\ref{eq3.1}) if there exists a stochastic basis $(\Omega
,\mathcal
{F},\{\mathcal{F}_t\}_{t\in[0,T]},P)$, a cylindrical Wiener process
$W$ on the space $U$ and a progressively measurable process $\theta
\dvtx[0,T]\times\Omega\rightarrow H$, such that for $P$-a.e. $\omega
\in
\Omega$,
\[
\theta(\cdot,\omega)\in L^\infty\bigl([0,T];H\bigr)\cap L^2
\bigl([0,T];H^\alpha \bigr)\cap C\bigl([0,T];H^{-\beta}\bigr),
\]
where $\beta$ in Hypothesis~\ref{hypG.1}, and such that $P$-a.s.
\begin{eqnarray*}
&&\bigl\langle\theta(t),\phi\bigr\rangle+\int_0^t
\bigl\langle A_\alpha^{1/2}\theta (s),A_\alpha^{1/2}
\phi\bigr\rangle\,ds-\int_0^t \bigl\langle u(s)
\cdot\nabla \phi ,\theta(s)\bigr\rangle\,ds
\\
&&\qquad=\langle\theta_0,\phi
\rangle+\biggl\langle\int_0^tG\bigl(\theta(s)
\bigr)\, dW(s),\phi \biggr\rangle
\end{eqnarray*}
for $t\in[0,T]$ and all $\phi\in C^1(\mathbb{T}^2)$.
\end{Definition}

\begin{Remark*}
(i) Note that, because $\operatorname{div} u=0$ for smooth functions
$\theta$ and $\psi$, we have
\[
\bigl\langle u(s)\cdot\nabla\theta(s),\psi\bigr\rangle=-\bigl\langle u(s)\cdot
\nabla\psi,\theta(s)\bigr\rangle.
\]
Thus, the integral equation in Definition~\ref{def3.1}
corresponds to equation (\ref{eq3.1}).

(ii) Note that since the solution $\theta\in L^2(0,T;H^\alpha)$ we
only need $\theta,\theta_n\in H^\alpha$ instead of $\theta,\theta
_n\in
H$ in Hypothesis~\ref{hypG.1}(ii).

(iii) A typical example satisfying Hypothesis~\ref{hypG.1}
is the following:
For $y\in U$,
\[
G(\theta)y=\sum_{k=1}^\infty
\bigl(c_k\Lambda^\alpha\theta+b_kg(\theta )
\bigr)\langle y,f_k\rangle_U,\qquad\theta\in
H^\alpha,
\]
where $g$ is continuous function on $\mathbb{R}$ of at most linear
growth and $b_k,c_k\in C^\infty(\mathbb{T}^2)$ satisfy $\sum_kc^2_k(\xi
)<2\kappa,\sum_kb^2_k(\xi)\leq M,\xi\in\mathbb{T}^2$, and $\sum_k|\Lambda^{\alpha}c_k|^2\leq M$.
\end{Remark*}

It is standard to show that under Hypothesis~\ref{hypG.1} there exists
a weak solution
to (\ref{eq3.1}) by using the Galerkin approximation. However, as
mentioned in
the \hyperref[intro]{Introduction}, we also need $L^p$ norm estimates for the solutions,
more precise that they belong to $L^p(\Omega;L^\infty
([0,T]);L^p(\mathbb
{T}^2))$, provided so do their initial values. This will be essential
to the proof of pathwise uniqueness. For this, we have to use another
approximation instead of the Galerkin approximation and the following
theorem from \cite{SV79}, Theorem~6.1.2.

Let $\Omega_0:=C([0,\infty),H^1), \Omega_0^t:=C([t,\infty),H^1)$ for
$t>0$ and $\mathcal{P}(\Omega_0)$ denote the set of all probability
measures on $(\Omega_0,\mathcal{B})$ with $\mathcal{B}$ being the Borel
$\sigma$-algebra coming from the topology of locally uniform
convergence on $\Omega_0$. Define the canonical process $x\dvtx\Omega
_0\rightarrow H^1$ as
\[
x_t(\omega)=\omega(t).
\]
Also define the $\sigma$-algebra $\mathcal{B}_t:=\sigma\{ x(s),s\leq
t\}
$ and $\mathcal{B}^t:=\sigma\{ x(s),s\geq t\}$.

%th3.2 #&#
%
\begin{Theorem}\label{the3.2}
Fix $t>0$. Let $x\mapsto Q_x$ be a mapping from $\Omega
_0$ to $\mathcal{P}(\Omega_0^t)$ such that for any $A\in\mathcal
{B}^t$, $x\mapsto Q_x(A)$ is $\mathcal{B}_t$-measurable, and for any
$x\in\Omega_0$
\[
Q_x\bigl(y\in\Omega_0^t\dvtx y(t)=x(t)
\bigr)=1.
\]
Then for any $P\in\mathcal{P}(\Omega_0)$, there exists a unique
$P\otimes_t Q\in\mathcal{P}(\Omega_0)$ such that
\[
(P\otimes_t Q) (A)=P(A),\qquad\forall A\in\mathcal{B}_t,
\]
and for $P\otimes_t Q$-almost all $x\in\Omega_0$
\[
Q_x=(P\otimes_tQ) (\cdot|\mathcal{B}_t)
(x).
\]
\end{Theorem}

Now we will prove the existence of a martingale solution under
Hypothesis~\ref{hypG.1}.

%th3.3 #&#
%
\begin{Theorem}\label{the3.3}
Let $\alpha\in(0,1)$. If $G$ satisfies Hypothesis~\ref{hypG.1},
then there exists a weak solution $(\Omega,\mathcal{F},\{\mathcal
{F}_t\}, P,W, \theta)$ to (\ref{eq3.1}). Moreover, assume that $G$ satisfies
the following condition:
\begin{longlist}[(Gp.1)]
\item[(Gp.1)] There exists some $p\in(2,\infty)$ such that for all
$\theta\in
H^\alpha\cap L^p(\mathbb{T}^2)$,
%
%e3.2 #&#
%
\begin{equation}
\int\biggl(\sum_j\bigl |G(\theta) (f_j)\bigr |^2
\biggr)^{p/2}\,d\xi\leq C\biggl(\int |\theta|^p\,d\xi+1
\biggr), \qquad\forall t>0
\end{equation}
for some constant $C:=C(p)>0$ and $\theta_0\in L^p(\mathbb{T}^2)$. Then
\[
E\sup_{t\in[0,T]}\bigl \|\theta(t)\bigr \|^p_{L^p}<
\infty.
\]
\end{longlist}
\end{Theorem}

%re3.4 #&#
%
\begin{remark}\label{rem3.4}
Typical examples for $G$ satisfying (Gp.1) have the
following form: for $\theta\in H^\alpha$
\[
G(\theta)y=\sum_{k=1}^\infty
b_k\langle y,f_k\rangle_U g(\theta ),
\qquad y\in U,
\]
where $g$ is a continuous function on $\mathbb{R}$ of at most linear
growth and $b_k$ are $C^\infty$ functions on $\mathbb{T}^2$ satisfying
$\sum_{k=1}^\infty b_k^2(\xi)\leq M$.
\end{remark}

\begin{pf*}{Proof of Theorem~\ref{the3.3}}
\textit{Step} 1: We first establish the
existence of martingale solutions of the following equation:
%
%e3.3 #&#
%
\begin{eqnarray}
\label{eq3.3} &&d\theta(t)+A_\alpha\theta(t)\,dt+w(t)\cdot\nabla\theta
(t)\,dt= k_\delta*G(\theta)\,dW(t),
\nonumber
\\[-8pt]
\\[-8pt]
&&\theta(0)=\theta_0\in H^3,
\nonumber
\end{eqnarray}
with a given smooth function $w(t)$ which satisfies $\operatorname
{div} w(t)=0$ for
all $t\in[0,T]$ and
\[
\sup_{t\in[0,T]}\bigl \|w(t)\bigr \|_{C^3(\mathbb{T}^2)}\leq C.
\]
Here, $k_\delta*G(\theta)$ means for $y\in U$, $k_\delta*G(\theta
)(y):=k_\delta*(G(\theta)(y))$, where $k_\delta$ is the periodic
Poisson kernel in $\mathbb{T}^2$ given by $\widehat{k_\delta}(\zeta
)=e^{-\delta|\zeta|}$, $\zeta\in\mathbb{Z}^2$.
By \cite{GRZ09}, Theorem~4.7, this equation has a martingale solution
$P\in
\mathcal{P}(C([0,\infty);H^1))$ with initial value $\theta_0$ in the
following sense:
\begin{enumerate}[(M3)]
\item[(M1)] $P(x(0)=\theta_0)=1$ and for any $n\in\mathbb{N}$
\begin{eqnarray*}
&&P\biggl\{x\in C\bigl([0,\infty);H^1\bigr)\dvtx \int
_0^n\bigl \|\Lambda^{2\alpha
}x(s)+w(s)\cdot
\nabla x(s)\bigr \|_{H^1}\,ds
\\
&&\phantom{\hspace*{83pt}}{}+\int_0^n
\bigl \|k_\delta*G\bigl(x(s)\bigr)\bigr \| ^2_{L_2(U;H^3)}\,ds<+\infty
\biggr\}=1.
\end{eqnarray*}
\item[(M2)] For every $e_i$, the process
\[
\bigl\langle x(t),e_i\bigr\rangle-\int^t_0
\bigl\langle-w(s)\cdot\nabla x(s)-A_\alpha x(s),e_i\bigr
\rangle\,ds
\]
is a continuous square-integrable $\mathcal{B}_t$-martingale under $P$,
whose\break quadratic variation process is
given by
\[
\int_0^t\bigl \|(k_{\delta}*G)^*\bigl(x(s)
\bigr) (e_i)\bigr \|_U^2\,ds,
\]
where the asterisk denotes the adjoint operator of $k_\delta*G(x(s))$.
\item[(M3)] For any $q\in\mathbb{N}$ there exists a continuous positive
real function $t\rightarrow C_{t,q}$ such that
\begin{eqnarray*}
&&E^P\biggl(\sup_{r\in[0,t]}\bigl |\Lambda^3
x(r)\bigr |^{2q}+\int_0^t\bigl |
\Lambda^3 x(r)\bigr |^{2q-2}\bigl |\Lambda^{\alpha+3}x(r)\bigr |^2
\,dr\biggr)
\\
&&\qquad\leq C_{t,q}\bigl(\bigl |\Lambda ^3\theta
_0\bigr |^{2q}+1\bigr),
\end{eqnarray*}
where $E^P$ denotes the expectation under $P$.
\end{enumerate}

Indeed, we only need to check conditions (C1)--(C3) in \cite{GRZ09}. The
demi-continuity condition (C1) is obvious by Hypothesis~\ref
{hypG.1}(ii) and the
linearity of the equation. For (C2), we have that for $x\in H^4$
\[
\langle-w\cdot\nabla x-A_\alpha x,x\rangle_{H^3}\leq-\kappa \bigl |
\Lambda ^{3+\alpha}x\bigr |^2+\bigl |\bigl\langle\Lambda^3(w\cdot
\nabla x),\Lambda ^{3}x\bigr\rangle\bigr |.
\]
By Lemma~\ref{lem2.3} and because $\langle w\cdot\nabla\Lambda^{3}
x,\Lambda
^{3}x\rangle=0$ for $x\in H^4$ we have that for $x\in H^4$
\begin{eqnarray*}
\bigl |\bigl\langle\Lambda^3(w\cdot\nabla x),\Lambda^{3}x\bigr
\rangle\bigr |&=&\bigl |\bigl\langle \Lambda ^3(w\cdot\nabla x)-w\cdot\nabla
\Lambda^{3}x,\Lambda^{3}x\bigr\rangle \bigr |
\\
&\leq&\| w
\|_{C^3(\mathbb{T}^2)}\bigl |\Lambda^3x\bigr |\bigl |\Lambda^{3+\alpha}x\bigr |.
\end{eqnarray*}
Thus, the coercivity condition (C2) follows from the above two
inequalities and Young's inequality.
Also by Hypothesis~\ref{hypG.1}, we have for $x\in H^4$
\[
\bigl \|k_\delta*G(x)\bigr \|_{L_2(U,H^3)}^2\leq C(\delta)\bigl \|G(x)\bigr \|
_{L_2(U,H)}^2\leq C\bigl(\bigl |\Lambda^\alpha
x\bigr |^2+1\bigr),
\]
and by Lemma~\ref{lem2.1}
\[
\|w\cdot\nabla x+A_\alpha x\|^2_{H^1}\leq2
\|A_\alpha x\|^2_{H^1}+C\| w\| ^2_{C^3(\mathbb{T}^2)}\bigl |
\Lambda^3x\bigr |^2\leq C\bigl |\Lambda^3x\bigr |^2,
\]
which implies the growth condition (C3).

\textit{Step} 2: Now we construct an approximation of (\ref{eq3.1}).

We pick a smooth $\phi\geq0$, with $\operatorname{supp}\phi\subset[1,2]$, $\int_0^\infty
\phi=1$, and for $\delta>0$ let
\[
U_\delta[\theta](t):=\int_0^\infty\phi(
\tau) \bigl(k_\delta*R^\bot \theta \bigr) (t-\delta\tau)\,d
\tau,
\]
where $k_\delta$ is the periodic Poisson kernel in $\mathbb{T}^2$ given
by $\widehat{k_\delta}(\zeta)=e^{-\delta|\zeta|},\zeta\in\mathbb
{Z}^2$, and we set $\theta(t)=0$, $t<0$. We take a sequence $\delta
_n\rightarrow0$ and consider the equation
%
%e3.4 #&#
%
\begin{equation}
\label{eq3.4} d\theta_n(t)+A_\alpha\theta_n(t)
\,dt+u_n(t)\cdot\nabla \theta_n(t)\,dt=
k_{\delta_n}*G(\theta_n)\,dW(t),
\end{equation}
with initial data $\theta_n(0)=k_{\delta_n}*\theta_0$ and
$u_n=U_{\delta
_n}[\theta_n]$.
For a fixed $n$, this is a linear equation in $\theta_n$ on each
subinterval $[t^n_k,t^n_{k+1}]$ with $t^n_k=k\delta_n$, since $u_n$ is
determined by the values of $\theta_n$ on the two previous
subintervals. By Step 1, we obtain the existence of a martingale
solution to (\ref{eq3.4}) for fixed $n$. Indeed, we obtain the martingale
solution $P_n^1\in\mathcal{P}(C([0,\infty),H^1))$ with initial
condition $k_{\delta_n}*\theta_0$ on the subinterval $[0,t^n_1]$ by
Step 1. Also, by Step 1, we get that for $x_0\in B_0$ with $B_0:=\{
x\in\Omega_0\dvtx\sup_{0\leq t\leq t_1^n}\|x(t)\|_{H^3}<\infty\}$, there
exists a $Q_{x_0}\in\mathcal{P}(C([t^n_1,t_2^n],H^1))$ satisfying the
following:
\begin{enumerate}[(M3)]
\item[(M1)] $Q_{x_0}(x(t^n_1)=x_0(t^n_1))=1$
\begin{eqnarray*}
&&Q_{x_0}\biggl\{x\in C\bigl(\bigl[t_1^n,t_2^n
\bigr];H^1\bigr)\dvtx \int_{t_1^n}^{t_2^n}\bigl \|
\Lambda^{2\alpha}x(s)+U_{\delta_n}[x_0](s)\cdot\nabla x(s)
\bigr \|_{H^1}\,ds
\\
&&\phantom{\hspace*{120pt}} {}+\int_{t_1^n}^{t_2^n}
\bigl \|k_\delta *G\bigl(x(s)\bigr)\bigr \|^2_{L_2(U;H^3)}\,ds<+\infty
\biggr\} =1.
\end{eqnarray*}
\item[(M2)] For every $e_i$, $i\in\mathbb{N}$, the process
\begin{eqnarray*}
M_i\bigl(t\wedge t^n_2,x\bigr)&:=&\bigl
\langle x\bigl(t\wedge t^n_2\bigr),e_i\bigr
\rangle -\bigl\langle x_0\bigl(t_1^n
\bigr),e_i\bigr\rangle
\\
&&{}-\int^{t\wedge t^n_2}_{t_1^n}
\bigl\langle -U_{\delta_n}[x_0](s)\cdot\nabla
x(s)-A_\alpha x,e_i \bigr\rangle\, ds,\qquad t\geq
t_1^n
\end{eqnarray*}
is a continuous square-integrable $\mathcal{B}_t$-martingale under
$Q_{x_0}$, whose\break quadratic variation process is
given by
\[
\langle M_i\rangle\bigl(t\wedge t^n_2,x
\bigr):=\int_{t_1^n}^{t\wedge t^n_2}\bigl \| (k_{\delta_n}*G)^*
\bigl(x(s)\bigr) (e_i)\bigr \|_U^2\,ds,
\]
where the asterisk denotes the adjoint operator of $k_{\delta_n}*G(x(s))$.
\item[(M3)] For any $q\in\mathbb{N}$, there exists a constant $C_{q}$
depending on\break $\sup_{t\in[0,t_1^n]}\|x_0(t)\|_{H^1}$ such that
\begin{eqnarray*}
&&E^{Q_{x_0}}\biggl(\sup_{r\in[t_1^n,t_2^n]}\bigl |\Lambda^3
x(r)\bigr |^{2q}+\int_{t_1^n}^{t_2^n}\bigl |
\Lambda^3 x(r)\bigr |^{2q-2}\bigl |\Lambda^{\alpha
+3}x(r)\bigr |^2
\,dr\biggr)
\\
&&\qquad\leq C_{q}\bigl(\bigl |\Lambda^3x_0
\bigl(t_1^n\bigr)\bigr |^{2q}+1\bigr).
\end{eqnarray*}
\end{enumerate}
Now we extend $Q_{x_0}$ to a probability measure on $C([t_1^n,+\infty
),H^1)$ by $Q_{x_0}\circ\psi^{-1}$ with $\psi
\dvtx C([t_1^n,t_2^n],H^1)\rightarrow C([t_1^n,+\infty),H^1)$ by $\psi
x(s):=x(s\wedge t_2^n)$, $s\in[t_1^n,+\infty)$.
The set of all such martingale solutions is denoted by $\mathcal
{Q}_{x_0}$. Now we can find $Q_{x_0}\in\mathcal{Q}_{x_0}$ satisfying
(M1)--(M3) such that the map $x_0\mapsto Q_{x_0}$ from $B_0$ to
$\mathcal{P}(\Omega^{t_1^n}_0)$ is measurable with respect to
$\mathcal
{B}_{t_1^n}$. This will be proved in Lemma~\ref{lemB.1} in
Appendix~\hyperref[appB]{B}.

For $x_0\in B_0^c$, define $Q_{x_0}:=\delta_{x_0|_{[t_1^n,\infty)}}$.
Thus, by Theorem~\ref{the3.2} we get that there exists $P_n^1\otimes
_{t_1^n}Q\in
\mathcal{P}(C([0,\infty),H^1))$ such that
\[
\bigl(P_n^1\otimes_{t_1^n} Q\bigr)
(A)=P_n^1(A),\qquad\forall A\in\mathcal{B}_{t_1^n},
\]
and for $P_n^1\otimes_{t_1^n} Q$-almost all $x\in\Omega_0$
\[
Q_x=\bigl(P_n^1\otimes_{t_1^n}Q
\bigr) (\cdot|\mathcal{B}_{t_1^n}) (x).
\]
Here, $Q_{x_0}$ extends to a probability measure on $C([0,\infty),H^1)$
by the following: Let $\delta_{x_0}$ be the point-mass on
$C([0,t_1^n],H^1)$ at $x_0|_{[0,t_1^n]}$, that is,
\[
\delta_{x_0}\bigl(x\in C\bigl(\bigl[0,t_1^n
\bigr],H^1\bigr)\dvtx x(t)=x_0(t),0\leq t\leq
t_1^n\bigr)=1.
\]
Define $\tilde{Q}=\delta_{x_0}\times Q_{x_0}$ on $\tilde
{X}:=C([0,t_1^n],H^1)\times C([t_1^n,\infty),H^1)$ and set $X:=\{
(x_1,x_2)\in C([0,t_1^n],H^1)\times C([t_1^n,\infty
),H^1)\dvtx x_1(t_1^n)=x_2(t_1^n)\}$. Then $X$ is a measurable subset of
$\tilde{X}$ and $\tilde{Q}(X)=1$. Then $\tilde{Q}$ can be restricted to
$X$. Finally, $\Psi\dvtx X\rightarrow C([0,\infty),H^1)$ defined by
$\Psi
((x_1,x_2))(t):=x_1(t)$, if $0\leq t\leq t_1^n$, $\Psi
((x_1,x_2))(t):=x_2(t)$, if $ t>t_1^n$, is a measurable map form $X$
onto $C([0,\infty),H^1)$. Then $\tilde{Q}|_X\circ\Psi^{-1}$ is the
desired measure, which still be denoted $Q_{x_0}$.

By (M2), we have for every $e_i, i\in\mathbb{N}$, that the process
\begin{eqnarray*}
M_i\bigl(t\wedge t^n_2,x\bigr)&=&\bigl
\langle x\bigl(t\wedge t^n_2\bigr),e_i\bigr
\rangle -\bigl\langle x_0\bigl(t_1^n
\bigr),e_i\bigr\rangle
\\
&&{}-\int^{t\wedge t^n_2}_{t_1^n}
\bigl\langle -U_{\delta_n}[x_0](s)\cdot\nabla
x(s)-A_\alpha x,e_i \bigr\rangle\,ds
\\
&=&\bigl\langle x\bigl(t\wedge t^n_2
\bigr),e_i\bigr\rangle-\bigl\langle x_0
\bigl(t_1^n\bigr),e_i\bigr\rangle
\\
&&{}-\int
^{t\wedge t^n_2}_{t_1^n}\bigl\langle-U_{\delta_n}[x](s)\cdot
\nabla x(s)-A_\alpha x,e_i \bigr\rangle\,ds
\end{eqnarray*}
is a continuous square-integrable $\mathcal{B}_t$-martingale under $Q_{x_0}$.
Thus, by \cite{SV79}, Theorem~1.2.10, we obtain for every $e_i$, $i\in
\mathbb
{N}$, that the process
\[
\bigl\langle x\bigl(t\wedge t^n_2\bigr),e_i
\bigr\rangle-\int^{t\wedge t^n_2}_0\bigl\langle
-U_{\delta_n}[x](s)\cdot\nabla x(s)-A_\alpha x,e_i
\bigr\rangle\,ds
\]
is a continuous square-integrable $\mathcal{B}_t$-martingale under
$P_n^1\otimes_{t_1^n}Q$, whose\break quadratic variation process is
given by
\[
\int_0^{t\wedge t^n_2}\bigl \|(k_{\delta_n}*G)^*\bigl(x(s)
\bigr) (e_i)\bigr \|_U^2\,ds.
\]
Thus, we construct a martingale solution $P_n^1\otimes_{t_1^n}Q\in
\mathcal{P}(C([0,\infty),H^1))$ of (\ref{eq3.4}) on $[0,t_2^n]$. Then
step by
step we can construct a martingale solution $P_n\in\mathcal
{P}(C([0,\infty),H^1))$ of (\ref{eq3.4}) on $[0,T]$ for any given $T$
in the
following sense:
\begin{enumerate}[(M3$'$)]
\item[(M1$'$)]$P_n(x(0)=k_{\delta_n}*\theta_0)=1$ and
\begin{eqnarray*}
&&P_n\biggl\{x\in C\bigl([0,+\infty);H^1\bigr)\dvtx
\int_0^T\bigl \|\Lambda^{2\alpha
}x(s)+U_{\delta_n}[x](s)
\cdot\nabla x(s)\bigr \|_{H^1}\,ds
\\
&&\phantom{\hspace*{113pt}} {}+\int_0^T\bigl \|
k_{\delta_n}*G\bigl(x(s)\bigr)\bigr \|^2_{L_2(U;H^3)}\,ds<+\infty
\biggr\}=1.
\end{eqnarray*}
\item[(M2$'$)] For every $e_i$, the process
\[
\bigl\langle x(t\wedge T),e_i\bigr\rangle-\int^{t\wedge T}_0
\bigl\langle -U_{\delta_n}[x](s)\cdot\nabla x(s)-A_\alpha
x,e_i \bigr\rangle \,ds
\]
is a continuous square-integrable $\mathcal{B}_t$-martingale under
$P_n$, whose\break quadratic variation process is
given by
\[
\int_{0}^{t\wedge T}\bigl \|(k_{\delta_n}*G)^*\bigl(x(s)
\bigr) (e_i)\bigr \|_U^2\,ds,
\]
where the asterisk denotes the adjoint operator of $k_{\delta_n}*G(x(s))$.

\item[(M3$'$)]$P_n(L_{\mathrm{loc}}^\infty([0,+\infty),H^3)\cap
\Omega_0)=1$.
\end{enumerate}

Then by the martingale representation theorem (cf. \cite{On05},
Theorem~2,
\cite{DZ92}, Theorem~8.2) we can find a new probability space $(\Omega
^n,P^n,W_n)$ and $\theta_n$ such that $(\theta_n,W_n)$ is a weak
solution of (\ref{eq3.4}) and $\theta_n$ has the same law as $P_n$.

\textit{Step} 3: Now we show that $\theta_n$ converge to the solution
of (\ref{eq3.1}).
Since we have
\[
\bigl\langle u_n(t)\cdot\nabla\theta_n(t),
\theta_n(t)\bigr\rangle=0,
\]
by It\^{o}'s formula we have
\begin{eqnarray*}
d|\theta_n|^p+p\kappa|\theta_n|^{p-2}\bigl |
\Lambda^\alpha\theta _n\bigr |^2\,dt&\leq& p|
\theta_n|^{p-2}\bigl\langle k_{\delta_n}*G(\theta
_n)\,dW_n,\theta_n\bigr\rangle
\\
&&{}+\frac{p(p-1)}{2}|\theta_n|^{p-2}
\bigl \|k_{\delta
_n}*G(\theta_n)\bigr \|^2_{L_2(U,H)}\,dt.
\end{eqnarray*}
By classical arguments, we easily show that there exist positive
constants $C_1,C_2$ independent of $n$, such that (cf. \cite{FG95}, Appendix~1) for $2\leq p<1+\frac{2\kappa}{\rho_1}$ if $\rho_1>0$ and for
$2\leq
p<\infty$ if $\rho_1=0$, the following are satisfied:
%
%e3.5 #&#
%
\begin{equation}
\label{eq3.5} E^{P^n}\Bigl(\sup_{0\leq s\leq T}\bigl |
\theta_n(s)\bigr |^p\Bigr)\leq C_1
\end{equation}
and
%
%e3.6 #&#
%
\begin{equation}
\label{eq3.6} E^{P^n} \int^T_0\bigl \|
\theta_n(s)\bigr \|_{H^\alpha}^2\,ds\leq C_2.
\end{equation}
Now we prove that the family $\mathcal{D}(\theta_n),n\in\mathbb{N}$,
is tight in $C([0,T];H^{-\beta})$, for all $\beta>3$. Here, $\mathcal
{D}(\theta_n)$ means the law of $\theta_n$. By (\ref{eq3.5}) for
each $t\in
[0,T]$, $\mathcal{D}(\theta_n(t))$ is tight on $H^{-\beta}$. Then by
Aldous' criterion in \cite{Al78}, it suffices to check that for all stopping
times $\tau_n\leq T$ and $\eta_n\rightarrow0$,
%
%e3.7 #&#
%
\begin{equation}
\label{eq3.7} \lim_nE^{P^n}\bigl \|
\theta_n(\tau_n+\eta_n)-
\theta_n(\tau _n)\bigr \|_{H^{-\beta}}=0.
\end{equation}
We have $P^n$-a.s.
\begin{eqnarray*}
\theta_n(\tau_n+\eta_n)-
\theta_n(\tau_n)&=&-\int_{\tau_n}^{\tau
_n+\eta
_n}A_\alpha
\theta_n(s)\,ds-\int_{\tau_n}^{\tau_n+\eta
_n}u_n(s)
\cdot\nabla \theta_n(s)\,ds
\\
&&{}+\int_{\tau_n}^{\tau_n+\eta_n}k_{\delta
_n}*G
\bigl(\theta_n(s)\bigr)\,dW_n(s).
\end{eqnarray*}
It is easy to obtain the following:
%
%e3.8 #&#
%
\begin{equation}
\label{eq3.8} E^{P^n}\biggl\|\int_{\tau_n}^{\tau_n+\eta_n}A_\alpha
\theta _n(s)\,ds\biggr\|_{H^{-\beta}}\leq C\eta_nE^{P^n}
\sup_{t\in[0,T]}\bigl |\theta_n(t)\bigr |.
\end{equation}
And since $H^2\subset L^\infty$, we obtain that for $v\in H^3$,
\[
\bigl |\langle u_n\cdot\nabla\theta_n,v\rangle\bigr |=\bigl |\langle
u_n\cdot\nabla v,\theta_n\rangle\bigr |\leq|
\theta_n||u_n|\|\nabla v\|_{L^\infty}\leq |
\theta_n||u_n|\|v\|_{H^3}.
\]
Since $\sup_{[0,t]}|u_n|\leq C\sup_{[0,t]}|\theta_n|$, we get that
%
%e3.9 #&#
%
\begin{equation}
\label{eq3.9} E^{P^n}\biggl\|\int_{\tau_n}^{\tau_n+\eta_n}u_n(s)
\cdot\nabla \theta_n(s)\,ds\biggr\|_{H^{-\beta}}\leq C
\eta_nE^{P^n}\sup_{t\in
[0,T]}\bigl |
\theta_n(t)\bigr |^2.
\end{equation}
In addition by Hypothesis~\ref{hypG.1}, we have
%
%e3.10 #&#
%
\begin{eqnarray}\label{eq3.10}
&& E^{P^n}\biggl\|\int_{\tau_n}^{\tau_n+\eta
_n}k_{\delta_n}*G
\bigl(\theta_n(s)\bigr)\,dW(s)\biggr\|^2_{H^{-\beta}}
\nonumber
\\
&&\qquad\leq C E^{P^n}\int_{\tau_n}^{\tau_n+\eta_n}\bigl \|G
\bigl(\theta_n(s)\bigr)\bigr \|_{L_2(U,H^{-\beta})}^2\,ds
\\
&&\qquad\leq C\eta_n\Bigl(E^{P^n}\sup_{t\in[0,T]}\bigl |
\theta _n(t)\bigr |^2+1\Bigr)\rightarrow 0\qquad\mbox{as }
\eta_n\rightarrow0.
\nonumber
\end{eqnarray}
Thus, (\ref{eq3.7}) follows by (\ref{eq3.8}), (\ref{eq3.9}) and
(\ref{eq3.10}), which implies the
tightness of $\mathcal{D}(\theta_n)$ in $C([0,T],H^{-\beta})$. This
yields that for each $\eta>0$
\[
\lim_{\delta\rightarrow0}\sup_nP^n\Bigl(
\sup_{|s-t|\leq\delta,s,t\leq
T}\bigl |\theta_n(t)-\theta_n(s)\bigr |_{H^{-\beta}}>
\eta\Bigr)=0.
\]
By this and (\ref{eq3.5}), (\ref{eq3.6}), it is easy to get that
$\mathcal{D}(\theta_n)$
is tight in $L^2([0,T];H)\cap C([0,T], H^{-\beta})$ (cf. \cite{MR05},
Lemma~2.7).
Therefore, we find a subsequence, still denoted by $\theta_n$, such
that $\mathcal{D}(\theta_n)$ converges weakly in
\[
L^2\bigl([0,T];H\bigr)\cap C\bigl([0,T], H^{-\beta}\bigr).
\]

By Skorohod's representation theorem, there exist a stochastic basis
$(\tilde{\Omega},\tilde{\mathcal{F}},\allowbreak  \{\tilde{\mathcal{F}}_t\}
_{t\in
[0,T]}, \tilde{P})$ and, on this basis, $L^2([0,T];H)\cap C([0,T],
H^{-\beta})$-valued random variables $\tilde{\theta},\tilde{\theta
}_n,n\geq1$, such that $\tilde{\theta}_n$ has the same law as
$\theta
_n$ on $L^2([0,T];\break H)\cap C([0,T], H^{-\beta})$, and $\tilde{\theta
}_n\rightarrow\tilde{\theta}$ in $L^2([0,T];H)\cap C([0,T],
H^{-\beta
})$ $\tilde{P}$-a.s. For $\tilde{\theta}_n$ we also have (\ref
{eq3.5}) and
(\ref{eq3.6}). Hence, it follows that
\[
\tilde{\theta}(\cdot,\omega)\in L^2\bigl([0,T];H^\alpha
\bigr)\cap L^\infty \bigl([0,T];H\bigr)\qquad\mbox{for } \tilde{P}
\mbox{-a.e. } \omega\in\Omega.
\]
For each $\tilde{\theta}_n$ we define $\tilde{u}_n:=U_{\delta
_n}[\tilde
{\theta}_n]$ and for each $n\geq1$ we define the process
\[
\tilde{M}_n(t):=\tilde{\theta}_n(t)-k_{\delta_n}*
\theta_0+\int_0^tA_\alpha
\tilde{\theta}_n(s)\,ds+\int_0^t
\tilde{u}_n(s)\cdot \nabla \tilde{\theta}_n(s)\,ds.
\]
In fact $\tilde{M}_n$ is a square integrable martingale with respect to
the filtration
\[
\{\mathcal{G}_n\}_t=\sigma\bigl\{\tilde{
\theta}_n(s),s\leq t\bigr\}.
\]
For all $r\leq t\in[0,T]$, all bounded continuous functions $\phi$ on
$C([0,r];H^{-\beta})\cap L^2([0,r];H)$, and all $v\in C^\infty
(\mathbb
{T}^2)$, we have
\[
\tilde{E}\bigl(\bigl\langle\tilde{M}_n(t)-\tilde{M}_n(r),v
\bigr\rangle\phi (\tilde {\theta}_n|_{[0,r]})\bigr)=0
\]
and
\[
\tilde{E}\biggl(\biggl(\bigl\langle\tilde{M}_n(t),v\bigr
\rangle^2-\bigl\langle\tilde {M}_n(r),v\bigr
\rangle^2-\int_r^t
\bigl \|(k_{\delta_n}*G)^*(\tilde{\theta }_n)v\bigr \|
^2_U\,ds\biggr)\phi(\tilde{\theta}_n|_{[0,r]})\biggr)=0.
\]
By the B--D--G inequality, we have for $1<p<\frac{1}{2}+\frac{\kappa
}{\rho
_1}$ if $\rho_1>0$ and $1<p<\infty$ if $\rho_1=0$, that
\[
\sup_n \tilde{E}\bigl |\bigl\langle\tilde{M}_n(t),v
\bigr\rangle\bigr |^{2p}\leq C\sup_n\tilde {E}\biggl(
\int_0^t\bigl \|(k_{\delta_n}*G)^*(\tilde{
\theta}_n)v \bigr \|_U^2\, ds
\biggr)^p<\infty.
\]
Since $\tilde{\theta}_n\rightarrow\tilde{\theta}$ in
$L^2(0,T;H)\cap
C(0,T, H^{-\beta})$, we also have
\[
\lim_{n\rightarrow\infty}\tilde{E}\bigl |\bigl\langle\tilde{ M}_n(t)-M(t),v
\bigr\rangle\bigr |=0
\]
and
\[
\lim_{n\rightarrow\infty}\tilde{E}\bigl |\bigl\langle\tilde {M}_n(t)-M(t),v
\bigr\rangle\bigr |^2=0,
\]
where
\[
M(t):=\tilde{\theta}(t)-\theta_0+\int_0^t
\tilde{u}\cdot\nabla \tilde {\theta}+A_\alpha\tilde{\theta}\,ds.
\]
Here, $\tilde{u}$ is defined by (\ref{eq1.3}) with $\theta$ replaced
by $\tilde
{\theta}$.
Taking the limit, we obtain that for all $r\leq t\in[0,T]$, all bounded
continuous functions $\phi$ on $C([0,r];H^{-\beta})\cap L^2([0,r];H)$,
and $v\in C^\infty(\mathbb{T}^2)$,
\[
\tilde{E}\bigl(\bigl\langle M(t)-M(r),v\bigr\rangle\phi(\tilde{
\theta}|_{[0,r]})\bigr)=0
\]
and
\[
\tilde{E}\biggl(\biggl(\bigl\langle M(t),v\bigr\rangle^2-\bigl
\langle M(r),v\bigr\rangle^2-\int_r^t
\bigl \| G(\theta)^*v\bigr \|_U^2 \,ds\biggr)\phi(\tilde{
\theta}|_{[0,r]})\biggr)=0.
\]
Thus, the existence of a weak solution for (\ref{eq3.1}) follows by the
martingale representation theorem (cf. \cite{DZ92}, Theorem~8.2, \cite{On05},
Theorem~2).

\textit{Step} 4: Now we prove the last statement. It is sufficient to prove
that
\[
E^{P^n}\sup_{t\in[0,T]}\bigl \|\theta_n(t)
\bigr \|^p_{L^p}\leq C,
\]
where $C$ is a constant independent of $n$.
We write for simplicity $\theta(t)=\theta_n(t)$, $u(t)=u_n(t)$,
$W(t)=W_n(t)$, $P=P^n$. By \cite{Kr10}, Lemma~5.1, or \cite{BVVL08},
Theorem~2.4, we have
\begin{eqnarray*}
\bigl \|\theta(t)\bigr \|_{L^p}^p&=&\|k_{\delta_n}*
\theta_0\|_{L^p}^p
\\
&&{}+\int_0^t
\biggl[-p\int_{\mathbb{T}^2}\bigl |\theta(s)\bigr |^{p-2}\theta(s) \bigl(
\Lambda ^{2\alpha
}\theta(s)+u(s)\cdot\nabla\theta(s)\bigr)\,d\xi
\\
&&\phantom{+\int_0^t
\biggl[} {}+\frac{1}{2}p(p-1)\int_{\mathbb{T}^2}\bigl |
\theta(s)\bigr |^{p-2}\biggl(\sum_j\bigl |k_{\delta_n}*G
\bigl(\theta (s)\bigr) (f_j)\bigr |^2\biggr)\,d\xi\biggr]\,ds
\\
&&{}+p\int_0^t \int_{\mathbb{T}^2}\bigl |
\theta(s)\bigr |^{p-2}\theta(s) k_{\delta
_n}*G\bigl(\theta (s)\bigr)\,d\xi
\,dW(s)
\\
&\leq&\|k_{\delta_n}*\theta_0\|_{L^p}^p
\\
&&{}+
\int_0^t\frac{1}{2}p(p-1) \int
_{\mathbb{T}^2}\bigl |\theta(s)\bigr |^{p-2}\biggl(\sum
_j\bigl |k_{\delta_n}*G\bigl(\theta (s)\bigr)
(f_j)\bigr |^2\biggr)\,d\xi\,ds
\\
&&{}+p\int_0^t \int_{\mathbb{T}^2}\bigl |
\theta(s)\bigr |^{p-2}\theta(s)k_{\delta_n}*G\bigl(\theta (s)\bigr)\,d\xi
\,dW(s)
\\
&\leq&\|k_{\delta_n}*\theta_0\|_{L^p}^p
\\
&&{}+
\int_0^t\biggl(\varepsilon\int
_{\mathbb{T}^2}\bigl |\theta(s)\bigr |^p\,d\xi
\\
&&\phantom{+
\int_0^t\biggl(}{}+C(\varepsilon)\int
\biggl(\sum_j\bigl |k_{\delta_n}*G\bigl(\theta(s)
\bigr) (f_j)\bigr |^2\biggr)^{p/2} \,d\xi\biggr)\,ds
\\
&&{}+p\int_0^t\int_{\mathbb{T}^2}\bigl |
\theta(s)\bigr |^{p-2}\theta(s) k_{\delta
_n}*G\bigl(\theta(s)\bigr)\,d\xi
\,dW(s),
\end{eqnarray*}
where in the first inequality we used $\operatorname{div} u=0$ and
$\int|\theta
|^{p-2}\theta\Lambda^{2\alpha}\theta\geq0$ (cf. \cite{Re95},
Lemma~3.2) as
well as Young's inequality in the second inequality.
Then by the Burkholder--Davis--Gundy inequality and Minkowski's
inequality, we obtain\looseness=-1
%
%e3.11 #&#
%
\begin{eqnarray}
\label{eq3.11}
&&E\sup_{s\in[0,t]}\bigl \|\theta(s)\bigr \|_{L^p}^p\nonumber
\\[-1pt]
&&\qquad
\leq  E\| \theta_0\|_{L^p}^p\nonumber
\\[-1pt]
&&\qquad\quad  {}+E\int
_0^t\biggl(\varepsilon\int_{\mathbb{T}^2}\bigl |
\theta (s)\bigr |^p\,d\xi+C\int\biggl(\sum_j\bigl |k_{\delta_n}*G
\bigl(\theta (s)\bigr) (f_j)\bigr |^2\biggr)^{p/2}\,d
\xi \biggr)\,ds\hspace*{-60pt}
\nonumber
\\[-1pt]
&&\qquad\quad  {} +pE\biggl(\int_0^t\biggl(\int
_{\mathbb{T}^2}\bigl |\theta(s)\bigr |^{p-1}\biggl(\sum
_j\bigl |k_{\delta
_n}*G\bigl(\theta(s)\bigr)
(f_j)\bigr |^2\biggr)^{1/2}\,d\xi
\biggr)^2\,ds\biggr)^{1/2}
\nonumber
\\[-1pt]
&&\qquad  \leq  E\|\theta_0\|_{L^p}^p\nonumber
\\[-1pt]
&&\qquad\quad  {}+E\int
_0^t\biggl(\varepsilon\int_{\mathbb
{T}^2}\bigl |
\theta(s)\bigr |^p\,d\xi+C\int\biggl(\sum_j\bigl |k_{\delta_n}*G
\bigl(\theta (s)\bigr) (f_j)\bigr |^2\biggr)^{p/2}\,d
\xi\biggr)\,ds
\nonumber
\\[-1pt]
&&\qquad\quad  {} +pE\sup_{s\in[0,t]}\bigl \|\theta(s)\bigr \|_{L^p}^{p-1}
\\[-1pt]
&&\qquad\quad\hspace*{54pt}{}\times\biggl(\int_0^t\biggl(\int_{\mathbb
{T}^2}
\biggl(\sum_j\bigl |k_{\delta_n}*G\bigl(\theta(s)
\bigr) (f_j)\bigr |^2\biggr)^{p/2}\,d\xi
\biggr)^{2/p}\,ds\biggr)^{1/2}\nonumber
\\[-1pt]
&&\qquad  \leq  E\|\theta_0\|_{L^p}^p\nonumber
\\[-1pt]
 &&\qquad\quad {}+E\int
_0^t\biggl(\varepsilon\int_{\mathbb
{T}^2}\bigl |
\theta(s)\bigr |^p\,d\xi+C\int\biggl(\sum_j\bigl |G
\bigl(\theta (s)\bigr) (f_j)\bigr |^2\biggr)^{p/2}\,d
\xi\biggr)\,ds
\nonumber
\\[-1pt]
&&\qquad\quad  {} +C(T)E\sup_{s\in[0,t]}\bigl \|\theta(s)\bigr \|_{L^p}^{p-1}
\biggl(\int_0^t\biggl(\int_{\mathbb{T}^2}
\biggl(\sum_j\bigl |G\bigl(\theta(s)\bigr)
(f_j)\bigr |^2\biggr)^{p/2}\,d\xi\biggr)\, ds
\biggr)^{1/p}\hspace*{-40pt}
\nonumber
\\[-1pt]
&&\qquad  \leq  E\|\theta_0\|_{L^p}^p+\varepsilon E
\sup_{s\in[0,t]}\bigl \|\theta (s)\bigr \|_{L^p}^p+C_1E
\int_0^t\bigl \|\theta(s)\bigr \|_{L^p}^p
\,ds+C_2
\nonumber
\\[-1pt]
&&\qquad  \leq  E\|\theta _0\|_{L^p}^p+\varepsilon E
\sup_{s\in[0,t]}\bigl \|\theta(s)\bigr \| _{L^p}^p+C_1
\int_0^tE\sup_{s\in[0,\sigma]}\bigl \|
\theta(s)\bigr \|_{L^p}^p\, d\sigma +C_2.
\nonumber
\end{eqnarray}\looseness=0%
Here, in the fourth inequality, we used (Gp.1) and Young's inequality.
By Gronwall's lemma, the assertion follows.
\end{pf*}

%$\hfill\Box$

%s4 #&#
\section{Existence and uniqueness of probabilistically (strong)
solutions in the subcritical case}\label{sec4}

In this section, we assume $\alpha>1/2$ and prove pathwise uniqueness
for equation (\ref{eq3.1}), and hence by the Yamada--Watanabe theorem the
existence of a unique (probabilistically) strong solution to (\ref{eq3.1})
in the subcritical case.
Let us first give the definition of a (probabilistically) strong
solution to (\ref{eq3.1}).

%de4.1 #&#
%
\begin{Definition}\label{def4.1}
We say that there exists a (probabilistically)
strong solution to (\ref{eq3.1}) over the time interval $[0,T]$
if for every
probability space $(\Omega,\mathcal{F},\{\mathcal{F}_t\}_{t\in
[0,T]},P)$ with an $\mathcal{F}_t$-Wiener process $W$, there exists an
$\mathcal{F}_t$-adapted process $\theta\dvtx[0,T]\times\Omega
\rightarrow
H$ such that
for $P$-a.e. $\omega\in\Omega$
\[
\theta(\cdot,\omega)\in L^\infty(0,T;H)\cap L^2
\bigl(0,T;H^\alpha\bigr)\cap C\bigl([0,T];H^{-\beta}\bigr)
\]
and $P$-a.e.
%e4.1 #&#
%
\begin{eqnarray}
\label{eq4.1}
&&\bigl\langle\theta(t),\varphi\bigr\rangle+\int
_0^t\bigl\langle A_\alpha
^{1/2}\theta (s),A_\alpha^{1/2}\varphi\bigr\rangle\,
ds-\int_0^t \bigl\langle u(s)\cdot \nabla
\varphi,\theta(s)\bigr\rangle\,ds
\nonumber
\\[-8pt]
\\[-8pt]
&&\qquad=\langle\theta_0,\varphi\rangle +
\biggl\langle\int_0^tG\bigl(\theta(s)\bigr)
\,dW(s),\varphi\biggr\rangle\nonumber
\end{eqnarray}
for all $t\in[0,T]$ and all $\varphi\in C^1(\mathbb{T}^2)$ (assuming
also that all integrals in the equation are defined).
\end{Definition}

%th4.2 #&#
%
\begin{Theorem}\label{the4.2}
Assume $\alpha>\frac{1}{2}$. If $G$ satisfies the
following condition:
\renewcommand{\theequation}{GL.1}
%e4.2 #&#
%
\begin{eqnarray}
\label{eqGL.1} \bigl \| \Lambda^{-1/2}\bigl(G(u)-G(v)\bigr)
\bigr \|^2_{L_2(U,H)}&\leq&\beta\bigl |\Lambda ^{-1/2}(u-v)\bigr |^2
\nonumber
\\[-8pt]
\\[-8pt]
&&{}+
\beta_1\bigl |\Lambda^{\alpha-{1}/{2}}(u-v)\bigr |^2\nonumber
\end{eqnarray}
for all $u,v\in H^\alpha$, for some $\beta\in\mathbb{R}$ independent
of $u,v$, and $\beta_1<2\kappa$,
then (\ref{eq3.1}) admits at most one probabilistically strong
solution in the
sense of Definition~\ref{def4.1} such that
\[
\sup_{t\in[0,T]}\bigl \|\theta(t)\bigr \|_{L^p}<\infty,\quad P
\mbox{-a.s.}
\]
for some $p\in((\alpha-\frac{1}{2})^{-1},\infty)$, and
\[
E\sup_{t\in[0,T]}\bigl |\Lambda^{-1/2}\theta(t)\bigr |^2<
\infty.
\]
\end{Theorem}

\begin{Remark*}
The examples in Remark~\ref{rem3.4} with $g$ being a Lipschitz
function on $\mathbb{R}$ satisfy (\ref{eqGL.1}) since
\begin{eqnarray*}
\bigl \| \Lambda^{-1/2}\bigl(G(u)-G(v)\bigr)\bigr \|^2_{L_2(K,H)}&=&
\sum_k\bigl |\Lambda ^{-1/2}
\bigl(b_k\bigl(g(u)-g(v)\bigr)\bigr)\bigr |^2
\\
&\leq& \int_{\mathbb{T}^2}\sum_kb_k^2
\bigl(g(u)-g(v)\bigr)^2\,d\xi
\\
&\leq& C|u-v|^2
\\
&\leq& C\bigl |\Lambda ^{-1/2}(u-v)\bigr |^2+
\varepsilon\bigl |\Lambda^{\alpha-
{1}/{2}}(u-v)\bigr |^2.
\end{eqnarray*}
\end{Remark*}

\begin{pf*}{Proof of Theorem~\ref{the4.2}}
Let $\theta_1$, $\theta_2$ be two
solutions of (\ref{eq3.1}), and let $\{e_k\}_{k\in\mathbb{N}}$ be the
eigenbasis of $A_\alpha$ from above. Then their difference $\theta
=\theta_1-\theta_2$ satisfies for $\psi\in C^1(\mathbb{T}^2)$\vspace*{1pt}
\setcounter{equation}{1}
%e4.2 #&#
%
\begin{eqnarray}
\label{eq4.2}
&& \bigl\langle\psi,\theta(t)\bigr\rangle-\int_0^t
\langle u\cdot\nabla\psi ,\theta _1\rangle\,ds -\int
_0^t\langle u_2\cdot\nabla\psi,
\theta\rangle \,ds+\kappa\int_0^t\bigl\langle
\theta,\Lambda^{2\alpha}\psi\bigr\rangle\, ds
\nonumber
\\[-8pt]
\\[-8pt]
&&\qquad=\int_0^t
\bigl\langle\psi,\bigl(G(\theta_1)-G(\theta_2)\bigr)\,dW
\bigr\rangle.\nonumber
\end{eqnarray}
Here, $u_1,u_2,u$ satisfy (\ref{eq1.3}) with $\theta$ replaced by
$\theta
_1,\theta_2,\theta$, respectively.
Now set $\phi_k=\langle e_k,\theta(t)\rangle$, $\varphi_k=\langle
\Lambda
^{-1}e_k,\theta(t)\rangle$. It\^{o}'s formula and (\ref{eq4.2}) yield\vspace*{1pt}
%
%e4.3 #&#
%
\begin{eqnarray}
\label{eq4.3} \phi_k\varphi_k&=&\int
_0^t\phi_k\,d
\varphi_k+\int_0^t
\varphi_k\, d\phi _k+\langle\varphi_k,
\phi_k\rangle(t)
\nonumber
\\
&=&2\int_0^t\langle u\cdot\nabla
e_k,\theta_1\rangle \bigl\langle\Lambda^{-1}
\theta,e_k\bigr\rangle+\langle u_2\cdot\nabla
e_k,\theta\rangle \bigl\langle\Lambda^{-1}\theta
,e_k\bigr\rangle
\nonumber
\\
&&\phantom{2\int_0^t}{} -\kappa\bigl\langle\Lambda^{2\alpha}e_k,
\theta\bigr\rangle \bigl\langle\Lambda ^{-1}\theta,e_k
\bigr\rangle\,ds
\\
&&{}+2\int_0^t\bigl\langle
\Lambda^{-1}\theta ,e_k\bigr\rangle \bigl\langle
e_k,\bigl(G(\theta_1)-G(\theta_2)\bigr)
\,dW(s)\bigr\rangle
\nonumber
\\
&&{}+\int_0^t\bigl\langle\bigl(G(
\theta_1)-G(\theta_2)\bigr)^*e_k,\bigl(G(
\theta_1)-G(\theta _2)\bigr)^*\Lambda^{-1}e_k
\bigr\rangle_U \,ds.
\nonumber
\end{eqnarray}
Here, $\langle\varphi_k,\phi_k\rangle(t)$ denotes the covariation
process of $\varphi_k,\phi_k$.
The dominated convergence theorem implies\vspace*{1pt}
\begin{eqnarray*}
\sum_{k\leq N}\int_0^t
\langle u\cdot\nabla e_k,\theta_1\rangle \bigl\langle
\Lambda^{-1}\theta,e_k\bigr\rangle\,ds&\rightarrow&\int
_0^t { }_{H^{-1}} \bigl\langle u\cdot
\nabla\theta_1, \Lambda^{-1}\theta\bigr\rangle_{H^1}
\,ds,\qquad N\rightarrow\infty,
\\
\sum_{k\leq N}\int_0^t
\langle u_2\cdot\nabla e_k,\theta\rangle \bigl\langle
\Lambda^{-1}\theta,e_k\bigr\rangle\,ds&\rightarrow&\int
_0^t { }_{H^{-1}} \bigl\langle
u_2\cdot\nabla\theta, \Lambda^{-1}\theta\bigr
\rangle_{H^1} \,ds,\qquad N\rightarrow\infty
\end{eqnarray*}
and\vspace*{1pt}
\[
\sum_{k\leq N}\int_0^t
\bigl\langle\Lambda^{2\alpha}e_k,\theta\bigr\rangle \bigl
\langle \Lambda^{-1}\theta,e_k\bigr\rangle\,ds\rightarrow
\int_0^t \bigl\langle \theta,
\Lambda^{2\alpha-1}\theta\bigr\rangle\,ds,\qquad N\rightarrow\infty.
\]
Furthermore, since\vspace*{1pt}
\begin{eqnarray*}
&&\int_0^t\bigl |\Lambda^{-1/2}
\theta\bigr |^2\bigl \|\Lambda^{-1/2}\bigl( G(\theta _1)-G(
\theta _2)\bigr)\bigr \|^2_{L_2(U,H)}\,ds
\\
&&\qquad\leq C\sup
_{s\leq t}\bigl |\theta(s)\bigr |^2\int_0^t
\bigl \| \Lambda ^{-1/2}\bigl(G(\theta_1)-G(\theta_2)
\bigr)\bigr \|^2_{L_2(U,H)}\,ds<\infty,
\end{eqnarray*}
we obtain
\begin{eqnarray*}
&&\sum_{k\leq N}\int_0^t
\bigl\langle\Lambda^{-1}\theta,e_k\bigr\rangle \bigl
\langle e_k,\bigl(G(\theta_1)-G(\theta_2)
\bigr)\,dW(s)\bigr\rangle
\\
&&\qquad\rightarrow M_t:=\int_0^t
\bigl\langle\Lambda^{-1/2}\theta,\Lambda^{-1/2}\bigl(G(\theta
_1)-G(\theta _2)\bigr)\,dW(s)\bigr\rangle,\qquad N
\rightarrow\infty,
\end{eqnarray*}
in probability.
Finally, the following inequality holds:
\begin{eqnarray*}
&&\sum_{k\leq N}\int_0^t
\bigl\langle\bigl(G(\theta_1)-G(\theta _2)
\bigr)^*e_k,\bigl(G(\theta _1)-G(\theta_2)
\bigr)^*\Lambda^{-1}e_k\bigr\rangle_U \,ds
\\
&&\qquad\leq\int_0^t\bigl \| \Lambda ^{-1/2}\bigl(G(
\theta_1)-G(\theta_2)\bigr)\bigr \|^2_{L_2(U,H)}
\,ds.
\end{eqnarray*}
Thus, summing up over $k\leq N$ in (\ref{eq4.3}) and letting
$N\rightarrow\infty
$, we obtain
\begin{eqnarray*}
&&\bigl |\Lambda^{-1/2}\theta\bigr |^2+2\kappa\int_0^t
\bigl |\Lambda^{\alpha-
{1}/{2}}\theta\bigr |^2 \,ds
\\
&&\qquad\leq 2M(t)+2\int_0^t {
}_{H_{-1}} \bigl\langle u\cdot \nabla\theta_1,
\Lambda^{-1}\theta\bigr\rangle_{H_1}+{ }_{H_{-1}} \bigl
\langle u_2\cdot\nabla\theta, \Lambda^{-1}\theta\bigr
\rangle_{H_1}\,ds
\\
&&\phantom{\qquad\leq} {} +\int_0^t\bigl \|
\Lambda^{-1/2}\bigl(G(\theta_1)-G(\theta_2)\bigr)
\bigr \|^2_{L_2(U,H)}\,ds.
\end{eqnarray*}
By \cite{Re95}, we have
\[
{ }_{H^{-1}} \bigl\langle u\cdot\nabla\theta_1,
\Lambda^{-1}\theta \bigr\rangle_{H^1}=0
\]
and
\begin{eqnarray*}
&&\bigl |{}_{H_{-1}} \bigl\langle u_2\cdot\nabla\theta,
\Lambda^{-1}\theta \bigr\rangle_{H_1}\bigr |
\\
&&\qquad\leq  \|u_2
\|_{L^p}\|\theta\|_{L^{p_1}}\bigl \|\nabla \Lambda ^{-1}\theta
\bigr \|_{L^{p_1}}\leq C\|u_2\|_{L^p}\|\theta
\|_{H^{1/p}}\bigl \| \nabla \Lambda^{-1}\theta\bigr \|_{H^{1/p}}
\\
&&\qquad\leq C\|\theta_2\|_{L^p}\bigl \|\Lambda ^{-1}\theta
\bigr \|_{H^{1+{1}/{p}}}^2\leq C\|\theta_2\|_{L^p}\bigl \|
\Lambda ^{-1}\theta\bigr \|^{2/r}_{H^{1/2}}\bigl \|
\Lambda^{-1}\theta\bigr \|^{2(1-
{1}/{r})}_{H^{{1}/{2}+\alpha}}
\\
&&\qquad\leq  \varepsilon\bigl |\Lambda^{\alpha
-{1}/{2}}\theta\bigr |^2+ C\|
\theta_2\|_{L^p}^r\bigl |\Lambda^{-1/2}\theta
\bigr |^2,
\end{eqnarray*}
where $\frac{1}{p}+\frac{2}{{p_1}}=1$ for $p\in((\alpha-\frac
{1}{2})^{-1},+\infty), r=\frac{\alpha}{\alpha-{1}/{2}-{1}/{p}}$.
Here we use $\operatorname{div} u_2=0$ in the first inequality, that
$H^{1/p}\hookrightarrow L^{p_1}$ continuously in the second inequality,
the interpolation inequality (\ref{eq2.3}) in the fourth inequality
and Young's
inequality in the last equality.

Now by (\ref{eqGL.1}) we have
\[
\bigl |\Lambda^{-1/2}\theta\bigr |^2\leq2M(t)+\int_0^tC
\|\theta_2\| _{L^p}^r\bigl |\Lambda^{-1/2}
\theta\bigr |^2\,ds+\beta\int_0^t\bigl |\Lambda
^{-1/2}(\theta _1-\theta_2) \bigr |^2\,ds.
\]

Let
\[
\tau_n^1:=\inf\bigl\{t>0, \bigl \|\theta_2(t)
\bigr \|_{L^p}> n\bigr\}.
\]
Then by the weak continuity of $\theta_2$, $\tau_n^1$ are
stopping times with respect to $\mathcal{F}_{t+}, ( \mathcal
{F}_{t+}:=\bigcap_{s>t}\mathcal{F}_s )$ and $\|\theta_2(t\wedge\tau
_n^1)\|
_{L^p}\leq n$ for large $n$. Furthermore, let $\tau_n^2$ be a
localizing sequence of stopping times for $M$ and $\tau_n:=\tau
_n^1\wedge\tau_n^2$. Then, since $M(t\wedge\tau_n)$ is a martingale
with respect to $\mathcal{F}_{t+}$, we get
\begin{eqnarray*}
E\bigl |\Lambda^{-1/2}\theta(t\wedge\tau_n)\bigr |^2&
\leq&Cn^rE\int_0^{t\wedge\tau_n}\bigl |
\Lambda^{-1/2}\theta\bigr |^2\,ds+\beta E\int_0^{t\wedge
\tau_n}
\bigl |\Lambda^{-1/2}\theta \bigr |^2\,ds
\\
&=& C(n)\int_0^tE\bigl |\Lambda^{-1/2}
\theta(s\wedge\tau_n)\bigr |^2\,ds
\\
&&{}+\beta \int
_0^t E\bigl |\Lambda^{-1/2}\theta(s\wedge
\tau_n) \bigr |^2\,ds.
\end{eqnarray*}
By Gronwall's inequality, we get $|\Lambda^{-1/2}\theta(t\wedge\tau
_n)|^2=0$ $P$-a.s., and recalling that $\tau_n\rightarrow T$
$P$-a.s. as $n\rightarrow\infty$, we obtain that $\theta(t)=0$
$P$-a.s. for $t\leq T$. By the weak continuity of $\theta$, we obtain
the zero set does not depend on $t$, thus completing the proof.
\end{pf*}
%
%$\hfill\Box$

\begin{Remark*}
From the proof of Theorem~\ref{the4.2}, we immediately obtain that
if there exists a probabilistically strong solution $\theta$ in the
sense of Definition~\ref{def3.1} satisfying
\[
\sup_{t\in[0,T]}\bigl \|\theta(t)\bigr \|_{L^p}<\infty,\qquad P
\mbox{-a.s.}
\]
for some $p\in((\alpha-\frac{1}{2})^{-1},+\infty)$ and $G$ satisfies
(\ref{eqGL.1}), then for any other solution $\tilde{\theta}$ such that
\[
E\sup_{t\in[0,T]}\bigl |\Lambda^{-1/2}\tilde{
\theta}(t)\bigr |^2<\infty,
\]
it follows that $\tilde{\theta}=\theta$, which implies that
\[
\sup_{t\in[0,T]}\bigl \|\tilde{\theta}(t)\bigr \|_{L^p}<\infty.
\]
\end{Remark*}

%th4.3 #&#
%
\begin{Theorem}\label{the4.3}
Assume $\alpha>\frac{1}{2}$ and that $G$ satisfies
Hypothesis~\ref{hypG.1},\break (\ref{eqGL.1}) and \textup{(Gp.1)} for some $p\in
((\alpha-\frac
{1}{2})^{-1},+\infty)$.
Then for each initial condition $\theta_0\in L^p$, there exists a
pathwise unique probabilistically strong solution $\theta$ of equation
(\ref{eq3.1}) over $[0,T]$ with initial condition $\theta(0)=\theta
_0$ such
that
\[
E\sup_{t\in[0,T]}\bigl |\Lambda^{-1/2}\theta(t)\bigr |^2<
\infty.
\]
Moreover, the solution satisfies
\[
E\sup_{t\in[0,T]}\bigl \|\theta(t)\bigr \|^p_{L^p}+E\int
_0^T\bigl |\Lambda^\alpha \theta
(t)\bigr |^2\,dt<\infty.
\]
\end{Theorem}

\begin{pf}
By Theorem~\ref{the4.2}, Theorem~\ref{the3.3} and the
Yamada--Watanabe theorem
(cf. \cite{RSZ08} or \cite{Ku07,PR07}), we get that for each initial condition $\theta
_0\in L^p$,
there exists a pathwise unique probabilistically strong solution
$\theta$ of equation (\ref{eq3.1}) over $[0,T]$ with initial
condition $\theta
(0)=\theta_0$ such that
\[
\sup_{t\in[0,T]}\bigl \|\theta(t)\bigr \|_{L^p}<\infty,\qquad P
\mbox{-a.s.},
\]
and
\[
E\sup_{t\in[0,T]}\bigl |\Lambda^{-1/2}\theta(t)\bigr |^2<
\infty.
\]
By the remark before Theorem~\ref{the4.3}, the first result follows.
By Theorem~\ref{the3.3} and (\ref{eq3.6}), the last part of the
assertion follows.
\end{pf}
%
%$\hfill\Box$

%th4.4 #&#
%
\begin{Theorem}[(Markov property)]\label{the4.4}
Assume $\alpha>\frac{1}{2}$ and that
$G$ satisfies Hypothesis~\ref{hypG.1}, (\ref{eqGL.1}) and
\textup{(Gp.1)} for some $p\in((\alpha-\frac
{1}{2})^{-1},+\infty)$.
If $\theta_0\in L^p$, then for every bounded, $\mathcal
{B}(H)$-measurable $F\dvtx H\rightarrow\mathbb{R}$, and all $s,t\in[0,T]$,
$s\leq t$
\[
E\bigl(F\bigl(\theta(t)\bigr)|\mathcal{F}_s\bigr) (\omega)=E\bigl(F
\bigl(\theta\bigl(t,s,\theta (s) (\omega )\bigr)\bigr)\bigr) \qquad\mbox{for } P
\mbox{-a.s. } \omega\in\Omega.
\]
Here, $\theta(t,s,\theta(s)(\omega))$ denotes the solution to (\ref{eq3.1})
starting from $\theta(s)$ at time $s$ satisfying
\[
E\sup_{t\in[s,T]}\bigl |\Lambda^{-1/2}\theta(t)\bigr |^2<
\infty.
\]
\end{Theorem}

\begin{pf}
By Theorem~\ref{the4.3}, we have $\theta(t)=\theta(t,s,\theta(s))$
$P$-a.s. Then by the Yamada--Watanabe theorem in \cite{RSZ08}, we have $P$-a.s.
\begin{eqnarray*}
E\bigl(F\bigl(\theta(t)\bigr)|\mathcal{F}_s\bigr) (\omega)&=&E
\bigl(F\bigl(\theta\bigl(t,s,\theta (s)\bigr)\bigr)|\mathcal{F}_s
\bigr) (\omega)
\\
&=&E\bigl(F\bigl(\mathbf{H}\bigl(\theta(s),W(\cdot+s)-W(s)\bigr)\bigr)|
\mathcal{F}_s\bigr) (\omega)
\\
&=&E\bigl(F\bigl(\mathbf{H}\bigl(\theta(s) (\omega),W(\cdot+s)-W(s)\bigr)\bigr)
\bigr)
\\
&=&E\bigl(F\bigl(\theta\bigl(t,s,\theta(s) (\omega)\bigr)\bigr)\bigr),
\end{eqnarray*}
where $\mathbf{H}$ is the functional obtained by the Yamada--Watanabe
theorem such that $\mathbf{H}(\theta(0),W)$ is a strong solution to
(\ref{eq3.1}).
\end{pf}
%
%$\hfill\Box$

We set for $\mathcal{B}(H)$-measurable $F\dvtx H\rightarrow\mathbb
{R}$, and
$t\in[0,T], x\in L^p$
\[
P_tF(x):=EF\bigl(\theta(t,x)\bigr).
\]
Here, and in the following, we use $\theta(t,x)$ to denote a solution
with initial value $x$.
Then by Theorem~\ref{the4.4}, we have for $F\dvtx H\rightarrow\mathbb
{R}$, bounded
and $\mathcal{B}(H)$-measurable, $s,t\geq0$,
\[
P_s(P_tF) (x)=P_{s+t} F(x),\qquad x\in
L^p, p\in\bigl(\bigl(\alpha-\tfrac
{1}{2}\bigr)^{-1},+
\infty\bigr).
\]

%s5 #&#
\section{Ergodicity in the subcritical case}\label{sec5}

Now fix $\alpha>\frac{1}{2}$ and we assume $U=H$, $W(t)$ is a
cylindrical Wiener process in $H$ defined on a filtered probability
space $(\Omega,\mathcal{F},\{\mathcal{F}_t\}_{t\geq0},P)$. We make the
following assumptions on $G$.

\renewcommand{\theHypothesis}{E.1}
%hy2 #&#
%
\begin{Hypothesis}\label{hypE.1}
$G$ does not depend on $\theta$ and there
exists $\sigma>0$ such that $G\in L_2(H;H^{2-\alpha+\sigma})$ that is,
\[
\mathcal{E}_0:=\operatorname{Tr}\bigl(\Lambda^{4-2\alpha+2\sigma} GG^*
\bigr)<\infty.
\]
\end{Hypothesis}

\renewcommand{\theHypothesis}{E.2}
%hy3 #&#
%
\begin{Hypothesis}\label{hypE.2}
There exist $N\in\mathbb{N}$ and $g\in L(H)$
such that $Gg=P_N$.
\end{Hypothesis}

For $\varepsilon_0>0$ and any $\overline{W}\in C(\mathbb
{R}^+,H^{-1-\varepsilon_0})$, we define
\[
z(\overline{W}) (t):=\sum_{i,j=1}^\infty
\biggl(g_{ij}\beta_i(t)-\lambda _j\int
_0^te^{-\lambda_j(t-s)}g_{ij}
\beta_i(s) \,ds\biggr)e_j,
\]
if the convergence of the sum is uniformly with respect to $t$ in every
bounded time interval, otherwise set $z(\overline{W}):=+\infty$. Here,
$\beta_i(t):={}_{H^{1+\varepsilon_0}} \langle e_i, \overline
{W}(t)\rangle_{H^{-1-\varepsilon_0}}$, $g_{ij}=\langle
Ge_i,e_j\rangle
$. Under Hypothesis~\ref{hypE.1}, there exists $\Omega'\subset\Omega
$ such that $P(\Omega
')=1$ and for $\omega\in\Omega'$, $z(W(\omega))\in C([0,\infty),
H^{2+\varepsilon})$ for some $0<\varepsilon<\sigma$, and on $(\Omega
,\mathcal{F},\mathcal{F}_t,P)$,
$z(W)$ is the mild solution of the equation:
$dz+A_\alpha z=G\,dW$ with initial condition $z(0)=0$.

Now for $v_0\in H^1, \overline{W}\in C(\mathbb{R}^+,H^{-1-\varepsilon
_0}) $ we define
\[
v(t,\overline{W},v_0):=\cases{ %
v
\bigl(t,v_0,z(\overline{W})\bigr),& $\mbox{if } z(\overline {W})\in C
\bigl(\mathbb{R}^+,H^{m}\bigr) \mbox{ for } m<2+\sigma$,
\cr
0,& $
\mbox{otherwise}$, %\end{array}
}
\]
where $v(t,v_0,z(\overline{W}))$ is the solution to (\ref{eqA.1}) we obtained
in Theorem~\ref{theA.1}.
Then by Theorem~\ref{theA.4} in Appendix~\hyperref[appA]{A}, $v$ is a
measurable mapping from
$\mathbb{R}^+\times C(\mathbb{R}^+,H^{-1-\varepsilon_0})\times H^1$
into $H^1$, $(t,\overline{W},\theta_0)\mapsto v(t,\overline
{W},\theta_0)$.
We can now define
\[
\theta(t,\overline{W},\theta_0):=v(t,\overline{W},\theta
_0)+z(t,\overline{W}),
\]
which is a measurable map from $\mathbb{R}^+\times C(\mathbb
{R}^+,H^{-1-\varepsilon_0})\times H^1$ into $H^1$. Then for the
cylindrical Wiener process $W$, $\theta(t,W,\theta_0)$ is a solution to
(\ref{eq3.1}), whose laws $P_{\theta_0},\theta_0\in H^1$ form a
Markov process
on $ H^1$, since $H^1$ is an invariant space for (\ref{eq3.1}) under
assumption Hypothesis~\ref{hypE.1}.
Let $(P_t)_{t\geq0}$ be the associated transition semigroup on
$\mathcal{B}_b( H^1)$.
Now we want to study the long time behavior of the semigroup $P_t$.

%re5.1 #&#
%
\begin{remark}\label{rem5.1}
(i) Hypothesis~\ref{hypE.1} obviously implies Hypothesis~\ref
{hypG.1},\break (Gp.1) for
all $p\in((\alpha-\frac{1}{2})^{-1},\infty)$ and (\ref{eqGL.1}).
For $x:=\theta
_0\in L^p$,
let $P_x$ denote the law of the corresponding solution $\theta$ to
(\ref{eq3.1}). Then by Theorems~\ref{the4.3} and~\ref{the4.4}, the
measures $P_x$, $x\in L^p$ form
a Markov process.

(ii) The existence of a map $g$ such that $Gg=P_N$ is equivalent to the
following property:
\[
P_NH\subset\operatorname{Im}(G).
\]

(iii) Hypothesis~\ref{hypE.1}
is to make sure that the associated O--U process
has a version $z\in C([0,\infty);H^{1,\infty}(\mathbb{T}^2))$ (see,
e.g., \cite{DZ92}, the proof of Theorem~5.16, and use Sobolev
embedding). If
we consider the stochastic integral taking values in a Banach space
[e.g., $L^p(\mathbb{T}^2)$, $p>1$] and use the theory developed in
\cite{B97}, we can change Hypothesis~\ref{hypE.1}
to the following condition:
$G\in L_2(H;H^{1-\alpha+{\varepsilon_1}/{2}})$ and for some
$\varepsilon_1,q$ satisfying $\varepsilon_1 q>2$
\[
\biggl\|\biggl[\sum_k \bigl(\Lambda^{1-\alpha+\varepsilon_1}Ge_k
\bigr)^2\biggr]^{1/2}\biggr\|_{L^q}+\biggl\| \biggl[\sum
_k (Ge_k)^2
\biggr]^{1/2}\biggr\|_{L^{({\alpha+1})/({\alpha-
{1}/{2}})}}<\infty.
\]
By this and similar arguments as in \cite{B97}, we obtain for
$\varepsilon
<\varepsilon_1$ and $\varepsilon q>2$ that the O--U process has a
version $z\in C([0,\infty);H^{1+\varepsilon,q})\subset C([0,\infty
);H^{1,\infty}(\mathbb{T}^2))$, but in this paper we stay in the
Hilbert space framework for simplicity.

(iv) For more general noise, we do not know how to obtain
Proposition~\ref{pro5.2.1} since we cannot control $E\exp\|\theta\|
_{L^p}^p$. Therefore, we
restrict ourselves to additive noise.
\end{remark}

%s5.1 #&#
\subsection{Preliminaries and some useful estimates}\label{sec5.1}

First, we want to collect some useful and fundamental results about
coupling from \cite{L92} and \cite{M02} which we will use later. Let
$(\Lambda
_1,\Lambda_2)$ be two probability measures on a Polish space $E$. Let
$(Z_1,Z_2)$ be a couple of random variables $(\Omega,\mathcal
{F})\rightarrow E\times E$. We say that $(Z_1,Z_2)$ is a coupling of
$(\Lambda_1,\Lambda_2)$ if $\Lambda_i=\mathcal{D}(Z_i)$ for $i=1,2$,
where we use $\mathcal{D}(Z_i)$ to denote the distribution of $Z_i$.

%le5.2 #&#
%
\begin{Lemma}\label{lem5.1.1}
Let $(\Lambda_1,\Lambda_2)$ be two probability
measures on a Polish space $(E,\mathcal{B}(E))$. Then
\[
\|\Lambda_1-\Lambda_2\|_{\mathrm{var}}=\min
P(Z_1\neq Z_2),
\]
where the minimum is taken over all couplings $(Z_1,Z_2)$ of $(\Lambda
_1,\Lambda_2)$. There exists a coupling for which the minimum value is
attained and it is called a maximal coupling. Moreover, the maximal
coupling has the following property:
\[
P(Z_1=Z_2,Z_1\in\Gamma)=(
\Lambda_1\wedge\Lambda_2) (\Gamma),\qquad \Gamma\in
\mathcal{B}(E).
\]
\end{Lemma}

%le5.3 #&#
%
\begin{Lemma}[(cf. \cite{M02}, Lemma C.1)]\label{lem5.1.2}
Let $\Lambda_1$ and $\Lambda_2$
be two equivalent probability measures on $E$. Then for any $p>1$ and
any measurable subset $A\subset E$
\[
I_p(A):=\int_A\biggl(\frac{d\Lambda_1}{d\Lambda_2}
\biggr)^p\,d\Lambda_1<\infty
\]
implies
\[
(\Lambda_1\wedge\Lambda_2) (A)\geq\biggl(1-
\frac{1}{p}\biggr) \biggl(\frac{\Lambda
_1(A)^p}{pI_p(A)}\biggr)^{{1}/({p-1})}.
\]
\end{Lemma}

%pr5.4 #&#
%
\begin{Proposition}[(cf. \cite{O08}, Proposition~1.4)]\label{pro5.1.3}
Let $E$ and $F$ be
two Polish spaces, $f_0\dvtx E\rightarrow F$ be a measurable map and
$(\Lambda_1,\Lambda_2)$ be two probability measures on $E$. Set
$\lambda
_i=f_0^*\Lambda_i$, $i=1,2$. Then there exists a coupling $(V_1,V_2)$
of $(\Lambda_1,\Lambda_2)$ such that $(f_0(V_1),f_0(V_2))$ is a maximal
coupling of $(\lambda_1,\lambda_2)$.
\end{Proposition}

Now we give some useful estimates which will be used in the next two
subsections. Let $\theta_n$ denote the approximation in the proof of
Theorem~\ref{the3.3}. As will be seen below, we shall need uniform
$L^p$-estimates, and a crucial ingredient to prove them is the
following improved version of the ``positivity lemma,'' that is, Lemma~3.2 in \cite{Re95}.

%le5.5 #&#
%
\begin{Lemma}[(Improved positivity lemma)]\label{lem5.1.4}
For $\alpha\in(0,1)$, and
$\theta\in L^p$ with $\Lambda^{2\alpha}\theta\in L^p$, for some
$2<p<\infty$,
\[
\int|\theta|^{p-2}\theta\biggl(\kappa\Lambda^{2\alpha}-
\frac{2\lambda
_1}{p}\biggr)\theta\geq0.
\]
\end{Lemma}

\begin{pf}
Denote the semigroup with respect to $-\kappa\Lambda^{2\alpha
}+\frac{2\lambda_1}{p}$ and $-\kappa\Lambda^{2\alpha}$ in $L^2$ by
$P^0_t$ and $P_t^1$, respectively.
Then we have $P^0_tf=e^{2t\lambda_1/p}P_t^1f$. Since
\[
\bigl \|P_t^1f\bigr \|_{L^2}\leq e^{-\lambda_1t}\|f
\|_{L^2}
\]
and
\[
\bigl \|P_t^1f\bigr \|_{L^\infty}\leq\|f\|_{L^\infty},
\]
by the interpolation theorem, we have
\[
\bigl \|P_t^1f\bigr \|_{L^p}\leq e^{-2\lambda_1t/p}\|f
\|_{L^p},
\]
which implies that
\[
\bigl \|P^0_tf\bigr \|_{L^p}\leq\|f\|_{L^p}.
\]
Then we get that
\[
\frac{d}{dt}\bigl \|P^0_t\theta\bigr \|_{L^p}^p=
\int\bigl |P^0_t\theta \bigr |^{p-2}\bigl(P^0_t
\theta \bigr) \biggl(P^0_t\biggl(-\kappa
\Lambda^{2\alpha}+\frac{2\lambda_1}{p}\biggr)\theta\biggr)\, dx\leq0.
\]
Letting $t\rightarrow0$, we obtain the result.
\end{pf}
%
%$\hfill\Box$

%pr5.6 #&#
%
\begin{Proposition}\label{pro5.1.5}
Let $\alpha>\frac{1}{2}$. Suppose Hypothesis~\ref{hypE.1} holds.
For $x\in L^p$, let $\theta$ denote the solution of equation (\ref
{eq3.1}) with
the initial value $x$. Then for $2<p<\infty$
\[
E\bigl \|\theta(t)\bigr \|_{L^p}^p\leq\|x\|_{L^p}^pe^{-\lambda_1t}+C_S^p
\bigl[\tfrac
{1}{2}p(p-1)\bigr]^{p/2}\lambda_1^{-{p}/{2}}
\mathcal {E}_0^{p/2}\bigl(1-e^{-\lambda_1t}\bigr),
\]
where $C_S$ is the constant for the Sobolev embedding.
\end{Proposition}

\begin{pf}
Using \cite{Kr10}, Lemma~5.1, or \cite{BVVL08}, Theorem~2.4, for
$\theta_n$,
we obtain
%
%e5.1 #&#
%
\begin{eqnarray}\label{eq5.1}
\bigl \|\theta(t)\bigr \|_{L^p}^p&=& \bigl \|\theta(s)\bigr \|
_{L^p}^p\nonumber
\\
&&{}+\int_s^t
\biggl[-p\int_{\mathbb{T}^2}\bigl |\theta(l)\bigr |^{p-2}\theta(l) \bigl(
\kappa\Lambda^{2\alpha}\theta(l)+u(l)\cdot\nabla\theta(l)\bigr)\, d\xi
\nonumber
\\
&&\qquad\ \ {} +\frac{1}{2}p(p-1)\int_{\mathbb
{T}^2}\bigl |
\theta(l)\bigr |^{p-2}\biggl(\sum_j\bigl |k_{\delta
_n}*G(e_j)\bigr |^2
\biggr)\,d\xi\biggr]\,dl \nonumber
\\
&&{}+p\int_s^t \int
_{\mathbb{T}^2}\bigl |\theta(l)\bigr |^{p-2}\theta(l) k_{\delta_n}*G
\,d\xi \,dW(l)
\nonumber
\\
&\leq& \bigl \|\theta(s)\bigr \|_{L^p}^p-2 \lambda_1\int
_s^t\int_{\mathbb{T}^2}\bigl |
\theta(l)\bigr |^p\,d\xi\,dl
\nonumber
\\
&&{}+\int_s^t
\frac
{1}{2}p(p-1)\int_{\mathbb{T}^2}\bigl |\theta(l)\bigr |^{p-2}
\biggl(\sum_j\bigl |k_{\delta
_n}*G(e_j)\bigr |^2
\biggr)\,d\xi\,dl
\\
&&{}+p\int_s^t \int_{\mathbb{T}^2}\bigl |
\theta(l)\bigr |^{p-2}\theta(l)k_{\delta_n}*G\,d\xi \,dW(l)\nonumber
\\
&\leq& \bigl \|\theta(s)\bigr \|_{L^p}^p-2\lambda_1\int
_s^t\int_{\mathbb
{T}^2}\bigl |\theta
(l)\bigr |^p\,d\xi\,dl\nonumber
\\
&&{} +\int_s^t\biggl(
\lambda_1\int_{\mathbb{T}^2}\bigl |\theta(l)\bigr |^p\,d\xi
\nonumber
\\
&&\phantom{+\int_s^t\biggl(}{}+\biggl[\frac
{1}{2}p(p-1)\biggr]^{p/2}\lambda_1^{-({p-2})/{2}}
\int\biggl(\sum_j\bigl |k_{\delta
_n}*G(e_j)\bigr |^2
\biggr)^{p/2}\,d\xi\biggr)\,dl
\nonumber
\\
&&{}+p\int_s^t \int_{\mathbb{T}^2}\bigl |
\theta(l)\bigr |^{p-2}\theta(l) k_{\delta_n}*G\,d\xi \,dW(l),
\nonumber
\end{eqnarray}
where we used Lemma~\ref{lem5.1.4} to get the first inequality and Young's
inequality to get the last inequality. Here, for simplicity, we write
$\theta(t)=\theta_n(t,x)$. Taking expectation, we obtain
\begin{eqnarray*}
E\bigl \|\theta_n(t)\bigr \|_{L^p}^p&\leq& E\bigl \|
\theta_n(s)\bigr \|_{L^p}^p-E\lambda_1
\int_s^t\int_{\mathbb{T}^2}\bigl |
\theta_n(l)\bigr |^p\,d\xi\,dl
\\
&&{}+C_S^p
\biggl[\frac
{1}{2}p(p-1)\biggr]^{p/2}\lambda_1^{- ({p-2})/{2}}
\mathcal{E}_0^{p/2}(t-s).
\end{eqnarray*}
Here, we use
$\int_{\mathbb{T}^2}(\sum_j|G(e_j)|^2)^{p/2}\,d\xi\leq(\sum_j(\int_{\mathbb{T}^2}|G(e_j)|^p\,d\xi)^{2/p})^{{p}/{2}}\leq C_S^p\mathcal
{E}_0^{p/2}$.
Then Gronwall's lemma yields that
\[
E\bigl \|\theta_n(t)\bigr \|_{L^p}^p\leq \bigl \|
\theta_n(0)\bigr \| _{L^p}^pe^{-\lambda_1t}+C_S^p
\bigl[\tfrac{1}{2}p(p-1)\bigr]^{p/2}\lambda_1^{-
{p}/{2}}
\mathcal{E}_0^{p/2}\bigl(1-e^{-\lambda_1t}\bigr).
\]
Letting $n\rightarrow\infty$ in the above inequality, we deduce
\[
E\bigl \|\theta(t)\bigr \|_{L^p}^p\leq\|x\|_{L^p}^pe^{-\lambda_1t}+C_S^p
\bigl[\tfrac
{1}{2}p(p-1)\bigr]^{p/2}\lambda_1^{-{p}/{2}}
\mathcal {E}_0^{p/2}\bigl(1-e^{-\lambda_1t}\bigr).
\]
\upqed
\end{pf}
%
%$\hfill\Box$

%s5.2 #&#
\subsection{Uniqueness of the invariant measure}

In this subsection, we assume conditions Hypotheses~\ref{hypE.1}
and~\ref{hypE.2}
to hold. To
prove uniqueness of invariant measure is much harder and in this section we
first concrete on proving this. Existence will be shown in the next
subsection. In addition, we shall prove polynomial convergence of the
semigroup to the invariant measure in Section~\ref{sec5.3} below. If the
dissipation term is strong enough (i.e., $\alpha>\frac{2}{3})$ we
actually obtain exponential convergence (see Section~\ref{sec6}).

Now we build an auxiliary process $\tilde{\theta}$. The aim is to find
a shift $h$ belonging to Cameron--Martin space of the driving process
such that $E\|\theta(t)-\tilde{\theta}(t)\|_{H^{-1/2}}\rightarrow0$ as
$t\rightarrow\infty$. Fix $\theta$, and consider
%
%e5.2 #&#
%
\begin{equation}
\label{eq5.2} \cases{ %
d\tilde{\theta}(t)+A_\alpha
\tilde{\theta}(t)\,dt+\tilde {u}(t)\cdot\nabla\tilde{\theta}(t)
\,dt+K_0P_N\bigl(\tilde{\theta }-\theta (t,W,
\theta_0)\bigr)\,dt
\cr
\qquad= G\,dW(t),\vspace*{3pt}
\cr
\tilde{\theta}(0)=\tilde{\theta
}_0\in H^1, %\end{array}
}
\end{equation}
where $\tilde{u}$ satisfies (\ref{eq1.3}) with $\theta$ replaced by
$\tilde
{\theta}$ and $K_0$ is a constant to be determined later.
Since $\|P_N\tilde{\theta}\|_{L^p}\leq C_N \|\tilde{\theta}\|
_{L^p}$ for
$p\geq2$, by a similar argument as in the proof of Theorems~\ref
{theA.4} in
Appendix~\hyperref[appA]{A} we obtain
that there exists a measurable mapping from $\mathbb{R}^+\times
C(\mathbb{R}^+,H^{-1-\varepsilon})\times H^1\times H^1$ into $H^1$,
$(t,\overline{W},\theta_0,\tilde{\theta}_0)\mapsto\tilde{\theta
}(t,\overline{W},\theta_0,\tilde{\theta}_0)$, such that $\tilde
{\theta
}(t,W,\theta_0,\tilde{\theta}_0)$ is the solution of (\ref{eq5.2}).
Moreover, by the $\omega$-wise uniqueness of (\ref{eq3.1}) and (\ref
{eq5.2}) (which
can be easily checked by a similar argument as the proof of
Theorem~\ref{the4.2}), we have
\[
\bigl(\theta(t,\theta_0),\tilde{\theta}(t,\theta_0,
\tilde{\theta}_0)\bigr) =\bigl(\theta\bigl(t,s,\theta(s)\bigr),
\tilde{\theta}\bigl(t,s,\theta(s,\theta _0),\tilde {\theta}(s,
\theta_0, \tilde{\theta}_0)\bigr)\bigr)\qquad P
\mbox{-a.s.},
\]
which implies that $(\theta(t),\tilde{\theta}(t))=(\theta
(t,W,\theta
_0),\tilde{\theta}(t,W,\theta_0,\tilde{\theta}_0))$ defines a Markov
process. Here, for simplicity, we omit $W$ and $\theta(t,s,\theta(s))$,
$\tilde{\theta}(t,s,\theta(s,\theta_0),\allowbreak  \tilde{\theta}(s,\theta_0,
\tilde{\theta}_0))$ denote the solutions to (\ref{eq3.1}), (\ref
{eq5.2}) starting from
$\theta(s),\tilde{\theta}(s) $ at time $s$, respectively.

Now we derive a uniform $|\cdot|^4$ estimate for $\tilde{\theta}$.
Here, we give formal calculations which can be made rigorous by using
Galerkin approximations:
\begin{eqnarray*}
&&d\bigl |\tilde{\theta}(t)\bigr |^4+4\kappa\bigl |\tilde{\theta}(t)\bigr |^{2}
\| \tilde {\theta}\|_{H^\alpha}^2\,dt+4K_0\bigl |\tilde{
\theta}(t)\bigr |^{2}|P_N\tilde {\theta }|^2\,dt
\\
&&\qquad\leq 4\bigl |\tilde{\theta}(t)\bigr |^{2}\bigl\langle G\,dW(t),\tilde {
\theta }\bigr\rangle+4K_0\bigl |\tilde{\theta}(t)\bigr |^{2}|P_N
\tilde{\theta}||\theta |\,dt+6|\tilde{\theta}|^{2}\|G
\|^2_{L_2(H,H)}\,dt
\\
&&\qquad\leq4\bigl |\tilde{\theta}(t)\bigr |^{2}\bigl\langle G\,dW(t),\tilde{
\theta }\bigr\rangle +\varepsilon\bigl |\tilde{\theta}(t)\bigr |^{4}\,dt+C(
\varepsilon) \bigl(|\theta |^4+1\bigr)\,dt .
\end{eqnarray*}
Taking expectation and by Proposition~\ref{pro5.1.5}, we obtain
%
%e5.3 #&#
%
\begin{equation}
\label{eq5.3} E\bigl |\tilde{\theta}(t)\bigr |^4\leq C,\qquad\forall t\geq0,
\end{equation}
where $C$ is a constant independent of $t$.

Define $h(\theta,\tilde{\theta}):=-gK_0P_N(\tilde{\theta}-\theta
)$ for
$g$ in Hypothesis~\ref{hypE.2}. Then for any $(t,\theta_0,\tilde
{\theta}_0)\in{\mathbb
{R}}^+\times H^1\times H^1$ and the cylindrical Wiener process $W$ we
have for $\omega\in\Omega'$ that
$z(W(\omega)), z(W(\omega)+\int_0^\cdot h(\theta(s,W(\omega
),\theta
_0)$, $\tilde{\theta}(s,W(\omega),\theta_0,\tilde{\theta}_0))\,
ds)\in
C([0,\infty),H^{2+\varepsilon}), \varepsilon<\sigma$.
Then for $\omega\in\Omega'$,
\[
\theta\biggl(t,W(\omega)+\int_0^\cdot h\bigl(
\theta\bigl(s,W(\omega),\theta _0\bigr),\tilde {\theta}\bigl(s,W(
\omega),\theta_0,\tilde{\theta}_0\bigr)\bigr)\,ds,\tilde
{\theta }_0\biggr)-z\bigl(W(\omega)\bigr)
\]
is a solution to the following equation:
\[
d\tilde{v}(t)+A_\alpha\tilde{v}(t)\,dt+u_{\tilde{v}+z}(t)\cdot \nabla
(\tilde{v}+z) (t)\,dt+K_0P_N\bigl(\tilde{v}-v(t,W,
\theta_0)\bigr)\,dt=0,
\]
where $u_{\tilde{v}+z}$ satisfies (\ref{eq1.3}) with $\theta$
replaced by
$\tilde{v}+z$.
Since for every $\omega\in\Omega'$ the above equation admits at most
one solution, for $\omega\in\Omega'$ we have
%
%e5.4 #&#
%
\begin{eqnarray}
\label{eq5.4} &&\tilde{\theta}\bigl(t,W(\omega),\theta_0,\tilde{
\theta }_0\bigr)
\nonumber
\\[-8pt]
\\[-8pt]
&&\qquad=\theta\biggl(t,W(\omega)+\int_0^\cdot
h\bigl(\theta\bigl(s,W(\omega),\theta _0\bigr),\tilde{\theta}
\bigl(s,W(\omega),\theta_0,\tilde{\theta}_0\bigr)\bigr)
\, ds,\tilde {\theta}_0\biggr).\nonumber
\end{eqnarray}

Now for $\rho=\tilde{\theta}(t,W,\theta_0,\tilde{\theta
}_0)-\theta
(t,W,\theta_0)$, we have the following results. Here, we want to
emphasize that although the initial value $\theta_0\in H^1$, we can
only obtain that $\rho$ converges to $0$ in $H^{-1/2}$ norm.

\begin{Proposition}\label{pro5.2.1}
Fix $\alpha>1/2$. Let $\delta_0:=\lambda
_{N+1}-2^{p/2}C_R^p C_S^{2p}\kappa^{1-p}[p(p-1)] ^{p/2}\lambda
_1^{-p/2}\mathcal{E}_0^{{p} /2}>0$ for $p=\frac{\alpha+1}{\alpha-
{1}/{2}}$, where $N$ is as in Hypothesis~\ref{hypE.2},
and $C_S$, $C_R$ are the
constants for the Sobolev embedding and Riesz transform, respectively.
Then for $\|\theta_0\|_{L^{2m(p-1)}}^{2m(p-1)}+\|\tilde{\theta}_0\|
_{L^{2m(p-1)}}^{2m(p-1)}\leq2C_0$ for some $m>5$, $K_0>\lambda_{N+1}$
and $1<q<\frac{m-1}{4}$, there exists a positive constant $\overline
{C}$ such that for any $t>0$
\[
E\bigl |\Lambda^{-1/2}\rho(t)\bigr |^2\leq\frac{\overline{C}}{(t+1)^{2q}}
\]
(where we can choose $C_0$ large enough such that $C_0>4C_S^p[\frac
{1}{2}p(p-1)]^{p/2}\*\lambda_1^{-{p}/{2}}\mathcal{E}_0^{p/2}$).
\end{Proposition}

\begin{Remark}\label{rem5.2.2}
From the condition $\lambda_{N+1}-2^{p/2}C_R^p
C_S^{2p}\kappa^{1-p}[p(p-1)] ^{p/2}\*\lambda_1^{-p/2}\mathcal{E}_0^{{p}
/2}>0$, which also appears in the main theorem, we know that if the
viscosity constant $\kappa$ is large enough or $\mathcal{E}_0$ is small
enough we could even take $N=0$.
\end{Remark}

\begin{pf*}{Proof of Proposition~\ref{pro5.2.1}}
In the proof, we omit $W$ for
simplicity. From (\ref{eq3.1}) and (\ref{eq5.2}), we obtain that
$\rho$ satisfies the
following equation in the weak sense:
\begin{eqnarray*}
\frac{d\rho(t)}{dt}&=&-A_\alpha\rho-K_0P_N
\rho-\tilde{u}\cdot \nabla\tilde{\theta}+u\cdot\nabla\theta
\\
&=&-A_\alpha\rho-K_0P_N\rho-u\cdot\nabla
\rho-u_\rho\cdot\nabla \tilde {\theta},
\end{eqnarray*}
where $u_\rho$ satisfies (\ref{eq1.3}) with $\theta$ replaced by
$\rho$.
Taking the inner product with $\Lambda^{-1}\rho$ in $H$, and using that
\[
{ }_{H^{-1}} \bigl\langle u_\rho\cdot\nabla\tilde{\theta},
\Lambda ^{-1}\rho \bigr\rangle_{H^1}=0
\]
(cf. \cite{Re95}),
we obtain
\[
\frac{1}{2}\frac{d}{dt}\bigl |\Lambda^{-{1}/{2}}\rho\bigr |^2=-
\kappa\bigl |\Lambda ^{\alpha-{1}/{2}} \rho\bigr |^2-K_0\bigl |P_N
\Lambda^{-{1}/{2}}\rho\bigr |^2-{ }_{H^{-1}} \bigl\langle u\cdot
\nabla\rho, \Lambda^{-1}\rho\bigr\rangle_{H^1}.
\]

We have
\begin{eqnarray*}
&&\bigl |{ }_{H^{-1}} \bigl\langle u\cdot\nabla\rho, \Lambda^{-1}\rho
\bigr\rangle _{H^1}\bigr |
\\
&&\qquad \leq  \|u\|_{L^{p} }\|\rho\|_{L^{p_1} }
\bigl \|\nabla\Lambda ^{-1}\rho \bigr \|_{L^{p_1} }\leq C_S\|u
\|_{L^{p} }\|\rho\|_{H^{1/{p} }}\bigl \|\nabla \Lambda ^{-1}\rho
\bigr \|_{H^{1/{p} }}
\\
&&\qquad \leq C_SC_R\|\theta\|_{L^{p} }\bigl \|\Lambda
^{-1}\rho\bigr \|_{H^{1+{1}/{{p}}}}^2\leq C_SC_R
\|\theta\|_{L^{p} }\bigl \| \Lambda^{-1}\rho\bigr \|^{2/r}_{H^{{1}/{2}}}
\bigl \|\Lambda^{-1}\rho\bigr \| ^{2(1-{1}/{r})}_{H^{{1}/{2}+\alpha}}
\\
&&\qquad \leq  \frac{\kappa
}{2}\bigl |\Lambda^{\alpha-{1}/{2}}\rho\bigr |^2+
C_1^r\biggl(\frac{\kappa
}{2}\biggr)^{1-r}\|
\theta\|_{L^{p} }^r\bigl |\Lambda^{-{1}/{2}}\rho
\bigr |^2,
\end{eqnarray*}
where $C_S$, $C_R$ are the constants for Sobolev embedding and Riesz
transform, respectively, and $C_1=C_SC_R$. Here, $\frac{1}{{p}}+\frac
{2}{{p_1}}=1$ for $p>\frac{1}{\alpha-{1}/{2}}$, $r=\frac{\alpha
}{\alpha-{1}/{2}-{1}/{{p} }} $ and we use H\"{o}lder's
inequality and that $\operatorname{div}u=0$ in the first inequality
and $H^{1/{p}
}\hookrightarrow L^{p_1} $ continuously in the second inequality, the
interpolation inequality (\ref{eq2.3}) in the fourth inequality and Young's
inequality in the last equality. Then we obtain
\begin{eqnarray*}
\frac{d}{dt}\bigl |\Lambda^{-{1}/{2}}\rho\bigr |^2&\leq&-{\kappa}\bigl |
\Lambda ^{\alpha
-{1}/{2}} \rho\bigr |^2-K_0\bigl |P_N
\Lambda^{-{1}/{2}}\rho\bigr |^2
\\
&&{}+2C_1^r\biggl(
\frac
{\kappa}{2}\biggr)^{1-r}\|\theta\|_{L^{p} }^r\bigl |
\Lambda^{-{1}/{2}}\rho\bigr |^2.
\end{eqnarray*}
Since, because $K_0>\lambda_{N+1}$, we have
\begin{eqnarray*}
\lambda_{N+1}\bigl |\Lambda^{-{1}/{2}}\rho\bigr |^2&\leq&
\kappa\bigl |Q_N\Lambda ^{\alpha-{1}/{2}}\rho\bigr |^2+K_0\bigl |P_N
\Lambda^{-{1}/{2}}\rho\bigr |^2
\\
&\leq& \kappa\bigl |\Lambda^{\alpha-{1}/{2}}
\rho\bigr |^2+K_0\bigl |P_N\Lambda^{-
{1}/{2}}
\rho\bigr |^2,
\end{eqnarray*}
it follows that
\[
\frac{d}{dt}\bigl |\Lambda^{-{1}/{2}}\rho\bigr |^2+\biggl(
\lambda_{N+1}-2C_1^r\biggl(\frac
{\kappa}{2}
\biggr)^{1-r}\|\theta\|_{L^{p} }^r\biggr)\bigl |
\Lambda^{-{1}/{2}}\rho \bigr |^2\leq0.
\]

Thus, by Gronwall's lemma, we obtain
\[
\bigl |\Lambda^{-{1}/{2}}\rho(t)\bigr |^2\leq e^{t\Gamma(t,\theta_0)}\bigl |\Lambda
^{-{1}/{2}}\rho(0)\bigr |^2,
\]
where
\[
\Gamma(t,\theta_0)=-\lambda_{N+1}+2C_1^r
\biggl(\frac{\kappa
}{2}\biggr)^{1-r}\frac
{1}{t}\int
_0^t\bigl \|\theta(s)\bigr \|_{L^{p} }^r
\,ds.
\]

By the same arguments as in the proof of Theorem~\ref{theA.1}
in Appendix~\hyperref[appA]{A},
we have $\theta_n\rightarrow\theta$ in $L^2([0,T],H^1)$ a.s. Letting
$n\rightarrow\infty$ in (\ref{eq5.1}), by (\ref{eq3.11}) we obtain
\begin{eqnarray*}
&&\bigl \|\theta(t)\bigr \|_{L^{p} }^{p} +\lambda_1\int
_0^t\int_{\mathbb
{T}^2}\bigl |
\theta(l)\bigr |^{p} \,d\xi\,dl
\\
&&\qquad\leq\|\theta_0
\|_{L^{p} }^{p} +C_S^p\biggl[
\frac{1}{2}p(p-1)\biggr]^{p/2}\lambda_1^{-({p-2})/{2}}
\mathcal{E}_0^{p/2}t
\\
&&\qquad\quad{} +p\int_0^t
\int_{\mathbb{T}^2}\bigl |\theta(l)\bigr |^{p-2}\theta(l) G\,d\xi\,dW(l).
\end{eqnarray*}
Here, we use that
$\int_{\mathbb{T}^2}(\sum_j|G(e_j)|^2)^{p/2}\,d\xi\leq(\sum_j(\int_{\mathbb{T}^2}|G(e_j)|^p\,d\xi)^{2/p})^{{p}/{2}}\leq C_S^p\mathcal
{E}_0^{p/2}$.

Since $p=\frac{\alpha+1}{\alpha-{1}/{2}}$ implies $p=r$, we get
\begin{eqnarray*}
\Gamma(t,\theta_0) &\leq&-\lambda_{N+1}
+2C_1^p\biggl(\frac{\kappa
}{2}\biggr)^{1-p}
\frac{1}{t}\int_0^t\bigl \|\theta(s)
\bigr \|_{L^{p} }^{p} \,ds
\\
&\leq&-\lambda_{N+1}+2C_1^p\biggl(
\frac{\kappa}{2}\biggr)^{1-p}\frac
{1}{t\lambda_1}\|
\theta_0\|_{L^{p} }^{p}
\\
&&{}+2^{p/2}C_1^p
C_S^p\kappa^{1-p}\bigl[p(p-1)\bigr]
^{p/2}\lambda _1^{-p/2}\mathcal
{E}_0^{{p} /2}
\\
&&{}+2C_1^p\biggl(\frac{\kappa}{2}
\biggr)^{1-p}\frac{p}{t\lambda_1}\int_0^t
\int_{\mathbb{T}^2}\bigl |\theta(l)\bigr |^{{p} -2}\theta(l) G\,d\xi
\,dW(l).
\end{eqnarray*}
For $M(t):={p} \int_0^t
\int_{\mathbb{T}^2}|\theta(l)|^{{p} -2}\theta(l) G\,d\xi\,dW(l)$,
we have
\[
\langle M\rangle_t\leq{p}^2\mathcal{E}_0C_S^2
\int_0^t\biggl(\int_{\mathbb
{T}^2}\bigl |
\theta(s)\bigr |^{{p} -1}\,d\xi\biggr)^2\,ds,
\]
where we use that $\sum_j|G(e_j)|^2(\xi)\leq\sum_j\|G(e_j)\|
_{L^\infty
}^2\leq C_S^2\mathcal{E}_0$.
Then for any $m>1$
\begin{eqnarray*}
\langle M\rangle_t^m&\leq& C_S^{2m}{p}
^{2m}\mathcal{E}_0^m \biggl(\int
_0^t\biggl(\int_{\mathbb{T}^2}\bigl |
\theta(s)\bigr |^{{p} -1}\,d\xi\biggr)^2\,ds\biggr)^m
\\
&\leq& C_S^{2m}{p} ^{2m}\mathcal{E}_0^mt^{m-1}
\int_0^t\biggl(\int_{\mathbb{T}^2}\bigl |
\theta (s)\bigr |^{2m(p -1)}\,d\xi\biggr)\,ds.
\end{eqnarray*}
Since $\|\theta_0\|_{L^{2m({p} -1)}}^{2m(p-1)}\leq2C_0$ by
Proposition~\ref{pro5.1.5}
there exists a constant $C_{{p}, m}(C_0)$ independent of $t$ such
that $E\|\theta(t)\|_{L^{2m({p} -1)}}^{2m({p} -1)}\leq C_{{p}, m}$ for
$t\geq0$. Thus, for $M_n=\sup_{n-1\leq t<n}M(t)$, we have
\[
P\biggl(|M_n|>\frac{\varepsilon\lambda_1}{4C_1^p({\kappa}/{2})^{1-p}} n\biggr)\leq\frac{{p} ^{2m}\mathcal{E}_0^mC_{{p}, m}n^m C_S^{2m}}{(
{\varepsilon\kappa^{p-1}\lambda_1}/({2^{p+1}C_1^p}))^{2m}n^{2m}}.
\]
Now define the following random times:
\[
T_{\mathrm{bound}}:=\sup\biggl\{n\dvtx |M_n|>\frac{\varepsilon\lambda
_1}{4C_1^p({\kappa}/{2})^{1-p}} n
\biggr\}.
\]
By \cite{M99}, Lemma~5, we have that if $m>1$, then $T_{\mathrm{bound}}$ is
finite almost surely.
Set
\[
\tau:=\max\biggl(T_{\mathrm{bound}}, \frac
{2^{p+1}C_0^{p/(2m(p-1))}C_1^p}{\kappa
^{p-1}\lambda_1\varepsilon}\biggr),
\]
then we have
\[
t>\tau\quad\Rightarrow\quad\Gamma(t,\theta_0)-(-\delta
_0)<\varepsilon,
\]
where $\delta_0=\lambda_{N+1}
-2^{p/2}C_1^p C_S^p\kappa^{1-p}[p(p-1)] ^{p/2}\lambda
_1^{-p/2}\mathcal
{E}_0^{{p} /2}$, which implies that for $\delta\in(0,\delta_0)$ and
$t>\tau$,
\[
\bigl |\Lambda^{-1/2}\rho(t)\bigr |^2\leq\bigl |\Lambda^{-1/2}(
\theta_0-\tilde {\theta }_0)\bigr |^2e^{-\delta t}.
\]
For $p_0\in(0,m-1)$, by \cite{M99}, Lemma~5, $E\tau^{p_0}$ is finite.
Moreover, we obtain that for $1<q<\frac{m-1}{4}$, there exists
$\overline{C}>0$ such that for any $t>0$
%
%e5.5 #&#
%
\begin{eqnarray}
\label{eq5.5} && E\bigl |\Lambda^{-1/2}\rho(t)\bigr |^2
\nonumber
\\
&&\qquad\leq Ce^{-\delta
t}+\bigl(E\bigl |\Lambda^{-1/2}
\rho(t)\bigr |^4\bigr)^{{1}/{2}}P(\tau>t)^{{1}/{2}}
\\
&&\qquad\leq \overline{C}\frac{1}{(t+1)^{2q}},
\nonumber
\end{eqnarray}
where we used (\ref{eq5.3}) in the last inequality.
\end{pf*}
%
%$\hfill\Box$

Now we fix $m>35$ and $8<q<\frac{m-3}{4}$. Proposition~\ref{pro5.2.1} still
holds for such $m$, $q$. Moreover, we also have for any $t_0\geq0$
%
%e5.6 #&#
%
\begin{eqnarray}
\label{eq5.6} &&P\biggl(\int_{t_0}^{\infty}\bigl |h(t)\bigr |^2
\,dt\geq\overline {C}\frac{1}{(t_0+1)^q}\biggr)
\nonumber
\\
&&\qquad\leq C\frac{(t_0+1)^{q}}{\overline{C}}\int_{t_0}^\infty
E\bigl |\Lambda ^{-1/2}\rho(t)\bigr |^2\,dt
\\
&&\qquad\leq\overline{C}\frac{1}{(t_0+1)^q},
\nonumber
\end{eqnarray}
where $h(t)=h(\theta(t,W,\theta_0),\tilde{\theta}(t,W,\theta
_0,\tilde
{\theta}_0))$ and we used Proposition~\ref{pro5.2.1} in the last inequality.
Moreover, by Theorem~\ref{the5.2.3}, we obtain that there exists
$p_2>0$ such that
%
%e5.7 #&#
%
\begin{eqnarray}\label{eq5.7}
P\biggl(\int_0^\infty\bigl |h(t)\bigr |^2\geq
\overline{C}\biggr)&\leq& \frac{C_1}{\overline{C}}E\int_0^\infty\bigl |
\Lambda^{-{1}/{2}}\rho (t)\bigr |^2\,dt
\nonumber
\\[-8pt]
\\[-8pt]
&\leq&1-p_2,
\nonumber
\end{eqnarray}
where $\overline{C}$ can be chosen large enough such that (\ref
{eq5.5}), (\ref{eq5.6})
and (\ref{eq5.7}) are satisfied.

Now we use a similar coupling method as in \cite{O08} to deduce the
uniqueness of the invariant measure. More precisely, we have the
following result.

\begin{Theorem}\label{the5.2.3}
Fix $\alpha>1/2$. Assume Hypotheses~\ref{hypE.1} and~\ref{hypE.2}
hold. Let
$\delta_0:=\lambda_{N+1}-2^{p/2}C_R^p C_S^{2p}\kappa^{1-p}[p(p-1)]
^{p/2}\lambda_1^{-p/2}\mathcal{E}_0^{{p} /2}>0$ for $p=\frac{\alpha
+1}{\alpha-{1}/{2}}$, where $N$ is as in Hypothesis~\ref{hypE.2},
and $C_S,
C_R$ are the constants for Sobolev embedding and Riesz transform,
respectively. Then there exists at most one invariant measure for the
Markov semigroup $P_t$ on $ H^1$.
\end{Theorem}

\begin{pf}\textit{Step} 1. Construction of a coupling of the solutions.

For $\theta_0^1,\theta_0^2\in H^1$ and $T>0$, we apply \cite{O08},
Corollary~1.5, to $(\theta(\cdot,W,\theta_0^1),\theta(\cdot,\allowbreak W,  \theta
_0^2),\tilde
{\theta}(\cdot,W,\theta_0^1,\theta_0^2))$ on $[0,T]$
and obtain $(\theta_1^0(\cdot,\theta_0^1,\theta_0^2),\theta
_2^0(\cdot
,\theta_0^1,\theta_0^2),\tilde{\theta}^0(\cdot,\theta_0^1,\allowbreak  \theta_0^2))$
on $[0,T]$ such that the law of $(\theta_1^0(\cdot,\theta_0^1$,
$\theta
_0^2),\tilde{\theta}^0(\cdot,\theta_0^1,\theta_0^2))$ is the same as
$(\theta(\cdot,W,\theta_0^1),\tilde{\theta}(\cdot,W,\theta
_0^1,\theta
_0^2))$ and $(\theta_2^0(\cdot,\theta_0^1,\theta_0^2),\tilde
{\theta
}^0(\cdot,\theta_0^1,\theta_0^2))$ is a maximal coupling of
$(\mathcal
{D}(\theta(\cdot,W,\theta_0^2)),\mathcal{D}(\tilde{\theta}(\cdot
,W,\theta_0^1,\theta_0^2)))$ on $[0,T]$.

Then we obtain a sequence of independent versions of the mapping
\[
\bigl(\theta_0^1,\theta_0^2
\bigr)\rightarrow\bigl(\theta_1^0\bigl(\cdot,\theta
_0^1,\theta _0^2\bigr),
\theta_2^0\bigl(\cdot,\theta_0^1,
\theta_0^2\bigr),\tilde{\theta }^0\bigl(
\cdot ,\theta_0^1,\theta_0^2
\bigr)\bigr).
\]
We denote this sequence by $(\theta_1^n,\theta_2^n,\tilde{\theta}^n)_n$
and define recursively
\[
\cases{ %
\theta_1\bigl(nT+\cdot,
\theta_0^1,\theta_0^2\bigr)=
\theta _1^n\bigl(\cdot,\theta_1(nT),
\theta_2(nT)\bigr),
\cr
\theta_2\bigl(nT+\cdot ,
\theta_0^1,\theta_0^2\bigr)=
\theta_2^n\bigl(\cdot,\theta_1(nT),
\theta_2(nT)\bigr),
\cr
\tilde{\theta}\bigl(nT+\cdot,
\theta_0^1,\theta_0^2\bigr)=
\tilde{\theta }^n\bigl(\cdot,\theta_1(nT),
\theta_2(nT)\bigr). %\end{array}
}
\]
Then $\theta_1(t,\theta_0^1,\theta_0^2),\theta_2(t,\theta
_0^1,\theta
_0^2),\tilde{\theta}(t,\theta_0^1,\theta_0^2)$ is defined for all
$t\in
[0,\infty)$ such that $(\theta_1(\cdot,\theta_0^1,\theta
_0^2),\theta
_2(\cdot,\theta_0^1,\theta_0^2))$ is a coupling of $(\mathcal
{D}(\theta
(\cdot,W,\theta_0^1)),\mathcal{D}(\theta(\cdot,W,\theta_0^2)))$. We
denote the associated probability space by $(\Omega,\mathcal
{F},\mathcal
{F}_t, P)$. Moreover, $(\theta_1(nT,\allowbreak\theta_0^1,  \theta_0^2),\theta
_2(nT,\theta_0^1,\theta_0^2),\tilde{\theta}(nT,\theta_0^1,\theta
_0^2))_n$ is a Markov chain and $\theta_1 (\cdot,\theta_0^1,\theta
_0^2),\theta_2(\cdot,\theta_0^1,\allowbreak  \theta_0^2),\tilde{\theta}(\cdot
,\theta
_0^1,\theta_0^2)$ satisfy the following property:
\[
E^{(\theta_0^1,\theta_0^2)}\bigl[f(\theta_1,\theta_2,\tilde{\theta
})\circ \Phi_{kT}|\mathcal{F}_{kT}\bigr]=E^{(\theta_1(kT),\theta_2(kT))}f(
\theta _1,\theta_2,\tilde{\theta}),
\]
where $\Phi_t$ is the shift operator.

\textit{Step} 2. Introduction of $l_0$.

We set
\[
l_0(k)=\min\{l\leq k|P_{l,k}\},
\]
where $\min\varnothing=\infty$ and
\[
(P_{l,k})\cases{ %
\tilde{\theta}\bigl(\cdot,
\theta_0^1,\theta_0^2\bigr)=
\theta _2\bigl(\cdot,\theta_0^1,
\theta_0^2\bigr)\qquad\mbox{on } (lT,kT),
\cr
\bigl \|\theta
_1\bigl(lT,\theta_0^1,\theta_0^2
\bigr)\bigr \|_{L^{2m(p-1)}}^{2m(p-1)}+\bigl \|\theta _2\bigl(lT,
\theta_0^1,\theta_0^2\bigr)
\bigr \|_{L^{2m(p-1)}}^{2m(p-1)}\leq2C_0. %\end{array}
}
\]
Then by (\ref{eq5.5}) and the Markov property of $\theta_1(\cdot
,\theta
_0^1,\theta_0^2),\theta_2(\cdot,\theta_0^1,\theta_0^2),\tilde
{\theta
}(\cdot,\theta_0^1,\theta_0^2)$ we have for $t> lT$
\begin{eqnarray*}
&& E\bigl(\bigl |\Lambda^{-{1}/{2}}\bigl(\theta_2\bigl(t,
\theta_0^1,\theta _0^2\bigr)-
\theta_1\bigl(t,\theta_0^1,
\theta_0^2\bigr)\bigr)\bigr | 1_{l_0(\infty)\leq l}\bigr)
\\
&& \qquad= \sum_{k=0}^lE\bigl(\bigl |
\Lambda^{-{1}/{2}}\bigl(\theta_2\bigl(t,\theta_0^1,
\theta _0^2\bigr)-\theta_1\bigl(t,
\theta_0^1,\theta_0^2\bigr)
\bigr)\bigr | 1_{l_0(\infty)=k}\bigr)
\\
&&\qquad= \sum_{k=0}^lE \bigl[E\bigl(\bigl |
\Lambda^{-{1}/{2}}\bigl(\theta_2\bigl(t-kT+kT,\theta
_0^1,\theta_0^2\bigr)
%&&\phantom{\hspace*{107pt}}
%{}
-
\theta_1\bigl(t-kT+kT,\theta_0^1,
\theta_0^2\bigr)\bigr)\bigr |
\\
&&\phantom{\hspace*{269pt}}{}\cdot 1_{l_0(\infty
)=k}|
\mathcal{F}_{kT}\bigr) \bigr]
\\
&&\qquad= \sum_{k=0}^lE\bigl[E^{(\theta_1(kT),\theta
_2(kT))}\bigl[\bigl |
\Lambda^{-{1}/{2}}\bigl(\theta_2\bigl(t-kT,\theta_1(kT),
\theta_2(kT)\bigr)
\\
&&\phantom{\hspace*{165pt}}{}- \theta_1\bigl(t-kT,
\theta_1(kT),\theta_2(kT)\bigr)\bigr)\bigr |
\\
&&\phantom{\hspace*{130pt}} {}\cdot1_{\{\theta_2(\cdot
-kT,\theta_1(kT),\theta_2(kT))
=\tilde{\theta}(\cdot-kT,\theta_1(kT),\theta_2(kT))\}}\bigr]
\\
&&\phantom{\hspace*{155pt}} {}\cdot1_{\{\|
\theta_1(kT)\|^{2m(p-1)}_{L^{2m(p-1)}}+\|\theta_2(kT)\|
^{2m(p-1)}_{L^{2m(p-1)}}\leq2C_0\}}\bigr]
\\
&&\qquad \leq\sum_{k=0}^lE
\bigl[E^{(\theta_1(kT),\theta_2(kT))}\bigl[\bigl |\Lambda ^{-
{1}/{2}}\bigl(\tilde{\theta}
\bigl(t-kT,W,\theta_1(kT),\theta_2(kT)\bigr)
\\
&&\phantom{\hspace*{190pt}}{}- \theta
\bigl(t-kT,W,\theta_1(kT)\bigr)\bigr)\bigr |\bigr]
\\
&&\phantom{\hspace*{137pt}} {}\cdot1_{\{\|\theta_1(kT)\|
^{2m(p-1)}_{L^{2m(p-1)}}+\|\theta_2(kT)\|^{2m(p-1)}_{L^{2m(p-1)}}\leq
2C_0\}}\bigr]
\\
&&\qquad\leq\overline{C}\sum_{k=0}^l(t-kT+1)^{-q}
\leq C(t-lT+1)^{-q+1},
\end{eqnarray*}
where we used $\theta_i(kT)$ to denote $\theta_i(kT,\theta
_0^1,\theta
_0^2)$ for simplicity.

\textit{Step} 3. Construction of Wiener processes.

Now we want to estimate $P(l_0(k+1)=0|l_0(k)=0)$. As in most papers
using coupling methods for SPDEs, our tool is the Girsanov transform.
Set
\[
\cases{ %
h(t,W)=h\bigl(\theta\bigl(t-kT,W,
\theta_1\bigl(kT,\theta_0^1,\theta
_0^2\bigr)\bigr),
\cr
\phantom{\hspace*{58pt}}\tilde{\theta}\bigl(t-kT,W,
\theta_1\bigl(kT,\theta_0^1,
\theta_0^2\bigr), \theta_2\bigl(kT,
\theta_0^1,\theta_0^2\bigr)
\bigr)\bigr),
\cr
\tau_1(W)=\inf\biggl\{t\in \bigl(kT,(k+1)T\bigr]\Big|\displaystyle\int
_{kT}^t\bigl |h(t,W)\bigr |^2\,dt>
\overline{C}(kT+1)^{-q}\biggr\}. %\end{array}
}
\]

Then by Proposition~\ref{pro5.1.3}, we obtain cylindrical Wiener processes
$W_1,W_2$ on $(\tilde{\Omega},\tilde{\mathcal{F}},\tilde{P})$ such that
\[
\biggl(W_2,W_1+\int_{kT}^{\tau_1(W_1)\wedge\cdot}h(t,W_1)
\,dt\biggr),
\]
is a maximal coupling of $(\mathcal{D}(W),\mathcal{D}(W+\int_{kT}^{\tau
_1(W)\wedge\cdot}h(t,W)\,dt))$ on $[kT,(k+1)T]$.

If $l_0(k)=0$, by construction in Step 1, we have
%
%e5.8 #&#
%
\begin{eqnarray}
\label{eq5.8} %\begin{array}{cc}
&& P\bigl(l_0(k+1)=0|
\mathcal{F}_{kT}\bigr)
\nonumber
\\
&&\qquad= P\bigl(\tilde {\theta}\bigl(t,\theta_0^1,
\theta_0^2\bigr)=\theta_2\bigl(t,
\theta_0^1,\theta_0^2\bigr)
\mbox{ for }t\in\bigl[kT,(k+1)T\bigr]|\mathcal{F}_{kT}\bigr)
\nonumber
\\
&&\qquad\geq \tilde{P}\bigl(\tilde{\theta}\bigl(\cdot-kT,W_1,\theta
_1\bigl(kT,\theta _0^1,\theta_0^2
\bigr),\theta_2\bigl(kT,\theta_0^1,
\theta_0^2\bigr)\bigr)
\\
&&\qquad=\theta\bigl(\cdot
-kT,W_2,\theta_2\bigl(kT,\theta_0^1,
\theta_0^2\bigr)\bigr)\mbox{ for }t\in \bigl[kT,(k+1)T
\bigr]\bigr)
\nonumber
\\
&&\qquad\geq\tilde{P}\biggl(W_2=W_1+\int
_{kT}^{\tau_1(W_1)\wedge\cdot
}h(t,W_1)\,dt \mbox{ and }
\tau_1(W_1)=(k+1)T\biggr),
\nonumber
\end{eqnarray}
where we used that $(\tilde{\theta}(\cdot,\theta_0^1,\theta_0^2),
\theta
_2(\cdot,
\theta_0^1,\theta_0^2))$ is a maximal coupling of $(\tilde{\theta
}(\cdot
-kT,W_1,\theta_1(kT,\theta_0^1,\theta_0^2) , \theta_2(kT,\theta
_0^1,\theta_0^2)) , \theta(\cdot-kT,W_2,\theta_2(kT,\theta
_0^1,\theta
_0^2)))$ in the first inequality and (\ref{eq5.4}) in the last inequality.

Now set $A:=\{W|\tau_1(W)=(k+1)T\}, \Lambda_1:=\mathcal{D}(W),
\Lambda
_2:=\mathcal{D}(W+\int_{kT}^{\tau_1(W)\wedge\cdot}\* h(t,W)\,dt)$.
Then the Novikov condition is satisfied for $\Lambda_1$ and $\Lambda
_2$, which by the Girsanov transform implies that
\[
\biggl(\frac{d\Lambda_1}{d\Lambda_2}\biggr) (W)=\exp\biggl(-\int_{kT}^{\tau
_1(W)}h(t,W)
\,dW(t)-\frac{1}{2} \int_{kT}^{\tau_1(W)}\bigl |h(t,W)\bigr |^2
\,dt\biggr).
\]
Thus, we have
\[
\int\biggl(\frac{d\Lambda_1}{d\Lambda_2}\biggr)^2\,d\Lambda_1\leq E
\bigl(M_2e^{\int
_{kT}^{\tau_1(W)}|h(t,W)|^2\,dt}\bigr)\leq e^{\overline{C}(kT+1)^{-q}},
\]
where $M_2=\exp(-2\int_{kT}^{\tau_1(W)}h(t,W)\,dW(t)-2
\int_{kT}^{\tau_1(W)}\!|h(t,W)|^2\,dt)$ and \mbox{$EM_2\leq1$}.
By this, (\ref{eq5.7}), (\ref{eq5.8}) and Lemmas~\ref{lem5.1.1}
and~\ref{lem5.1.2},
we obtain
%
%e5.9 #&#
%
\begin{eqnarray}
\label{eq5.9} P\bigl(l_0(1)=0\bigr)&\geq&(\Lambda_1\wedge
\Lambda_2) (A)
\nonumber
\\[-8pt]
\\[-8pt]
&\geq&\frac
{1}{4}\biggl(\int\biggl(
\frac{d\Lambda_1}{d\Lambda_2}\biggr)^2\,d\Lambda _1
\biggr)^{-1}\Lambda _1(A)^2\geq
\frac{p_2^2}{4}e^{-\overline{C}}.\nonumber
\end{eqnarray}

\textit{Step} 4. Estimate for $P(l_0(k+1)\neq0, l_0(k)=0)$.

By (\ref{eq5.8}), we obtain
\begin{eqnarray*}
&&P\bigl(l_0(k+1)\neq0|\mathcal{F}_{kT}\bigr)
\\
&&\qquad
\leq \tilde{P}\biggl(W_2=W_1+\int_{kT}^{\tau_1(W_1)\wedge\cdot}h(t,W_1)
\,dt \mbox{ and }\tau _1(W_1)<(k+1)T\biggr)
\\
&&\qquad\quad{}+\tilde{P}\biggl(W_2\neq W_1+\int
_{kT}^{\tau_1(W_1)\wedge
\cdot}h(t,W_1)\,dt\biggr).
\end{eqnarray*}
Since $(W_2,W_1+\int_0^{\tau_1(W_1)\wedge\cdot}h(t,W_1)\,dt)$ is a
maximal coupling, it follows from Lemma~\ref{lem5.1.1} and the
construction of
$\tau_1$ that
%
%e5.10 #&#
%
\begin{eqnarray}\label{eq5.10}
&&\tilde{P}\biggl(W_2\neq W_1+\int
_{kT}^{\tau
_1(W_1)\wedge\cdot}h(t,W_1)\,dt\biggr)
\nonumber
\\
&&\qquad=\|\Lambda_1-\Lambda_2\|_{\mathrm{var}}
\nonumber
\\[-8pt]
\\[-8pt]
&&\qquad\leq\frac{1}{2}\sqrt{\int\biggl(\frac{d\Lambda_1}{d\Lambda
_2}
\biggr)^2\,d\Lambda _2-1}\leq\frac{1}{2}\sqrt{\int
\biggl(\biggl(\frac{d\Lambda_1}{d\Lambda
_2}\biggr)^2\,d\Lambda _1
\biggr)^{1/2}-1}
\nonumber
\\
&&\qquad\leq e^{{\overline{C}}/{4}}(kT+1)^{-
{q}/{2}}.
\nonumber
\end{eqnarray}
Since by the Markov property of $(\theta_1(\cdot,\theta_0^1,\theta
_0^2),\tilde{\theta}(\cdot,\theta_0^1,\theta_0^2))$, we have
\begin{eqnarray*}
&&\tilde{P}\biggl(W_2=W_1+\int
_{kT}^{\tau_1(W_1)\wedge\cdot}h(t,W_1)\,dt \mbox{ and }
\tau_1(W_1)<(k+1)T\biggr)
\\
&&\qquad\leq\tilde{P}\Bigl(\theta\bigl(\cdot -kT,W_2,
\theta_2\bigl(kT,\theta_0^1,
\theta_0^2\bigr)\bigr)
\\
&&\phantom{\qquad\leq\tilde{P}\bigl(}=\tilde{\theta}\bigl(\cdot
-kT,W_1,\theta_1\bigl(kT,\theta_0^1,
\theta_0^2\bigr),\theta_2\bigl(kT,\theta
_0^1,\theta _0^2\bigr)\bigr)
\\
&&\phantom{\hspace*{151pt}}\mbox{and }\tau_1(W_1)<(k+1)T\Bigr)
\\
&&\qquad\leq\tilde{P}\biggl(\int_{kT}^{(k+1)T}\bigl |h(t,W_1)\bigr |^2
\,dt>\overline {C}(kT+1)^{-q}\biggr)
\\
&&\qquad\leq P\biggl(\int_{kT}^{(k+1)T}\bigl |h\bigl(
\theta\bigl(t,W,\theta_0^1\bigr),\tilde {\theta }
\bigl(t,W,\theta_0^1,\theta_0^2
\bigr)\bigr)\bigr |^2\,dt>\overline{C}(kT+1)^{-q}\biggr),
\end{eqnarray*}
by (\ref{eq5.6}) and (\ref{eq5.10}), we obtain
%
%e5.11 #&#
%
\begin{equation}
\label{eq5.11} P\bigl(l_0(k+1)\neq0 \mbox{ and }
l_0(k)=0\bigr)\leq C(kT+1)^{-{q}/{2}},
\end{equation}
where $C$ depends on $\overline{C}$.

\textit{Step} 5. Estimate for $El_0(\infty)^q$.

Since $l_0(k)=0$ implies $l_0(l)=0$ for any $0\leq l\leq k\leq\infty$,
\[
P\bigl(l_0(\infty)\neq0\bigr)\leq\sum_{k=0}^\infty
P\bigl(l_0(k+1)\neq0 \mbox{ and } l_0(k)=0\bigr).
\]
By (\ref{eq5.9}) and (\ref{eq5.11}), we obtain
\[
P\bigl(l_0(\infty)\neq0\bigr)\leq1-\frac{p_2^2}{4}e^{-\overline{C}}+C
\sum_{k=1}^\infty(kT+1)^{-{q}/{2}}.
\]
Then there exists $T_0$ such that for $T\geq T_0$ we have
%
%e5.12 #&#
%
\begin{equation}
\label{eq5.12} P\bigl(l_0(\infty)=0\bigr)\geq p_0=
\frac{p_2^2}{8}e^{-\overline{C}}.
\end{equation}
Now fix $T=T_0$. Define
\[
\sigma:=\inf\bigl\{n\in\mathbb{N}|l_0(n)>0\bigr\}.
\]
It follows from (\ref{eq5.11}) that
\[
P(\sigma=k+1)\leq C(kT+1)^{-{q}/{2}}.
\]
Now for $1<q_1<\frac{q}{2}-1$,
%
%e5.13 #&#
%
\begin{equation}
\label{eq5.13} E\sigma^{q_1}1_{\sigma<\infty}\leq K_1,
\end{equation}
where $K_1$ is a constant.
For $\delta{:=}\min\{n\in\mathbb{N}|\|\theta_1(nT)\|
_{L^{2m(p-1)}}^{2m(p-1)}+\break\|\theta_2(nT)\|_{L^{2m(p-1)}}^{2m(p-1)}\leq
2C_0\}$, by Proposition~\ref{pro5.1.5}
we obtain that there exist $\gamma>0$ and
$c>0$ such that
%
%e5.14 #&#
%
\begin{equation}
\label{eq5.14} E\bigl(e^{\gamma\delta}\bigr)\leq c\bigl(1+\|
\theta_1^0\| _{L^{2m(p-1)}}^{2m(p-1)}+\|
\theta_2^0\|_{L^{2m(p-1)}}^{2m(p-1)}\bigr)
\end{equation}
(cf. \cite{M02}, \cite{O06}, (1.56)), where we used $C_0>4C_S^p[\frac
{1}{2}p(p-1)]^{p/2}\lambda_1^{-{p}/{2}}\mathcal{E}_0^{p/2}$. Set
\[
\cases{ %
\delta_0:=\delta, &
\cr
\sigma_{k+1}:=\infty&\quad$\mbox{if } \delta_k=\infty;
\hspace *{23.1pt} \sigma _{k+1}:=\sigma\circ\Phi_{\delta_kT}+
\delta_k \qquad\mbox{else}$,
\cr
\delta_{k}:=\infty&
\quad$\mbox{if } \sigma_k=\infty;\qquad\delta _{k}:=
\delta\circ \Phi_{\sigma_kT}+\sigma_k \hspace*{34pt}
\mbox{else}$, %\end{array}
}
\]
where $\Phi_t$ is the shift operator. Set $\eta:=\sigma+\delta\circ
\Phi
_{\sigma T}$.
If $l_0(0)=0$, by the Markov property, (\ref{eq5.13}) and (\ref{eq5.14})
\begin{eqnarray*}
E\bigl(\eta^{q_1}1_{\eta<\infty}\bigr)&\leq& C\bigl(E
\bigl(\sigma^{q_1}1_{\sigma
<\infty}\bigr)+E\bigl((\delta\circ
\Phi_{\sigma T})^{q_1}1_{\delta\circ\Phi
_\sigma
<\infty}1_{\sigma<\infty}\bigr)
\bigr)
\\
&\leq& C\bigl(E\bigl(\sigma^{q_1}1_{\sigma<\infty}\bigr)
\\
&&\phantom{\hspace*{10pt}}{}+cE\bigl(1+\bigl \|
\theta_1(\sigma T)\bigr \| _{L^{2m(p-1)}}^{2m(p-1)}+\bigl \|
\theta_2(\sigma T)\bigr \| _{L^{2m(p-1)}}^{2m(p-1)}
\bigr)1_{\sigma<\infty}\bigr)
\\
&\leq&C\bigl(1+\bigl \|\theta_1^0\bigr \| _{L^{2m(p-1)}}^{2m(p-1)}+
\bigl \|\theta_2^0\bigr \| _{L^{2m(p-1)}}^{2m(p-1)}\bigr),
\end{eqnarray*}
where we used Proposition~\ref{pro5.1.5} in the last inequality.
Since $\delta_k=\delta_{k-1}+\eta\circ\Phi_{\delta_{k-1}T}$, we obtain
for $1<q_1<\frac{q}{2}-1$,
%
%e5.15 #&#
%
\begin{eqnarray}\label{eq5.15}
E\bigl(\delta_k^{q_1}1_{\delta_k<\infty}
\bigr)&\leq& (k+1)^{q_1-1}\Biggl(E\delta^{q_1}+\sum
_{n=0}^{k-1}E(\eta\circ\Phi _{\delta_n
T})^{q_1}1_{\eta\circ\Phi_{\delta_n T}<\infty}
\Biggr)
\nonumber
\\[-8pt]
\\[-8pt]
&\leq&C(k+1)^{q_1} \bigl(1+\bigl \|\theta_1^0
\bigr \|_{L^{2m(p-1)}}^{2m(p-1)}+\bigl \|\theta_2^0\bigr \|
_{L^{2m(p-1)}}^{2m(p-1)}\bigr).
\nonumber
\end{eqnarray}
Moreover, if $\delta_k<\infty$, then $\sigma_{k+1}=\infty$ deduces that
$l_0(\infty)=\delta_k$. Define
\[
k_0:=\inf\bigl\{k\in\mathbb{Z}^+|\sigma_{k+1}=\infty\bigr
\}.
\]
Then (\ref{eq5.12}) implies that
%
%e5.16 #&#
%
\begin{equation}\label{eq5.16}
P(k_0\geq n)\leq(1-p_0)^n.
\end{equation}
By (\ref{eq5.16}), we obtain $k_0<\infty$ a.s., which implies
$l_0(\infty
)<\infty$ a.s. Moreover,
we have for $1<q_2<\frac{q}{2}-1$,
\[
E\bigl(l_0(\infty)^{q_2}\bigr)\leq\sum
_{n=0}^\infty E\bigl(\delta_n^{q_2}1_{\delta
_n<\infty}1_{k_0=n}
\bigr).
\]
Then by H\"{o}lder's inequality, we have for $\frac{1}{p_1}+\frac
{1}{p_1'}=1$, $p_1,p_1'>1$, satisfying $p_1q_2<\frac{q}{2}-1$\vspace*{1pt}
\[
E\bigl(l_0(\infty)^{q_2}\bigr)\leq\sum
_{n=0}^\infty\bigl(E\delta _n^{p_1{q_2}}1_{\delta_n<\infty}
\bigr)^{{1}/{p_1}}P(k_0=n)^{{1}/{p_1'}}.
\]
By (\ref{eq5.15}) and (\ref{eq5.16}),
we obtain\vspace*{1pt}
\begin{eqnarray*}
E\bigl(l_0(\infty)^{q_2}\bigr)&\leq& C\Biggl(\sum
_{n=0}^\infty(n+1)^{q_2}(1-p_0)^{
{n}/{p_1'}}
\Biggr) \bigl(1+\bigl \|\theta_1^0\bigr \|_{L^{2m(p-1)}}^{2m(p-1)}+
\bigl \|\theta_2^0\bigr \| _{L^{2m(p-1)}}^{2m(p-1)}\bigr)
\\
&<&
\infty.
\end{eqnarray*}

\textit{Step} 6. Conclusion.

By Step 2 and Step 5, we have for $t>0$ and $1<q_2<\frac{q}{2}-1$\vspace*{1pt}
\begin{eqnarray*}
&&E\bigl |\Lambda^{-{1}/{2}}\bigl(\theta_2\bigl(t,
\theta_0^1,\theta _0^2\bigr)-
\theta_1\bigl(t,\theta_0^1,
\theta_0^2\bigr)\bigr)\bigr |
\\
&&\qquad\leq E\bigl(\bigl |\Lambda^{-
{1}/{2}}\bigl(\theta_2\bigl(t,
\theta_0^1,\theta_0^2\bigr)-
\theta_1\bigl(t,\theta _0^1,\theta
_0^2\bigr)\bigr)\bigr | 1_{l_0(\infty)\leq l}\bigr)
\\
&&\qquad\quad{}+CP
\bigl(l_0(\infty)\geq l+1\bigr)^{1/2}
\\[1pt]
&&\qquad\leq C\bigl(1+\bigl \|\theta_1^0\bigr \|_{L^{2m(p-1)}}^{2m(p-1)}+
\bigl \|\theta_2^0\bigr \| _{L^{2m(p-1)}}^{2m(p-1)}\bigr)
\bigl[(t+1-lT)^{-q+1}+(l+1)^{-q_2/2}\bigr], %\end{array}
\end{eqnarray*}
where we used Proposition~\ref{pro5.1.5} in the first inequality.
Choosing $l=[\frac{t+1}{2T}]$, we obtain for $1<q_3<\frac{q}{4}-1$\vspace*{1pt}
%
%e5.17 #&#
%
\begin{eqnarray}\label{eq5.17}
&&E\bigl |\Lambda^{-{1}/{2}}\bigl(\theta_2\bigl(t,\theta
_0^1,\theta_0^2\bigr)-
\theta_1\bigl(t,\theta_0^1,
\theta_0^2\bigr)\bigr)\bigr |
\nonumber
\\[-8pt]
\\[-8pt]
&&\qquad\leq C\bigl(1+\bigl \|\theta _1^0
\bigr \|_{L^{2m(p-1)}}^{2m(p-1)}+\bigl \|\theta_2^0\bigr \|
_{L^{2m(p-1)}}^{2m(p-1)}\bigr) (t+1)^{-q_3}.
\nonumber
\end{eqnarray}
Thus, for $\psi\in C(H^1)$ with $C_\psi:=\sup_{x,y\in H^1}\frac
{|\psi
(x)-\psi(y)|}{|\Lambda^{-{1}/{2}}(x-y)|}<\infty$, we have\vspace*{1pt}
%
%e5.18 #&#
%
\begin{eqnarray}\label{eq5.18}
&&\bigl |P_t\psi(x)-P_t\psi(y)\bigr |
\nonumber
\\
&&\qquad\leq C_\psi E\bigl |\Lambda^{-{1}/{2}}\bigl(
\theta_2(t,x,y)-\theta_1(t,x,y)\bigr)\bigr |
\\
&&\qquad\leq C C_\psi\bigl(1+\|x\|_{L^{2m(p-1)}}^{2m(p-1)}+
\|y\| _{L^{2m(p-1)}}^{2m(p-1)}\bigr) (t+1)^{-q_3}.
\nonumber
\end{eqnarray}
By Proposition~\ref{pro5.1.5}, we obtain that for $2<p_2<\infty$\vspace*{1pt}
\begin{eqnarray*}
E\bigl \|\theta(t)\bigr \|_{L^{p_2}}^{p_2}&\leq&\|x\|_{L^{p_2}}^{p_2}e^{-\lambda
_1t}
\\[1pt]
&&{}+C_S^{p_2}
\biggl[\frac{1}{2}p_2(p_2-1)\biggr]^{p_2/2}
\lambda_1^{-
{p_2}/{2}}\mathcal{E}_0^{p_2/2}
\bigl(1-e^{-\lambda_1t}\bigr).
\end{eqnarray*}
Since for any invariant measure $\mu$ on $H^1$ and any $\varepsilon>0$,
there exists $b_\varepsilon>0$ such that $\mu(x\in H^1\dvtx\|x\|
_{L^{p_2}}^{p_2}> b_\varepsilon)\leq\varepsilon$,
we obtain that for any $L>0$
\begin{eqnarray*}
\int\bigl(\|x\|_{L^{p_2}}^{p_2}\wedge L\bigr)\,d\mu&
\leq&\int_{\{x:\|x\|
_{L^{p_2}}^{p_2}\leq b_\varepsilon\}}\bigl(E^x\bigl \|\theta(t)\bigr \|
_{L^{p_2}}^{p_2}\wedge L\bigr)\,d\mu+L\varepsilon
\\
&\leq& b_\varepsilon e^{-\lambda_1t}+C_S^{p_2}
\biggl[\frac{1}{2}p_2(p_2-1)\biggr]^{p_2/2}
\lambda _1^{-{p_2}/{2}}\mathcal{E}_0^{p_2/2}
\bigl(1-e^{-\lambda_1t}\bigr)
\\
&&{}+L\varepsilon . %\end{array}
\end{eqnarray*}
Letting $t\rightarrow\infty, \varepsilon\rightarrow0$ and
$L\rightarrow
\infty$,
we obtain that for any invariant measure~$\mu$
%
%e5.19 #&#
%
\begin{equation}
\label{eq5.19} \int\|x\|^{p_2}_{L^{p_2}}\,d\mu(x)\leq
C_S^{p_2}\biggl[\frac
{1}{2}p_2(p_2-1)
\biggr]^{p_2/2}\lambda_1^{-{p_2}/{2}}\mathcal{E}_0^{p_2/2}.
\end{equation}

Then by (\ref{eq5.18}), (\ref{eq5.19}) for any invariant measures
$\mu_1$, $\mu_2$ we
obtain for $\psi\in C(H^1)$ with $C_\psi<+\infty$ and $1<q_3<\frac{q}{4}-1$,
\begin{eqnarray*}
&&\biggl|\int\psi(x)\mu_1(dx)-\int\psi(x)\mu_2(dx)\biggr|
\\
&&\qquad\leq C C_\psi \biggl(1+\int\|x\|_{L^{2m(p-1)}}^{2m(p-1)}
\mu_1(dx)+\int\|x\| _{L^{2m(p-1)}}^{2m(p-1)}
\mu_2(dx)\biggr) (t+1)^{-q_3}. %\end{array}
\end{eqnarray*}
Letting $t\rightarrow\infty$, we get that $\mu_1=\mu_2$.
\end{pf}
%
%$\hfill\Box$

\begin{remark}\label{rem5.2.4}
(i) The coupling method has been introduced, for
example, in \cite{KS01,KPS02,KS02,M02,DO05} to study ergodicity for
stochastic partial differential equations. In these papers, they
decompose the process into the sum of a strongly dissipative process
$h$ and another finite dimensional dynamics $l$ driven by a
nondegenerate noise. The process is uniquely determined by the
nondegenerate part $l$ which can be treated by probabilistic arguments.
However, in our case, we cannot decompose the process into the two
desired parts since the uniqueness of the process $h$ depends on the
$L^p$-norm estimate, which cannot be obtained for $h$.

(ii) It is not clear how to directly use the results in \cite{O08} for the
following two reasons: Although we consider the semigroup in $H^1$, the
convergence we used in Theorem~\ref{the5.2.3} is in $H^{-{1}/{2}}$. In
\cite{O08}, only one state space has been considered. If we choose the
general Hilbert space in \cite{O08} as $H^1$, we cannot get the estimate
(\ref{eq5.5}) for the $H^1$-norm. If we choose the general Hilbert
space in
\cite{O08} as $H^{-1/2}$, the estimate (\ref{eq5.5}) does also not
hold for rough
initial values in $H^{-1/2}$. The second reason is that, since
Theorem~\ref{the5.2.3} depends on the $L^p$-norm estimate, we can only
prove $E\|\theta
_1(t,\theta_0^1,\theta_0^2)-\theta_2(t,\theta_0^1,\theta_0^2)\|
_{H^{-1/2}}$ converges to zero polynomially fast instead of
exponentially fast, when time goes to infinity, where $(\theta
_1(t,\theta_0^1,\theta_0^2),\theta_2(t,\theta_0^1,\theta_0^2))$ denotes
a coupling of two solutions to (\ref{eq3.1}) with different initial values
$\theta_0^i\in H^1$, $i=1,2$.

(iii) In the situation of Theorem~\ref{the5.2.3}, we also obtain that
$P_t$ on $
H^1$ is asymptotically strong Feller. In fact, for $x,y\in H^1$, define
$d_n(x,y):=1\wedge n|\Lambda^{-1/2}(x-y)|$. For any two probabilities
on $ H^1$ $\mu_1$, $\mu_2$, we denote the set of positive measures on $
H^1\times H^1$ with marginals $\mu_1$ and $\mu_2$ by $\mathcal
{C}(\mu
_1,\mu_2)$. Define the Wasserstein distance
\[
\|\mu_1-\mu_2\|_d:=\inf
_{\mu\in\mathcal{C}(\mu_1,\mu_2)}\int_{
H^1\times H^1}d(x,y)\mu(dx,dy).
\]
By definition and (\ref{eq5.17}), we obtain
\begin{eqnarray*}
\bigl \|P_n(x,\cdot)-P_n(y,\cdot)
\bigr \|_{d_n}&\leq&nE\bigl |\Lambda ^{-1/2}\bigl(\theta_2(n,x,y)-
\theta_1(n,x,y)\bigr)\bigr |
\\
&\leq&C\bigl(\|x\|_{L^{2m(p-1)}},\|y\|_{L^{2m(p-1)}}\bigr)n n^{-q_3}.
\end{eqnarray*}
Then we have
\[
\lim_{\gamma\rightarrow0}\limsup_{n\rightarrow\infty}\sup
_{y\in
B(x,\gamma)}\bigl \|P_n(x,\cdot)-P_n(y,\cdot)
\bigr \|_{d_n}=0,
\]
where $B(x,\gamma)$ denotes the ball in $ H^1$ with center $x$ and
radius $\gamma$, which implies that $P_t$ on $ H^1$ is asymptotically
strong Feller.

(iv) It seems difficult to directly verify the gradient estimate for
the semigroup as \cite{HM06} did for the 2D Navier--Stokes equation. By
their method, we need to consider an infinitesimal perturbation to the
initial condition and to estimate the derivative of the solution
$D\theta$ with respect to the initial value, which requires a good
estimate for $E\exp\|\theta\|_{L^p}^p$. However, this cannot be
obtained for $\alpha>\frac{1}{2}$. Even if the noise is nondegenerate
and we use the Bismut--Elworthy--Li formula to compute the gradient of
the semigroup, the ergodicity results only holds for $\alpha>\frac
{2}{3}$ by delicate estimates (see Section~\ref{sec6}). We cannot
directly use
the criterion in \cite{KPS10}, since it is not clear how to verify the
e-property in \cite{KPS10} for the semigroup associated with the 2D
stochastic quasi-geostrophic equation.
\end{remark}

%s5.3 #&#
\subsection{Existence of invariant measures for \texorpdfstring{$\alpha>\frac{1}{2}$}{$alpha>\frac{1}{2}$}}\label{sec5.3}

Assume that $G$ satisfies condition Hypothesis~\ref{hypE.1}.

\begin{Lemma}\label{lem5.3.1}
Let $\alpha>\frac{1}{2}$. If $\theta_0\in H^1, t>0$, then:
\begin{enumerate}[(iii)]
\item[(i)]$ E(|\theta(t)|^2)+E\int_0^t|\Lambda^\alpha\theta
(r)|^2\,dr\leq
|\theta_0|^2+t\operatorname{Tr}[GG^*]$,

\item[(ii)]for $\delta\leq1$ and $q\geq\frac{2\alpha+2}{2\alpha
-1}, p\geq
1$, we have
\[
E\int_0^t\frac{|\Lambda^{\delta+\alpha}\theta(r)|^2}{(1+|\Lambda
^\delta
\theta(r)|^2)^{p+1}}\,dr\leq C\biggl(
\int_0^t E\bigl \|\theta(r)\bigr \|_{L^q}^q
\, dr+1\biggr)\leq Ct\bigl(\|\theta_0\|_{L^q}^q+1
\bigr),
\]

\item[(iii)]for $q\geq\frac{2\alpha+2}{2\alpha-1}$, there exist
$0<\delta
_1<1-\alpha$ and $0<\gamma_0<1$ such that
\[
E\biggl[\int_0^t\bigl |A_\alpha^{\delta_1}
\theta(r)\bigr |_{H^1}^{2\gamma_0} \, dr\biggr]\leq C(1+t) \bigl(\|
\theta_0\|_{L^q}^q+1\bigr).
\]
\end{enumerate}
\end{Lemma}

\begin{pf}
(i) is well known and follows from It\^{o}'s formula applied\vspace*{1pt} to
$|\theta(t)|^2$. By Theorems~\ref{theA.1}, \ref{theA.2} in
Appendix~\hyperref[appA]{A}, we obtain
$\theta\in C([0,\infty),H^1)\cap L^2_{\mathrm{loc}}([0,\infty
),\allowbreak  H^{1+\alpha
})$ $P$-a.s. By a similar argument as in the proof of Theorem~\ref
{the4.2}, we
obtain for $\delta\leq1$
\begin{eqnarray*}
&&\tfrac{1}{2}d\bigl |\Lambda^\delta\theta\bigr |^2+
\kappa\bigl |\Lambda^{\delta
+\alpha}\theta\bigr |^2\,dt+\bigl\langle
\Lambda^{\delta-\alpha}(u\cdot \nabla\theta), \Lambda^{\delta+\alpha}\theta\bigr
\rangle\,dt
\\
&&\qquad=\bigl\langle\Lambda^{\delta
}\theta,\Lambda^\delta G
\,dW_t\bigr\rangle+\tfrac{1}{2}\operatorname{Tr}\bigl[GG^*
\Lambda ^{2\delta}\bigr]\,dt. %\end{array}
\end{eqnarray*}
Then we apply It\^{o}'s formula to the function $(1+|\Lambda^\delta
\theta|^2)^{-p}$ and get
\begin{eqnarray*}
&&\frac{1}{(1+|\Lambda^\delta\theta(t)|^2)^{p}}-\frac
{1}{(1+|\Lambda^\delta\theta_0|^2)^{p}}
\\
&&\qquad= 2p\kappa\int_0^t\frac{|\Lambda^{\delta+\alpha}\theta
|^2}{(1+|\Lambda
^\delta\theta|^2)^{p+1}}
\,dr+2p\int_0^t\frac{\langle\Lambda
^{\delta
-\alpha}(u\cdot\nabla\theta),\Lambda^{\delta+\alpha}\theta
\rangle
}{(1+|\Lambda^\delta\theta|^2)^{p+1}}\,dr
\\
&&\phantom{\qquad=} {}-2p\int_0^t
\frac{\langle
\Lambda^{\delta}\theta,\Lambda^\delta G\,dW_r\rangle}{(1+|\Lambda
^\delta
\theta|^2)^{p+1}}-p\int_0^t\frac{\operatorname{Tr}[GG^*\Lambda
^{2\delta
}]}{(1+|\Lambda^\delta\theta|^2)^{p+1}}
\,dr
\\
&&\phantom{\qquad=} {}+2p(p+1)\int_0^t
\frac
{|G^*\Lambda^{2\delta}\theta|^2}{(1+|\Lambda^\delta\theta
|^2)^{p+2}}\,dr, %\end{array}
\end{eqnarray*}
where the last term is meaningful since $|G^*\Lambda^{2\delta}\theta
|^2\leq|\Lambda^\delta\theta|^2\|\Lambda^\delta G\|_{L_2(H,H)}^2$.
For $q\geq\frac{2\alpha+2}{2\alpha-1}$ and $\sigma:=\frac
{2}{q}<2\alpha
-1$, we have
\begin{eqnarray*}
\bigl |\bigl\langle\Lambda^{\delta-\alpha}(u\cdot\nabla\theta),\Lambda
^{\delta+\alpha}\theta\bigr\rangle\bigr |&=&\bigl |\bigl\langle\Lambda^{\delta-\alpha
}\nabla
\cdot(u\theta),\Lambda^{\delta
+\alpha}\theta\bigr\rangle\bigr |
\\
&\leq& C\bigl |\Lambda^{\delta-\alpha+1+\sigma}\theta\bigr |\cdot\|\theta\| _{L^q} \bigl |
\Lambda^{\delta+\alpha}\theta\bigr |
\\
&\leq& C\|\theta\|_{L^q}^{
{2\alpha}/({2\alpha-1-\sigma})}\bigl |\Lambda^\delta
\theta\bigr |^2+\kappa \bigl |\Lambda ^{\delta+\alpha}\theta\bigr |^2,
\end{eqnarray*}
where we used $\operatorname{div}u=0$ in the first equality and
Lemmas~\ref{lem2.1} and~\ref{lem2.2}
in
the first inequality and Young's together with the interpolation
inequality (\ref{eq2.3}) in the last inequality.

Hence, we obtain
\[
E\int_0^t\frac{|\Lambda^{\delta+\alpha}\theta|^2}{(1+|\Lambda
^\delta
\theta|^2)^{p+1}}\,dr\leq C\biggl(
\int_0^t E\|\theta\|_{L^q}^q
\,dr+t\biggr)\leq Ct\bigl(\| \theta_0\|_{L^q}^q+1
\bigr),
\]
where we used Proposition~\ref{pro5.1.5} in the last step.

(iii) Since by Young's inequality for some $\gamma_0>0$, we have
\[
\bigl |\Lambda^{\delta+\alpha}\theta\bigr |^{2\gamma_0}\leq c\biggl[
\frac{|\Lambda
^{\delta
+\alpha}\theta|^2}{(1+|\Lambda^\delta\theta|^2)^{p+1}}+1+\bigl |\Lambda ^\delta\theta\bigr |^2\biggr],
\]
we obtain for $\delta+\alpha>1$
\[
E\biggl[\int_0^t\bigl |
\Lambda^{\delta+\alpha} \theta\bigr |^{2\gamma_0} \,dr\biggr]\leq C(1+t) \bigl(\|
\theta_0\|_{L^q}^q+1\bigr). %\end{array}
\]
\upqed
\end{pf}
%
%$\hfill\Box$

\begin{Theorem}\label{the5.3.2}
Let $\alpha>\frac{1}{2}$ and suppose Hypothesis~\ref{hypE.1} holds.
Then $(P_t)_{t\geq0}$ is $ H^1$-Feller, that is, for every $t>0$ and
$\psi\in\mathcal{C}_b( H^1)$, $P_t\psi\in C_b( H^1)$. Furthermore,
there exists an invariant measure $\nu$ on $ H^1$ of the transition
semigroup $(P_t)_{t\geq0}$. Moreover, there are $0<\delta_1<1-\alpha$
and $0<\gamma_0<1$ such that
\[
\int\bigl |A_\alpha^{\delta_1} x\bigr |^{2\gamma_0}_{H^1}\,d\nu<
\infty.
\]
\end{Theorem}

\begin{pf}
Choose $x_0\in H^1$ and define for $t>0$
\[
\mu_t:=\frac{1}{t}\int_0^tP_r^*
\delta_{x_0} \,dr.
\]
By Lemma~\ref{lem5.3.1}(iii), we have for $t>1$ that
\[
\int\bigl |A_\alpha^{\delta_1} x\bigr |^{2\gamma_0}_{H^1}
\mu_t(dx)\leq C.
\]
This implies that $\{\mu_t|t>0\}$ is tight on $ H^1$.
By Theorem~\ref{theA.3} in Appendix \hyperref[appA]{A}, we obtain that $(P_t)_{t\geq
0}$ is $
H^1$-Feller.
Hence, any limit point of $\mu_t$ is an invariant measure for
$(P_t)_{t\geq0}$.
\end{pf}
%
%$\hfill\Box$

Combining Theorem~\ref{the5.2.3} and Theorem~\ref{the5.3.2}, we
obtain the following results.

\begin{Theorem}\label{the5.3.3}
Fix $\alpha>1/2$. Assume Hypotheses~\ref{hypE.1} and~\ref{hypE.2}
hold.\vspace*{1pt} Let
$\delta_0=\lambda_{N+1}-2^{p/2}C_R^p C_S^{p+1}\kappa^{1-p}[p(p-1)]
^{p/2}\lambda_1^{-p/2}\mathcal{E}_0^{{p} /2}>0$ for $p=\frac{\alpha
+1}{\alpha-{1}/{2}}$, where $N$ is as in Hypothesis~\ref{hypE.2}, $C_S,
C_R$ are the constants for Sobolev embedding and Riesz transform,
respectively. Then there exists exactly one invariant probability
measure $\nu$ for $P_t$.

Moreover, for $\psi\in C(H^1)$ with $C_\psi:=\sup_{x,y\in H^1}\frac
{|\psi(x)-\psi(y)|}{|\Lambda^{-{1}/{2}}(x-y)|}<\infty$ and any
initial distribution $\mu_0$ on $H^1$ with $\int\|x\|
_{L^{2m(p-1)}}^{2m(p-1)}\,d\mu_0<\infty$ for some $m>35$, the following
polynomial bound is satisfied for $1<q_3<\frac{m-19}{16}$:
%
%e5.20 #&#
%
\begin{eqnarray}\label{eq5.20}
&&\biggl|\int P_t\psi(x)\mu_0(dx)-\int\psi(x)\nu
(dx)\biggr|
\nonumber
\\[-8pt]
\\[-8pt]
&&\qquad\leq C C_\psi\biggl(1+\int\|x\|_{L^{2m(p-1)}}^{2m(p-1)}
\mu _0(dx)\biggr) (t+1)^{-q_3}.
\nonumber
\end{eqnarray}
\end{Theorem}

\begin{pf}
(\ref{eq5.20}) can be easily deduced from (\ref{eq5.18}) and (\ref{eq5.19}).
\end{pf}
%
%$\hfill\Box$

%s5.4 #&#
\subsection{Law of large numbers}\label{sec5.4}

In this section, we establish the law of large numbers for the solution
of the stochastic quasi-geostrophic equation. The proof is mainly
inspired by the approach used in \cite{KW}.

\begin{Theorem}\label{the5.4.1}
Fix $\alpha>1/2$. Assume Hypotheses~\ref{hypE.1} and~\ref{hypE.2}
hold. Set
$\delta_0:=\lambda_{N+1}-2^{p/2}C_R^p C_S^{2p}\kappa^{1-p}[p(p-1)]
^{p/2}\lambda_1^{-p/2}\mathcal{E}_0^{{p} /2}>0$ for $p=\frac{\alpha
+1}{\alpha-{1}/{2}}$, where $N$ is as in Hypothesis~\ref{hypE.2}, $C_S,
C_R$ are the constants for the Sobolev embedding and Riesz transform,
respectively.
Then for $\psi\in C(H^1)$ with $C_\psi:=\sup_{x,y\in H^1}\frac
{|\psi
(x)-\psi(y)|}{|\Lambda^{-{1}/{2}}(x-y)|}<\infty$ and any initial
distribution $\mu_0$ on $H^1$ with\break $\int\|x\|
_{L^{2m(p-1)}}^{2m(p-1)}\,d\mu_0<\infty$ for some $m>35$,
\[
\lim_{T\rightarrow\infty}\frac{1}{T}\int_0^T
\psi\bigl(\theta(s)\bigr)\, ds=\int\psi \,d\nu\qquad\mbox{in probability}.
\]
\end{Theorem}

\begin{pf}
(\ref{eq5.20}) implies that for $\psi\in C(H^1)$ with $C_\psi<\infty$
%
%e5.21 #&#
%
\begin{equation}
\label{eq5.21} %\begin{array}{cc}
\lim_{T\rightarrow\infty}\biggl|\frac{1}{T}\int
_0^TE\psi\bigl(\theta(t)\bigr)\,dt-\int\psi(x)
\nu(dx)\biggr|=0. %\end{array}
\end{equation}
Now we want to prove that for bounded $\psi\in C(H^1)$ with $C_\psi
<\infty$
%
%e5.22 #&#
%
\begin{equation}\label{eq5.22}
\lim_{T\rightarrow\infty}\biggl|\frac{1}{T^2}E\biggl(\int
_0^T\psi\bigl(\theta(t)\bigr)\,dt
\biggr)^2-\biggl(\int\psi(x)\nu(dx)\biggr)^2\biggr|=0.
\end{equation}

We have
\begin{eqnarray*}
\frac{1}{T^2}E\biggl(\int_0^T
\psi\bigl(\theta(t)\bigr)\,dt\biggr)^2&=&\frac
{1}{T^2}E\biggl(
\int_0^T\psi\bigl(\theta(t)\bigr)\,dt\int
_0^T\psi\bigl(\theta(s)\bigr)\,ds\biggr)
\\
&=&\frac{2}{T^2}\int_0^T\int
_0^tE\bigl[\psi\bigl(\theta(t)\bigr)\psi\bigl(
\theta (s)\bigr)\bigr]\,dt\,ds
\\
&=&\frac{2}{T^2}\int_0^T\int
_0^t\langle\mu_0P_s,
\psi P_{t-s}\psi \rangle\,dt\,ds. %\end{array}
\end{eqnarray*}
Moreover, we have that for $B:=\{\|x\|_{L^{2m(p-1)}}\leq R\}$,
\begin{eqnarray*}
&&\biggl|\frac{2}{T^2}\int_0^T\int
_0^t\biggl\langle\mu_0P_s,
\psi\biggl(P_{t-s}\psi -\int\psi(x)\nu(dx)\biggr)\biggr\rangle\,dt\,ds\biggr|
\\
&&\qquad\leq\biggl|\frac{2}{T^2}\int_0^T\int
_0^t\biggl\langle\mu _0P_s,1_{B}
\psi \biggl(P_{t-s}\psi-\int\psi(x)\nu(dx)\biggr)\biggr\rangle\,dt\,ds\biggr|
\\
&&\phantom{\qquad\leq} {}+\biggl|\frac{2}{T^2}\int_0^T
\int_0^t\biggl\langle\mu_0P_s,1_{B^c}
\psi\biggl(P_{t-s}\psi-\int\psi (x)\nu (dx)\biggr)\biggr\rangle\,dt\,ds\biggr|
\\
&&\qquad:=I_T+\mathit{II}_T. %\end{array}
\end{eqnarray*}
By (\ref{eq5.20}), we obtain that there exists $T_1>0$ such that for
any $T>T_1$
%
%e5.23 #&#
%
\begin{equation}
\label{eq5.23} \sup_{x\in B}\biggl|\frac{1}{T}\int
_0^TP_t\psi(x)-\int\psi \,d\nu\biggr|<
\varepsilon.
\end{equation}
Thus, for the first term we have the following:
\begin{eqnarray*}
I_T&=&\biggl|\frac{2}{T^2}\int_0^T(T-s)
\biggl\langle\mu_0P_s,1_{B}\psi \biggl[
\frac{1}{T-s}\int_0^{T-s}
\biggl(P_{t}\psi-\int\psi(x)\nu(dx)\biggr)\, dt\biggr]\biggr\rangle
\,ds\biggr|
\\
&\leq&\biggl|\frac{2}{T^2}\int_0^{T_1}s\biggl
\langle\mu_0P_{T-s},1_{B}\psi \biggl[
\frac
{1}{s}\int_0^{s}
\biggl(P_{t}\psi-\int\psi(x)\nu(dx)\biggr)\,dt\biggr]\biggr\rangle
\,ds\biggr|
\\
&&{}+\biggl|\frac
{2}{T^2}\int_{T_1}^Ts\biggl
\langle\mu_0P_{T-s},1_{B}\psi\biggl[
\frac{1}{s}\int_0^{s}
\biggl(P_{t}\psi-\int\psi(x)\nu(dx)\biggr)\,dt\biggr]\biggr\rangle
\,ds\biggr|
\\
&\leq& 4\|\psi\|_{L^\infty}^2\biggl(\frac{T_1}{T}
\biggr)^2+\varepsilon\|\psi\| _{L^\infty}, %\end{array}
\end{eqnarray*}
where we used (\ref{eq5.23}) in the last step.
For the second term by Proposition~\ref{pro5.1.5}, we have
\begin{eqnarray*}
\mathit{II}_T&\leq& \frac{4\|\psi\|_{L^\infty}^2}{T^2}\int
_0^T\int_0^t
\mu _0P_s\bigl(B^c\bigr)\,ds\,dt
\\
&\leq&\|\psi\|_{L^\infty}^2\frac{C}{R}. %\end{array}
\end{eqnarray*}
Choosing $R$ large enough, we obtain for any $\varepsilon>0$ that there
exists $T_0$ such that for $T\geq T_0$
\[
\biggl|\frac{2}{T^2}\int_0^T\int
_0^t\biggl\langle\mu_0P_s,
\psi\biggl(P_{t-s}\psi -\int\psi (x)\nu(dx)\biggr)\biggr\rangle\,dt\,ds\biggr|
\leq\varepsilon.
\]
The latter implies
\begin{eqnarray*}
&&\lim_{T\rightarrow\infty}\biggl|\frac{1}{T^2}E\biggl(\int
_0^T\psi\bigl(\theta (t)\bigr)\,dt
\biggr)^2-\biggl(\int\psi(x)\nu(dx)\biggr)^2\biggr|
\\
&&\qquad\leq\lim_{T\rightarrow\infty}\biggl|\frac{2}{T^2}\int\psi (x)\nu(dx)
\int_0^T\int_0^t
\langle\mu_0P_s,\psi\rangle\,dt\,ds-\biggl(\int\psi(x)
\nu(dx)\biggr)^2\biggr|
\\
&&\qquad=\biggl|\int\psi(x)\nu(dx)\biggr|\lim_{T\rightarrow\infty}\biggl|\frac
{2}{T^2}
\int_0^Tt\,dt\biggl[\frac{1}{t}\int
_0^t\langle\mu_0P_s,
\psi\rangle\,ds-\int \psi(x)\nu (dx)\biggr]\biggr|
\\
&&\qquad=0. %\end{array}
\end{eqnarray*}

Now by (\ref{eq5.21}) and (\ref{eq5.22}) we obtain for bounded $\psi
$ with $C_\psi
<\infty$,
\[
\lim_{T\rightarrow\infty}\frac{1}{T}\int_0^T
\psi\bigl(\theta(s)\bigr)\, ds=\int\psi \,d\nu\qquad\mbox{in probability}.
\]

In general, we can remove the restriction of the boundedness of $\psi$
by defining $\psi_L=\psi\wedge L\vee(-L)$ for $L\in\mathbb{R}^+$. Since
for $x,y\in H^1$
\[
\bigl |\psi_L(x)-\psi_L(y)\bigr |\leq\bigl |\psi(x)-\psi(y)\bigr |\leq
C_\psi\bigl |\Lambda^{-
{1}/{2}}(x-y)\bigr |,
\]
we have
%
%e5.24 #&#
%
\begin{equation}
\label{eq5.24} \lim_{T\rightarrow\infty}\frac{1}{T}\int
_0^T\psi _L\bigl(\theta(s)\bigr)
\,ds=\int\psi_L \,d\nu\qquad\mbox{in probability}.
\end{equation}
Since $\int|\psi|\,d\nu<\infty$, it is clear that
\[
\lim_{L\rightarrow\infty}\int\psi_L \,d\nu=\int\psi\,d\nu.
\]
Applying (\ref{eq5.21}) for $|\psi-\psi_L|$, we have
\[
\lim_{L\rightarrow\infty}\lim_{T\rightarrow\infty}E\frac
{1}{T}\int
_0^T\bigl |\psi_L\bigl(\theta(s)\bigr)-
\psi\bigl(\theta(s)\bigr)\bigr |\,ds=\lim_{L\rightarrow
\infty}\int |
\psi_L -\psi|\,d\nu=0.
\]
Now the result follows by taking the limit on both sides of (\ref{eq5.24}).
\end{pf}
%
%$\hfill\Box$

%s6 #&#
\section{Exponential convergence for \texorpdfstring{$\alpha>\frac{2}{3}$}{$alpha>\frac{2}{3}$}}\label{sec6}

Under the conditions (Hypotheses~\ref{hypE.1}, \ref{hypE.2}) on $G$
we only obtain the semigroup
converges to the invariant measure polynomially fast [see (\ref
{eq5.20})]. In
this section, we prove that the convergence is exponentially fast,
however, under stronger conditions for $\alpha$ and $G$.
We assume that $\alpha>\frac{2}{3}$, and that $G$ satisfies:

\renewcommand{\theHypothesis}{E.3}
\begin{Hypothesis}\label{hypE.3}
There are an isomophism $Q_0$ of $H$ and a
number $ s\geq1$ such that $G=A_\alpha^{-({s+\alpha})/({2\alpha
})}Q_0^{1/2}$, and furthermore, $G$ satisfies \textup{(Gp.1)} for some fixed
$p\in
((\alpha-\frac{1}{2})^{-1},\infty)$ (which is, e.g., always the case if
$Q_0=I$).
\end{Hypothesis}

For $x:=\theta_0\in L^p$,
let $P_x$ denote the law of the corresponding solution $\theta(\cdot
,x)$ to (\ref{eq3.1}). Since Hypothesis~\ref{hypG.1},
(Gp.1) and (\ref{eqGL.1}) are satisfied under Hypothesis~\ref{hypE.3},
by Theorems~\ref{the4.3} and~\ref{the4.4}
the measures $P_x$, $x\in L^p$, form a Markov process.
Let $(P_t)_{t\geq0}$ be the associated transition semigroup on
$\mathcal{B}_b(H)$, defined as
%
%e6.1 #&#
%
\begin{equation}
\label{eq6.1} P_t(\varphi) (x):=E\bigl[\varphi\bigl(\theta(t,x)
\bigr)\bigr],\qquad x\in L^p, \varphi\in\mathcal{B}_b(H).
\end{equation}

\begin{Remark}\label{rem6.1}
If Hypothesis~\ref{hypE.3} is satisfied with $s>3-2\alpha$,
then Hypotheses~\ref{hypE.1}, \ref{hypE.2}
hold for $G$ and (Gp.1) holds for any $p\in(0,\infty)$.
\end{Remark}

%s6.1 #&#
\subsection{The strong Feller property for \texorpdfstring{$\alpha>\frac{2}{3}$}{$alpha>\frac{2}{3}$}}

In this subsection, we prove that its transition semigroup has the
strong Feller property under Hypothesis~\ref{hypE.3}.

\begin{Remark}\label{rem6.1.1}
(i) Since in our case $\alpha<1$, the linear part
$(-\Delta)^\alpha$ in (\ref{eq1.1}) is less regularizing. As
$G=A_\alpha^{-
({s+\alpha})/({2\alpha})}Q_0^{1/2}$,
we get the trajectories $z$ of the associated O--U process to be in $
C([0,\infty), H^{s+2\alpha-1-\varepsilon_0})$ for every $\varepsilon
_0>0$ (cf. \cite{DZ92}, Theorem~5.16, \cite{DO06}, Proposition~3.1). However,
in order to prove the weak-strong uniqueness principle
(see Theorem~\ref{the6.1.3} below) and the strong Feller property
of the semigroup associated
with the solution of the cutoff equation (see Proposition~\ref{pro6.1.4}
below),
we need $z\in C([0,\infty), H^{s+1-\alpha+\sigma_1})$ for some
$\sigma
_1>0$. Therefore, we need $s+2\alpha-1>s+1-\alpha$, that is, $\alpha
>\frac{2}{3}$. The situation of the 3D Navier--Stokes equation is
different. While in our case the needed regularity of $z$ is higher
than the regularity of our solution space
$C((0,\infty),H^s)$ for the cutoff equation (\ref{eq6.2}), for the
3D Navier--Stokes equation
the needed regularity of $z$ is the same as for the solution of the
cutoff equation.

(ii) Since $\alpha<1$, we cannot apply the same type of estimate as in
\cite{FR08} (cf. \cite{FR08}, Lemma D.2). Instead, we use Lemma~\ref
{lem2.1} and choose
suitable parameters ($s,\sigma_1,\sigma_2$) such that the approach in
\cite{FR08} can be modified to apply here [see (\ref{eq6.6})--(\ref
{eq6.10}) and so on].

(iii) It seems difficult to use the Kolmogorov equation method as in
\cite{DD03,DO06} or a coupling approach as in \cite{O07} in our
situation. In
fact, to get a uniform $H^s$-norm estimate for the solutions of the
Galerkin approximations of equation (\ref{eq1.1}) for some $s>0$,
the regularity, needed for the trajectories of the associated
Ornstein--Uhlenbeck (O--U)
process $z$ is higher than $H^s$, which is entirely different from the
situation of the 3D Navier--Stokes equation. According to the method
in \cite{DD03,DO06} and \cite{O07}, we should use the solutions'
$H^{s+\alpha
}$-norm to control the $H^{s+\alpha}$-norm of the derivative of the
solutions as required for the Bismut--Elworthy--Li formula. In
particular, the associated O--U process $z$ should be also in
$H^{s+\alpha}$. However, under Hypothesis~\ref{hypE.3} for the noise,
the O--U
process $z$ is only in $L^2([0,T],H^{s+2\alpha-1})$. As a consequence,
for their method to apply here, we need even $\alpha\geq1$.
\end{Remark}

Fix $s>1$ as in Hypothesis~\ref{hypE.3} and set $\mathcal{W}:=H^s$ and
$|x|_\mathcal{W}:=\|x\|_{H^s}$. In this subsection, we choose
\[
\Omega:=C\bigl([0,\infty);H^{-\beta}\bigr)
\]
for some $\beta>3$
and let $\mathcal{B}$ denote the Borel $\sigma$-algebra on $\Omega$.

Now we state the main result of this section.

\begin{Theorem}\label{the6.1.2}
Fix $\alpha>\frac{2}{3}$. Under Hypothesis~\ref{hypE.3},
$(P_t)_{t\geq0}$ is $\mathcal{W}$-strong Feller, that is, for every
$t>0$ and $\psi\in\mathcal{B}_b(H)$, $P_t\psi\in C_b(\mathcal{W})$.
\end{Theorem}

We shall use \cite{FR08}, Theorem~5.4, which is an abstract result to prove
the strong Feller property. In order to use \cite{FR08}, Theorem~5.4, we
follow the idea of \cite{FR08}, Theorem~5.11, to construct $P_x^{(R)}$.
We introduce an equation which differs from the original one by a
cut-off only, so that with large probability they have the same
trajectories on a small random time interval [see (\ref{eq6.3})
below]. We
consider the equation
%
%e6.2 #&#
%
\begin{equation}
\label{eq6.2} d\theta(t)+A_\alpha\theta(t)\,dt+\chi_{R}
\bigl(|\theta |_\mathcal{W}^2\bigr)u(t)\cdot\nabla\theta(t)
\,dt= G\,dW(t),
\end{equation}
where $\chi_R\dvtx\mathbb{R}\rightarrow[0,1]$ is of class $C^\infty
$ such
that $\chi_R(|\theta|)=1$ if $|\theta|\leq R$, $\chi_R(|\theta
|)=0$ if
$|\theta|>R+1$ and with its
first derivative bounded by $1$. Then, if we can prove the following
Theorem~\ref{the6.1.3} and Proposition~\ref{pro6.1.4}, Theorem~\ref
{the6.1.2} follows.

\begin{Theorem}[(Weak--strong uniqueness)]\label{the6.1.3}
Fix $\alpha>\frac{2}{3}$.
Suppose Hypothesis~\ref{hypE.3} holds. Then for every $x\in\mathcal{W}$,
equation (\ref{eq6.2}) has a unique martingale solution $P_x^{(R)}$, with
\[
P_x^{(R)}\bigl[C\bigl([0,\infty);\mathcal{W}\bigr)\bigr]=1.
\]
Let $\tau_R\dvtx\Omega\rightarrow[0,\infty]$ be defined by
\[
\tau_R(\omega):=\inf\bigl\{t\geq0\dvtx\bigl |\omega(t)\bigr |^2_\mathcal{W}
\geq R\bigr\},
\]
and $\tau_R(\omega):=\infty$ if this set is empty. If $x\in\mathcal
{W}$ and $|x|_{\mathcal{W}}^2<R$, then
%
%e6.3 #&#
%
\begin{equation}
\label{eq6.3} \lim_{\varepsilon\rightarrow0}P_{x+h}^{(R)}[
\tau_R\geq \varepsilon]=1, \qquad\mbox{uniformly in } h\in
\mathcal{W}, |h|_\mathcal{W}<1.
\end{equation}
Moreover,
%
%e6.4 #&#
%
\begin{equation}
\label{eq6.4} E^{P_x^{(R)}}\bigl[\varphi(\omega_t)1_{[\tau_R\geq
t]}
\bigr]=E^{P_x}\bigl[\varphi(\omega_t)1_{[\tau_R\geq t]}\bigr]
\end{equation}
for every $t\geq0$ and $\varphi\in\mathcal{B}_b(H)$.
\end{Theorem}

\begin{pf}
Let $z$ denote the solution to
\[
dz(t)+A_\alpha z(t)\,dt=G\,dW(t),
\]
with initial data $z(0)=0$ and let $v^{(R)}_x$ be the solution to the
auxiliary problem
%
%e6.5 #&#
%
\begin{equation}
\label{eq6.5} \qquad\frac{dv^{(R)}(t)}{dt}+A_\alpha v^{(R)}(t)+u^{(R)}(t)
\cdot\nabla\bigl(v^{(R)}(t)+z(t)\bigr)\chi _R
\bigl(\bigl |v^{(R)}+z\bigr |^2_\mathcal{W}\bigr)=0,
\end{equation}
with $v^{(R)}(0)=x$. Here, $u^{(R)}(t)=u_{v^{(R)}}(t)+u_z(t)$,
$u_{v^{(R)}}$ and $u_z$ satisfy (\ref{eq1.3}) with $\theta$ replaced by
$v^{(R)}$ and $z$, respectively. Moreover, define $\theta
^{(R)}:=v^{(R)}+z$, which is a weak solution to equation (\ref
{eq6.2}). We
denote its law on $\Omega$ by $P_x^{(R)}$.
By Hypothesis~\ref{hypE.3}, the trajectories of the noise belong to
\[
\Omega^*:=\bigcap_{\beta\in(0,{1}/{2}),\eta\in[0, ({s+\alpha
})/({2\alpha})-{1}/({2\alpha}))} C^\beta\bigl([0,
\infty);D\bigl(A_\alpha^\eta\bigr)\bigr),
\]
with probability one. Hence, the analyticity of the semigroup generated
by $A_\alpha$ implies that for each $\omega\in\Omega^*$, $z(\omega
)\in
C([0,\infty), H^{s+2\alpha-1-\varepsilon_0})$ for every $\varepsilon_0>0$.

Now, for $\omega\in\Omega^*$ we prove that equation (\ref{eq6.5}) with
$z(\omega)$ replacing $z$ has a unique global weak solution in the
space $C([0,\infty);\mathcal{W})$. First, we obtain the following a
priori estimate for suitable $\sigma_1,\sigma_2>0$ with $\sigma
_2\leq
s, \sigma_2+\sigma_1=1,s+\sigma_1-\alpha+1<s+2\alpha-1<s+\alpha$, where
we used that $\alpha>\frac{2}{3}$ since $0<\sigma_1<3\alpha-2$:
%
%e6.6 #&#
%
\begin{eqnarray}\label{eq6.6}
&&\frac{1}{2}\frac{d}{dt}\bigl |\Lambda^s
v^{(R)}\bigr |^2+\kappa\bigl |\Lambda^{s+\alpha}v^{(R)}\bigr |^2\nonumber
\\
&&\qquad=
\chi_R\bigl(\bigl |\theta ^{(R)}\bigr |^2_\mathcal{W}
\bigr)\bigl\langle\Lambda^{s-\alpha}\nabla\cdot\bigl(u^{(R)}
\theta^{(R)}\bigr),\Lambda^{s+\alpha}v^{(R)}\bigr\rangle
\nonumber
\\
&&\qquad\leq C\chi_R\bigl(\bigl |\theta^{(R)}\bigr |^2_\mathcal{W}
\bigr)\bigl |\Lambda^{s-\alpha+1}\bigl(u^{(R)} \theta^{(R)}
\bigr)\bigr |\cdot\bigl |\Lambda^{s+\alpha}v^{(R)}\bigr |
\nonumber
\\
&&\qquad\leq C\chi_R\bigl(\bigl |\theta^{(R)}\bigr |^2_\mathcal{W}
\bigr)\bigl |\Lambda^{s-\alpha+1+\sigma
_1}\theta^{(R)}\bigr |\bigl |\Lambda^{\sigma_2}
\theta^{(R)}\bigr |\cdot\bigl |\Lambda ^{s+\alpha
}v^{(R)}\bigr |
\nonumber
\\
&&\qquad\leq C\chi_R\bigl(\bigl |\theta^{(R)}\bigr |^2_\mathcal{W}
\bigr) \bigl(\bigl |\Lambda^{s-\alpha+1+\sigma
_1}v^{(R)}\bigr |+\bigl |\Lambda^{s-\alpha+1+\sigma_1}z\bigr |
\bigr)\cdot\bigl |\Lambda ^{s+\alpha}v^{(R)}\bigr |
\\
&&\qquad\leq C\chi_R\bigl(\bigl |\theta^{(R)}\bigr |^2_\mathcal{W}
\bigr) \bigl(C\bigl |\Lambda ^{s}v^{(R)}\bigr |^{1-r_1}\bigl |
\Lambda^{s+\alpha}v^{(R)}\bigr |^{r_1}+\bigl |\Lambda
^{s-\alpha
+1+\sigma_1}z\bigr |\bigr)\nonumber
\\
&&\qquad\quad{}\cdot\bigl |\Lambda^{s+\alpha}v^{(R)}\bigr |
\nonumber
\\
&&\qquad\leq C\chi_R\bigl(\bigl |\theta^{(R)}\bigr |^2_\mathcal{W}
\bigr) \bigl(\bigl |\Lambda ^{s}v^{(R)}\bigr |^2+\bigl |
\Lambda^{{s-\alpha+1+\sigma_1}}z\bigr |^2\bigr)+ \frac{\kappa
}{2}\bigl |
\Lambda^{s+\alpha}v^{(R)}\bigr |^2
\nonumber
\\
&&\qquad\leq C\chi_R\bigl(\bigl |\theta^{(R)}\bigr |^2_\mathcal{W}
\bigr) \bigl(C(R)+\bigl |\Lambda ^{{s-\alpha
+1+\sigma_1}}z\bigr |^2\bigr)+
\frac{\kappa}{2}\bigl |\Lambda^{s+\alpha
}v^{(R)}\bigr |^2,
\nonumber
\end{eqnarray}
where $r_1:=\frac{1-\alpha+\sigma_1}{\alpha}$. Here, in the first
equality, we used $\operatorname{div} u=0$, and in the second
inequality we used
Lemmas~\ref{lem2.1} and~\ref{lem2.2}, and in the fourth inequality we
used the
interpolation inequality (\ref{eq2.3}) and that $s-\alpha+1+\sigma
_1<s+2\alpha
-1$, and in the fifth inequality we used Young's inequality and in the
last inequality we used $|\Lambda^sv^{(R)}|\leq|\Lambda^s\theta
^{(R)}|+|\Lambda^{s-\alpha+1+\sigma_1}z|$.
Then as in the proof of Theorem~\ref{theA.1} in Appendix~\hyperref[appA]{A}, we prove
(\ref{eq6.5})
has a weak solution in $L^\infty([0,T],\mathcal{W})\cap
L^2([0,T],H^{s+\alpha})$.

\textit{Continuity}. For each $\omega\in\Omega^*$, $\sigma_1$ and
$\sigma_2$
as in (\ref{eq6.6}),
since $s-\alpha+1+\sigma_1<s+2\alpha-1$, we have $z\in C([0,\infty
);H^{s-\alpha+1+\sigma_1})$. Since $s>3-3\alpha$, multiplying the
equations (\ref{eq6.5}) by $\frac{d}{dt}\Lambda^{2(s-\alpha
)}v^{(R)}$, we obtain
%
%e6.7 #&#
%
\begin{eqnarray}
&&\frac{\kappa}{2}\frac{d}{dt}\bigl |\Lambda^{s}
v^{(R)}\bigr |^2+\bigl |\Lambda^{s-\alpha}\dot{v}^{(R)}\bigr |^2
\nonumber
\\
&&\qquad =
C\chi_R\bigl(\bigl |\theta^{(R)}\bigr |^2_\mathcal{W}
\bigr)\bigl\langle\Lambda^{s-\alpha
}\nabla \cdot\bigl(u^{(R)}
\theta^{(R)}\bigr),\Lambda^{s-\alpha}\dot{v}^{(R)}\bigr
\rangle
\nonumber
\\
&&\qquad\leq  C\chi_R\bigl(\bigl |\theta^{(R)}\bigr |^2_\mathcal{W}
\bigr)\bigl |\Lambda^{s-\alpha+1}\bigl(u^{(R)} \theta^{(R)}
\bigr)\bigr |\cdot\bigl |\Lambda^{s-\alpha}\dot{v}^{(R)}\bigr |
\nonumber
\\[-8pt]
\\[-8pt]
&&\qquad\leq  C\chi_R\bigl(\bigl |\theta^{(R)}\bigr |^2_\mathcal{W}
\bigr)\bigl |\Lambda^{{s-\alpha}+1+\sigma
_1}\theta^{(R)}\bigr |\bigl |\Lambda^{\sigma_2}
\theta^{(R)}\bigr |\cdot\bigl |\Lambda ^{s-\alpha
}\dot{v}^{(R)}\bigr |\nonumber
\\
&&\qquad\leq C\chi_R\bigl(\bigl |\theta^{(R)}\bigr |^2_\mathcal{W}
\bigr) \bigl(\bigl |\Lambda^{s+\alpha}v^{(R)}\bigr |^2+\bigl |
\Lambda^{s}v^{(R)}\bigr |^2+\bigl |\Lambda^
{s-\alpha+1+\sigma_1}z\bigr |^2
\bigr)
\nonumber
\\
&&\qquad\quad{}+ \frac{1}{2}\bigl |\Lambda^{s-\alpha}\dot {v}^{(R)}\bigr |^2.
\nonumber
\end{eqnarray}
Here, $\dot{v}^{(R)}=\frac{d{v}^{(R)}}{dt}$ and in the first equality
we used $\operatorname{div} u=0$, in the second inequality we
used Lemmas~\ref{lem2.1} and~\ref{lem2.2},
and in the third inequality we used the interpolation inequality (\ref{eq2.3}),
that $s-\alpha+1+\sigma_1\leq s+\alpha$ and Young's inequality.

As $\int_0^T|\Lambda^{s+\alpha}v^{(R)}(t)|^2\,dt$ can be dominated by
(\ref{eq6.6}), we get an a priori estimate for the time derivative
$\frac
{d}{dt}v^{(R)}$ in $L^2(0,T;H^{s-\alpha})$. Then by \cite{Te84}, we obtain
$v^{(R)}\in C([0,T],\mathcal{W})$.

\textit{Uniqueness}. Let $v_1,v_2$ be two solutions of equation (\ref
{eq6.5}) in
$C([0,\infty);\mathcal{W})$ and set $w:=v_1-v_2$ and $u_w:=u_1-u_2$,
where $u_1,u_2$ satisfy (\ref{eq1.3}) with $\theta$ replaced by
$\theta
_1=v_1+z, \theta_2=v_2+z$. Then by a similar argument as in the proof
of Theorem~\ref{the4.2}, we have for small $0<\varepsilon_1<(2\alpha
-1-\sigma
_1)\wedge\sigma_1$ with $\sigma_1$ as in (\ref{eq6.6})
\begin{eqnarray*}
&&\frac{1}{2}\frac{d}{dt}\bigl |\Lambda^{s-\alpha}
w\bigr |^2+\kappa\bigl |\Lambda^{s}w\bigr |^2
\\
&&\qquad= -\bigl(
\chi_R\bigl(|\theta_1|^2_\mathcal{W}
\bigr)-\chi_R\bigl(|\theta_2|^2_\mathcal
{W}\bigr)\bigr)\bigl\langle\Lambda^{s+\varepsilon_1-2\alpha}(u_1\cdot\nabla
\theta_1), \Lambda^{s-\varepsilon_1}w \bigr\rangle
\\
&&\qquad\quad{}- \chi_R\bigl(|\theta_2|^2_\mathcal
{W}\bigr)\bigl\langle\Lambda^{s-2\alpha}(u_1\cdot\nabla w),
\Lambda^{s}w \bigr\rangle
\\
&&\qquad\quad{}- \chi_R\bigl(|
\theta_2|^2_\mathcal{W}\bigr)\bigl\langle
\Lambda^{s-2\alpha
}(u_w\cdot \nabla\theta_2),
\Lambda^{s}w \bigr\rangle
\\
&&\qquad= I+\mathit{II}+\mathit{III}. %\end{array}
\end{eqnarray*}
As
\[
\bigl |\chi_R\bigl(|\theta_1|^2_\mathcal{W}
\bigr)-\chi_R\bigl(|\theta_2|^2_\mathcal
{W}\bigr)\bigr |\leq C(R)|w|_\mathcal{W}\bigl[1_{[0,R+1]}\bigl(|
\theta_1|_\mathcal {W}^2\bigr)+1_{[0,R+1]}
\bigl(|\theta_2|_\mathcal{W}^2\bigr)\bigr],
\]
we get for $\sigma_1,\sigma_2$ as in (\ref{eq6.6}),
%
%e6.8 #&#
%
\begin{eqnarray}
I&=& -\bigl(\chi_R\bigl(|\theta_1|^2_\mathcal{W}
\bigr)-\chi_R\bigl(|\theta_2|^2_\mathcal
{W}\bigr)\bigr)\bigl\langle\Lambda^{s+\varepsilon_1-2\alpha}\nabla\cdot(u_1
\theta_1), \Lambda^{s-\varepsilon_1}w \bigr\rangle
\nonumber
\\
&\leq& C\bigl[1_{[0,R+1]}\bigl(|\theta _1|_\mathcal{W}^2
\bigr)+1_{[0,R+1]} \bigl(|\theta_2|_\mathcal{W}^2
\bigr)\bigr]\nonumber
\\
&&{}\times |w|_\mathcal{W} \bigl |\Lambda^{{s-2\alpha+\varepsilon_1}+1+\sigma_1}\theta_1\bigr |\bigl |
\Lambda ^{\sigma
_2}\theta_1\bigr |\bigl  |\Lambda^{s-\varepsilon_1}w\bigr |
\\
&\leq& C\bigl(R,|\theta_1|_\mathcal{W},|\theta
_2|_\mathcal{W}\bigr)|w|_\mathcal{W}\bigl |\Lambda^{s-\varepsilon
_1}w\bigr |
\nonumber
\\
&\leq& C\bigl(R,|\theta_1|_\mathcal{W},|\theta_2|_\mathcal{W}\bigr)\bigl |
\Lambda ^{s-\alpha
}w\bigr |^2+\frac{\kappa}{4}|w|_\mathcal{W}^2,
\nonumber
\end{eqnarray}
where in the first equality we used $\operatorname{div} u_1=0$ and in
the first
inequality we used Lemmas~\ref{lem2.1} and~\ref{lem2.2}, in the
second inequality we used
that $s-2\alpha+\varepsilon_1+1+\sigma_1<s$, that is, $\varepsilon
_1<2\alpha-1-\sigma_1$ and in the third inequality we used the
interpolation inequality (\ref{eq2.3}) and Young's inequality.
In a similar way, we obtain
\begin{eqnarray*}
\mathit{II}&\leq&\bigl |\Lambda^sw\bigr |\bigl |\Lambda^{s-2\alpha+1}(u_1w)\bigr |
\\
&\leq&C\bigl |\Lambda^sw\bigr |\bigl[\bigl |\Lambda^{s-2\alpha+1+\sigma_1}
\theta_1\bigr | \bigl |\Lambda ^{s-\varepsilon_1}w\bigr |+\bigl |\Lambda^{s-2\alpha+1+\sigma_1}w\bigr | \bigl |
\Lambda ^{s}\theta_1\bigr |\bigr]
\\
&\leq& C\bigl(R,|\theta_1|_\mathcal{W}\bigr)\bigl |\Lambda^{s-\alpha}w\bigr |^2+
\frac
{\kappa
}{4}|w|_\mathcal{W}^2, %\end{array}
\end{eqnarray*}
where in the first inequality we used $\operatorname{div} u_1=0$ and in the second
inequality we used Lemmas~\ref{lem2.1} and~\ref{lem2.2}
and $s-\varepsilon_1\geq1-\sigma
_1$, and in the third inequality we used the interpolation inequality
(\ref{eq2.3}) and Young's inequality.
Similarly,
\[
\mathit{III}\leq C\bigl(R,|\theta_2|_\mathcal{W}\bigr)\bigl |
\Lambda^{s-\alpha
}w\bigr |^2+\frac{\kappa
}{4}|w|_\mathcal{W}^2.
\]
Then we obtain
\begin{eqnarray*}
&&\frac{1}{2}\frac{d}{dt}\bigl |\Lambda^{s-\alpha} w\bigr |^2+
\kappa\bigl |\Lambda ^{s}w\bigr |^2
\\
&&\qquad\leq C\Bigl(R,\sup
_{t\in[0,T]}\bigl |\theta_1(t)\bigr |_\mathcal{W},\sup
_{t\in
[0,T]}\bigl |\theta_2(t)\bigr |_\mathcal{W}\Bigr)\bigl |
\Lambda^{s-\alpha}w\bigr |^2+\frac
{3\kappa
}{4}|w|_\mathcal{W}^2.
\end{eqnarray*}
Gronwall's lemma now yields that $|\Lambda^{s-\alpha} w|=0$, which
implies $w=0$.

So, equation (\ref{eq6.5}) has a unique global weak solution in the space
$C([0,\infty);\break \mathcal{W})$.

Next, we prove (\ref{eq6.3}). In order to do so, it is sufficient to
show that
$P_x^{(R)}[\tau_R<\varepsilon]\leq C(\varepsilon, R)$
with $C(\varepsilon, R)\downarrow0$ as $\varepsilon\downarrow0$, for
all $x\in\mathcal{W}$, with $|x|_\mathcal{W}^2\leq\frac{R}{8}$.
So, fix $\varepsilon>0$ small enough, let $\Theta_{\varepsilon
,R}:=\sup_{t\in[0,\varepsilon]}|\Lambda^{s-\alpha+1+\sigma_1}z(t)|$ and assume
that $\Theta^2_{\varepsilon,R}\leq\frac{R}{8}$.
Setting $\varphi(t):=|v^{(R)}|^2_\mathcal{W}+\Theta_{\varepsilon,R}^2$,
by (\ref{eq6.6}) we get $\dot{\varphi}\leq C(R)$.
This implies, together with the bounds on $x$ and $\Theta_{\varepsilon
,R}$, that
\[
\sup_{t\in[0,\varepsilon]}\bigl |\theta^{(R)}(t)\bigr |^2_\mathcal{W}
\leq R
\]
for $\varepsilon$ small enough.
It follows that $\tau_R\geq\varepsilon$.
Hence,
\[
P_x^{(R)}[\tau_R<\varepsilon]\leq
P_x^{(R)}\biggl[\sup_{t\in
[0,\varepsilon
]}\bigl |
\Lambda^{s+1+\sigma_1-\alpha}z(t)\bigr |^2>\frac{R}{8}\biggr].
\]
Letting $\varepsilon\downarrow0$, we have $P_x^{(R)}[\tau
_R<\varepsilon
]\rightarrow0$, and the claim is proved,
since the probability above is independent of $x$.

Finally, the same arguments as in the proof of Theorem~\ref{the4.2}
imply that
\[
\theta\bigl(t\wedge\tau_R\bigl(\theta^{(R)}\bigr)\bigr)=
\theta^{(R)}\bigl(t\wedge\tau _R\bigl(\theta
^{(R)}\bigr)\bigr)\qquad\forall t, P\mbox{-a.s.}
\]
Moreover, since $\theta$ is $H$-valued weakly continuous, we obtain
$\tau_R(\theta^{(R)})=\tau_R(\theta)$.
\end{pf}
%
%$\hfill\Box$

In order to apply \cite{FR08}, Theorem~5.4, we now only need the
following result.

\begin{Proposition}\label{pro6.1.4}
Fix $\alpha>\frac{2}{3}$. Suppose Hypothesis~\ref{hypE.3}
holds. For every $R>0$, the transition semigroup
$(P_t^{(R)})_{t\geq0}$ associated to equation (\ref{eq6.2}) is
$\mathcal
{W}$-strong Feller.
\end{Proposition}

\begin{pf}
We shall provide formal estimates, that can, however, be made
rigorous through Galerkin approximations.
Let $(\Sigma,\mathcal{F},(\mathcal{F}_t)_{t\geq0}, \mathbb{P})$ be a
filtered probability space, $(W_t)_{t\geq0}$ a cylindrical Wiener
process on $H$ and, for every $x\in\mathcal{W}$, let $\theta_x^{(R)}$
be the solution to equation (\ref{eq6.2}) with initial value $x\in
\mathcal{W}$.
By the Bismut--Elworthy--Li formula,
\[
D_y\bigl(P_t^{(R)}\psi\bigr) (x)=
\frac{1}{t}E^{\mathbb{P}} \biggl[\psi\bigl(\theta _x^{(R)}(t)
\bigr)\int_0^t\bigl\langle G^{-1}D_y
\theta_x^{(R)}(l), dW(l)\bigr\rangle \biggr],
\]
where $D_y(P_t^{(R)}\psi)$ denotes $\langle D(P_t^{(R)}\psi),
y\rangle$
for $y\in H$, $D_y\theta_x^{(R)}=D\theta_x^{(R)}\cdot y$ and $D\theta
_x^{(R)}$ denotes the derivative of $\theta_x^{(R)}$ with respect to
the initial value.
Then for $\|\psi\|_\infty\leq1$, by the B--D--G inequality
\begin{eqnarray*}
&&\bigl |\bigl(P_t^{(R)}\psi\bigr) (x_0+h)-
\bigl(P_t^{(R)}\psi\bigr) (x_0)\bigr |
\\
&&\qquad\leq
\frac{C}{t}\sup_{\eta
\in[0,1]}E^{\mathbb{P}} \biggl[\biggl(
\int_0^t\bigl |G^{-1}D_{h}
\theta_{x_0+\eta
h}^{(R)}(l)\bigr |^2\,dl\biggr)^{1/2}
\biggr].
\end{eqnarray*}
The proposition is proved once we prove that the right-hand side of the
above inequality converges to $0$ as $|h|_\mathcal{W}\rightarrow0$.

Fix $x\in\mathcal{W}$, $h\in H$ and write $\theta=\theta_x^{(R)},
v=v^{(R)}, u=u^{(R)}, D\theta=D_h\theta$ for simplicity. The term
$D\theta$ solves the following equation:
\begin{eqnarray*}
&&\frac{d}{dt}D\theta+\kappa\Lambda^{2\alpha}(D\theta)
\\
&&\qquad=-\bigl[\chi
_R\bigl(|\theta |_\mathcal{W}^2\bigr)[ Du\cdot
\nabla\theta+u\cdot\nabla D\theta]+2\chi _R'\bigl(|
\theta|_\mathcal{W}^2\bigr)\langle\theta,D\theta
\rangle_\mathcal {W}u\cdot\nabla\theta\bigr],
\end{eqnarray*}
with initial value $D\theta(0)=h$ and $Du$ satisfying (\ref{eq1.3})
with $\theta
$ replaced by $D\theta$.
Multiplying the above equation with $\Lambda^{2s}D\theta$ and taking
the inner product in $L^2$, we have
\begin{eqnarray*}
&&\frac{1}{2}\frac{d}{dt}\bigl |\Lambda^s D\theta\bigr |^2+
\kappa\bigl |\Lambda ^{s+\alpha
}(D\theta)\bigr |^2
\\
&&\qquad=-\bigl\langle \bigl[
\chi_R\bigl(|\theta|_\mathcal{W}^2\bigr)[ Du\cdot
\nabla\theta+u\cdot\nabla D\theta]+2\chi_R'\bigl(|
\theta|_\mathcal {W}^2\bigr)\langle\theta,D\theta
\rangle_\mathcal{W}u\cdot\nabla\theta \bigr],
\\
&&\phantom{\hspace*{293pt}}\Lambda^{2s}D\theta
\bigr\rangle.
\end{eqnarray*}
For the first term on the right-hand side, we have for $|\theta
|_\mathcal{W}^2\leq R$
%
%e6.9 #&#
%
\begin{eqnarray}
\label{eq6.9} %\begin{array}{cc}
\bigl |\bigl\langle Du\cdot\nabla\theta,\Lambda
^{2s}D\theta\bigr\rangle\bigr |&=& \bigl |\bigl\langle\Lambda^{s-\alpha}\nabla
\cdot(Du \theta),\Lambda ^{s+\alpha
}D\theta\bigr\rangle\bigr |
\nonumber
\\
&\leq&C\bigl |\Lambda^{{s-\alpha}+1+\sigma_1}\theta\bigr |\cdot \bigl |\Lambda^{\sigma_2} D\theta\bigr |
\cdot \bigl |\Lambda^{s+\alpha}D\theta\bigr |\nonumber
\\
&&{}+C\bigl |\Lambda^{{s-\alpha}+1+\sigma
_1}D\theta \bigr |\cdot\bigl |
\Lambda^{\sigma_2} \theta\bigr |\cdot \bigl |\Lambda^{s+\alpha}D\theta\bigr |
\\
&\leq& \varepsilon\bigl |\Lambda^{s+\alpha}D\theta\bigr |^2\nonumber
\\
&&{}+C\bigl(C(R)+\bigl |
\Lambda^{s+\alpha} v\bigr |^2+\bigl |\Lambda^{s-\alpha+1+\sigma_1}z\bigr |^2
\bigr)\bigl |\Lambda^s D\theta \bigr |^2
\nonumber
\end{eqnarray}
for $\sigma_1,\sigma_2$ as (\ref{eq6.6}), where we used $\operatorname{div} Du=0$
in the first
equality and Lemmas~\ref{lem2.1} and~\ref{lem2.2} in the first
inequality as well as the
interpolation inequality (\ref{eq2.3}) and Young's inequality in the
second inequality.

The second term can be estimated similarly. For the third term, by
Lemmas~\ref{lem2.1} and~\ref{lem2.2}, we have
%
%e6.10 #&#
%
\begin{eqnarray}\label{eq6.10}
\bigl |\bigl\langle u\cdot\nabla\theta,\Lambda^{2s}D\theta \bigr
\rangle\bigr |&=&\bigl |\bigl\langle\Lambda^{s-\alpha}\nabla\cdot( u \theta ),\Lambda
^{s+\alpha}D\theta\bigr\rangle\bigr |
\nonumber
\\
&\leq&C\bigl |\Lambda^{{s-\alpha}+1+\sigma_1}\theta \bigr |\bigl |\Lambda^{\sigma_2}\theta\bigr |\cdot\bigl |
\Lambda^{s+\alpha} D\theta\bigr |
\\
&\leq& C\bigl(\bigl |\Lambda^{s+\alpha} v\bigr |+\bigl |\Lambda^{s-\alpha+1+\sigma_1}z\bigr |\bigr)\bigl |
\Lambda ^s\theta\bigr |\bigl |\Lambda^{s+\alpha}D\theta\bigr |,
\nonumber
\end{eqnarray}
where in the first equality we used $\operatorname{div}u=0$.
Then we obtain
\begin{eqnarray*}
&&\frac{1}{2}\frac{d}{dt}\bigl |\Lambda^s D
\theta\bigr |^2+\kappa\bigl |\Lambda ^{s+\alpha}(D\theta)\bigr |^2
\\
&&\qquad\leq
\frac{\kappa}{2}\bigl |\Lambda^{s+\alpha
}(D\theta )\bigr |^2+C\bigl(C(R)+\bigl |
\Lambda^{s+\alpha} v\bigr |^2+\bigl |\Lambda^{s-\alpha+1+\sigma
_1}z\bigr |^2
\bigr)\bigl |\Lambda^s D\theta\bigr |^2. %\end{array}
\end{eqnarray*}
From Gronwall's inequality and (\ref{eq6.6}), we finally get
\begin{eqnarray*}
&&\int_0^t\bigl |\Lambda^{s+\alpha}
\bigl(D\theta(l)\bigr)\bigr |^2\,dl
\\
&&\qquad\leq C\bigl |\Lambda ^sh\bigr |^2+
\exp\biggl(C\int_0^t\bigl(C(R)+\bigl |
\Lambda^{s+\alpha} v\bigr |^2+\bigl |\Lambda ^{s-\alpha
+1+\sigma_1}z\bigr |^2
\,dl\bigr)\biggr)\bigl |\Lambda^s h\bigr |^2
\\
&&\qquad\leq  C\bigl |\Lambda^sh\bigr |^2+\exp \biggl(C\biggl(\bigl |
\Lambda^{s} x \bigr |^2+\int_0^t
\bigl(C(R)+\bigl |\Lambda^{s-\alpha+1+\sigma
_1}z\bigr |^2\,dl\bigr)\biggr)\biggr)\bigl |
\Lambda^s h\bigr |^2. %\end{array}
\end{eqnarray*}
Since by $s-\alpha+1+\sigma_1<s+2\alpha-1$, $z$ is a Gaussian random
variable in $C([0,\infty);H^{s-\alpha+1+\sigma_1})$ (cf. \cite{D04},
Proposition~2.15), by Fernique's theorem we could choose $t_0$ small
enough and obtain
\[
E\int_0^{t_0}\bigl |\Lambda^{s+\alpha}\bigl(D
\theta(l)\bigr)\bigr |^2\,dl\leq c(t_0,R)\bigl |\Lambda
^s h\bigr |^2,
\]
which, as $G^{-1}=Q_0^{-1/2}\Lambda^{s+\alpha}$, implies the assertion
for $t_0$. For general $t$, by the semigroup property the assertion
follows easily.
\end{pf}
%
%$\hfill\Box$

%s6.2 #&#
\subsection{A support theorem for \texorpdfstring{$\alpha>2/3$}{$alpha>2/3$}}

A Borel probability measure $\mu$ on $H$ is fully supported on
$\mathcal
{W}$ if $\mu(U)>0$ for every nonempty
open set $U\subset\mathcal{W}$. Set $\mathcal{W}_1:=H^{s-\alpha
+1+\sigma_1}$, where $\sigma_1$ is the same as
(\ref{eq6.6}) and we will use it below.

\begin{Lemma}[(Approximate controllability)]\label{lem6.2.1}
Let $R>0$, $T>0$. Let
$x\in\mathcal{W}$ and $y\in\mathcal{W}$, with $A_\alpha y\in
\mathcal
{W}_1$, such that
\[
|x|_\mathcal{W}^2\leq\frac{R}{2},\qquad|y|_{\mathcal{W}}^2
\leq \frac{R}{2}.
\]
Then there exist (a control function) $\omega\in\operatorname{Lip}([0,T];\mathcal{W}_1)$ and
\[
\theta\in C\bigl([0,T];\mathcal{W}\bigr)\cap L^2\bigl([0,T];
H^{s+\alpha}\bigr),
\]
such that $\theta$ solves the equation
%
%e6.11 #&#
%
\begin{eqnarray}
\label{eq6.11}
&&\theta(t)-x+\int_0^tA_\alpha
\theta(r)+\chi_R\bigl(|\theta |_\mathcal{W}^2
\bigr)u(r)\cdot\nabla\theta(r) \,dr
\nonumber
\\[-8pt]
\\[-8pt]
&&\qquad=\omega(t) \qquad dt\mbox{-a.e. } t\in[0,T],\nonumber
\end{eqnarray}
with $\theta(0)=x$ and $\theta(T)=y$, and
%
%e6.12 #&#
%
\begin{equation}
\label{eq6.12} \sup_{t\in[0,T]}\bigl |\theta(t)\bigr |_\mathcal{W}^2
\leq R.
\end{equation}
\end{Lemma}

\begin{pf}
First consider $\omega=0$. By similar arguments as in Theorems~\ref{theA.1}
and~\ref{theA.2},
there exist a unique solution $\theta\in C([0,T],\mathcal
{W})$. Then by a similar calculation as (\ref{eq6.6}), we get
\[
\frac{d}{dt}|\theta|_\mathcal{W}^2+\kappa\bigl |
\Lambda^{\alpha}\theta \bigr |_\mathcal{W}^2\leq C(R).
\]
Hence, $\theta(t)\in H^{s+\alpha}$ for almost every $t\in[0,T]$ and, by
solving again the equation with one of these regular points as initial
condition and using Lemmas~\ref{lem2.1} and~\ref{lem2.2} we have
\begin{eqnarray*}
\frac{d}{dt}\bigl |\Lambda^{\alpha+s}\theta\bigr |^2+
\kappa\bigl |\Lambda ^{2\alpha+s}\theta\bigr |_\mathcal{W}^2&=&
\chi_R\bigl(|\theta|_\mathcal{W}^2\bigr) \bigl
\langle\Lambda^{s}\nabla\cdot(u\theta),\Lambda^{2\alpha+s}\theta
\bigr\rangle
\\
&\leq& C\chi_R\bigl(|\theta|_\mathcal{W}^2
\bigr)\bigl |\Lambda^{2\alpha+s}\theta \bigr |\bigl |\Lambda^{s+1+\sigma_3}\theta\bigr |\|\theta
\|_{L^p}
\\
&\leq& C(R)\bigl |\Lambda ^{s+\alpha}\theta\bigr |^2+\frac{\kappa}{2}\bigl |
\Lambda^{s+2\alpha}\theta \bigr |^2, %\end{array}
\end{eqnarray*}
where $\sigma_3=\frac{2}{p}<2\alpha-1$ and we used $\operatorname{div}u=0$ in the
first equality and $H^s\subset L^p$ and the interpolation inequality
(\ref{eq2.3}), Young's inequality in the last step. Then\vspace*{1pt} by a boot
strapping argument,
we find a small $T_*\in(0,\frac{T}{2})$ such that $|\theta
(t)|^2_\mathcal{W}\leq R$ and $A_\alpha\theta(T_*)\in\mathcal{W}_1$
for all $t\leq T_*$. Define $\theta$ to be the solution above for
$t\in
[0,T_*]$ and extended by linear interpolation between $y$ and $\theta
(T_*)$ in $[T_*,T]$. Then obviously (\ref{eq6.12}) follows.

Next, if we set
\[
\eta:=\partial_t\theta+A_\alpha\theta+\chi_R
\bigl(|\theta|_\mathcal {W}^2\bigr)u\cdot\nabla\theta,\qquad
T_*\leq t\leq T,
\]
$\omega:=0$ for $t\leq T_*$ and $\omega(t)=\int_{T_*}^t\eta_s\,ds$ for
$t\in[T_*,T]$, we also have (\ref{eq6.11}). It remains to prove that
$\eta\in
L^\infty(0,T;\mathcal{W}_1)$. For
the first two terms of $\eta$, this is obvious. For the nonlinear term,
we have that
\[
|u\cdot\nabla\theta|_{\mathcal{W}_1}=\bigl |\nabla\cdot(u \theta )\bigr |_{\mathcal
{W}_1}
\leq C\bigl |\Lambda^{2\alpha}\theta\bigr |^2_{\mathcal{W}_1}
\]
for any $\theta\in\mathcal{W}_1$, where in the first equality we used
$\operatorname{div} u=0$ and in the last step we used Lemma~\ref{lem2.1}.
\end{pf}
%
%$\hfill\Box$

Let $l\in(0,\frac{1}{2})$ and $p>1$ such that $l-\frac{1}{p}>0$. Under
Hypothesis~\ref{hypE.3},
we see that for every $\alpha_1<\frac{s+\alpha-1}{2\alpha}$ the map
\[
\omega\mapsto z(\cdot,\omega)\dvtx W^{l,p}\bigl([0,T];D
\bigl(A_\alpha ^{\alpha
_1}\bigr)\bigr)\rightarrow C\bigl([0,T];D
\bigl(A_\alpha^{\alpha_1+l-{1}/{p}-\varepsilon}\bigr)\bigr)
\]
is continuous, for all $\varepsilon>0$ (cf. \cite{DZ96}), where $z$
is the
solution to the following equation:
%
%e6.13 #&#
%
\begin{equation}
\label{eq6.13} z(t)+\int_0^tA_\alpha
z(s)\,ds=\omega(t).
\end{equation}
In particular, it is possible to find $\alpha_1\in(0,\frac{s+\alpha
-1}{2\alpha})$ and $p$
such that the above map is continuous from $W^{l,p}([0,T];D(A_\alpha
^{\alpha_1}))$ to $C([0,T];H^{s-\alpha+1+\sigma_1})$ since $\alpha
>\frac
{2}{3}$ and $\sigma_1<3\alpha-2$.

\begin{Lemma}[(Continuity with respect to the control functions)]\label
{lem6.2.2}
Let
$l$, $p$ and $\alpha_1$ be chosen as above, and let $\omega
_n\rightarrow
\omega$ in $W^{l,p}([0,T];D(A_\alpha^{\alpha_1}))$. Let $\theta$ be the
solution to equation (\ref{eq6.11}) corresponding to $\omega$ and
some initial
condition $x\in\mathcal{W}$ (the solution exists by the same arguments
as the proof of Theorem~\ref{theA.1}), and let
\[
\tau=\inf\bigl\{t\geq0\dvtx\bigl |\theta(t)\bigr |_\mathcal{W}^2\geq R
\bigr\},
\]
where as usual we set $\inf\emptyset=\infty$.
For each $n\in\mathbb{N}$, define similarly $\theta_n$ and $\tau_n$
corresponding to $\omega_n$ with the same initial condition $x$. If
$\tau>T$, then $\tau_n>T$ for $n$ large enough and
\[
\theta_n\rightarrow\theta\qquad\mbox{in } C\bigl([0,T];\mathcal{W}
\bigr).
\]
\end{Lemma}

\begin{pf}
Set $v_n:=\theta_n-z_n$ for each $n\in\mathbb{N}$, and
$v:=\theta-z$, where $z_n, z$ are the solutions to (\ref{eq6.13})
corresponding
to $\omega_n, \omega$, respectively. Since $\omega_n\rightarrow
\omega$
in $W^{l,p}([0,T];D(A_\alpha^{\alpha_1}))$, we can find a common lower
bound for $(\tau_n)_{n\in\mathbb{N}}$ and $\tau$ by (\ref{eq6.6}).
For every
time smaller than this lower bound $t_0$, by (\ref{eq6.6}), we have
\[
\sup_{(0,t_0)}\bigl |\Lambda^{s}\theta_n\bigr |^2
\leq R,\qquad\sup_{(0,t_0)}\bigl |\Lambda ^{s}
\theta\bigr |^2\leq R,\qquad\sup_{(0,t_0)}\bigl |
\Lambda^{s-\alpha
+1+\sigma
_1}z_n\bigr |\leq C,
\]
and
\begin{eqnarray*}
\sup_{(0,t_0)}\bigl |\Lambda^{s-\alpha+1+\sigma_1}z\bigr |&\leq& C,
\\
\int
_0^{t_0}\bigl |\Lambda^{s+\alpha}v_n(l)\bigr |^2
\,dl&\leq& C(R),
 \qquad
\int_0^{t_0}\bigl |\Lambda
^{s+\alpha}v(l)\bigr |^2\,dl \leq  C(R),
\end{eqnarray*}
where $C(R)$ is a constant depending only on $R$.
Moreover, we obtain for $t\leq t_0$
\begin{eqnarray*}
&&\frac{d}{dt}|v-v_n|^2_\mathcal{W}+2
\kappa\bigl |\Lambda^\alpha (v_n-v)\bigr |_\mathcal{W}^2
\\
&&\qquad =
\bigl\langle u_n\cdot\nabla\theta_n, \Lambda
^{2s}(v-v_n)\bigr\rangle-\bigl\langle u\cdot\nabla\theta,
\Lambda ^{2s}(v-v_n)\bigr\rangle
\\
&&\qquad = \bigl[\bigl\langle(u_{v_n}-u_v)\cdot\nabla
\theta_n,\Lambda ^{2s}(v-v_n)\bigr\rangle+
\bigl\langle u\cdot\nabla(v_n-v),\Lambda ^{2s}(v-v_n)
\bigr\rangle
\\
&&\qquad\quad  {}+\bigl\langle(u_{z_n}-u_z)\cdot\nabla
\theta_n,\Lambda ^{2s}(v-v_n)\bigr\rangle+
\bigl\langle u\cdot\nabla(z_n-z),\Lambda ^{2s}(v-v_n)
\bigr\rangle\bigr], %\end{array}
\end{eqnarray*}
where $u_{v_n},u_{z_n}$ satisfy (\ref{eq1.3}) with $\theta$ replaced by
$v_n,z_n$, respectively.
For the first term on the right-hand side, we have
\begin{eqnarray*}
&&\bigl |\bigl\langle(u_{v_n}-u_v)\cdot\nabla
\theta_n,\Lambda ^{2s}(v-v_n)\bigr\rangle\bigr |
\\
&&\qquad = \bigl |
\bigl\langle\Lambda^{s-\alpha}\nabla\cdot \bigl((u_{v_n}-u_v)
\theta_n\bigr),\Lambda^{s+\alpha}(v-v_n)\bigr
\rangle\bigr |
\\
&&\qquad \leq  C\bigl |\Lambda^{s+\alpha}(v-v_n)\bigr |\bigl |\Lambda^{{s-\alpha}+1+\sigma
_1}(v-v_n)\bigr |\bigl |
\Lambda^{\sigma_2}\theta_n\bigr |
\\
&&\qquad\quad  {}+C\bigl |\Lambda^{s+\alpha}(v-v_n)\bigr |\bigl |\Lambda^{{s-\alpha}+1+\sigma
_1}
\theta _n\bigr |\bigl |\Lambda^{\sigma_2}(v-v_n)\bigr |
\\
&&\qquad \leq  \frac{\kappa}{4}\bigl |\Lambda^{s+\alpha
}(v-v_n)\bigr |^2+C
\bigl(C(R)+\bigl |\Lambda^{s+\alpha}v_n\bigr |^2\bigr)\bigl |
\Lambda^s(v-v_n)\bigr |^2
\\
&&\qquad\quad  {}+c\bigl |\Lambda^{{s-\alpha}+1+\sigma_1}z_n\bigr |^2\bigl |
\Lambda^s(v-v_n)\bigr |^2.
\end{eqnarray*}
Here, $\sigma_1,\sigma_2$ are as (\ref{eq6.6}) and we used
$\operatorname{div}(u_{v_n}-u_v)=0$
in the first equality and Lemmas~\ref{lem2.1} and~\ref{lem2.2} in the
first inequality
and the interpolation inequality (\ref{eq2.3}) and Young's inequality
in the
last step. The other term can be estimated similarly. Then we obtain
\begin{eqnarray*}
&&\frac{d}{dt}|v-v_n|^2_\mathcal{W}+2
\kappa\bigl |\Lambda^\alpha (v_n-v)\bigr |_\mathcal{W}^2
\\
&&\qquad\leq\kappa\bigl |\Lambda^\alpha(v_n-v)\bigr |_\mathcal
{W}^2
\\
&&\qquad\quad{}+C\bigl(C(R)+\bigl |\Lambda^\alpha v_n\bigr |^2_\mathcal{W}+\bigl |
\Lambda^\alpha v\bigr |^2_\mathcal{W}\bigr)
\bigl(|v-v_n|^2_\mathcal{W}+\bigl |\Lambda^{s-\alpha+1+\sigma
_1}(z-z_n)\bigr |^2
\bigr). %\end{array}
\end{eqnarray*}
Then Gronwall's lemma yields that
\begin{eqnarray*}
|v-v_n|^2_\mathcal{W}&\leq&\Theta_n
\exp\biggl(C\int_0^t\bigl(C(R)+\bigl |\Lambda
^\alpha v_n\bigr |_\mathcal{W}^2+\bigl |
\Lambda^\alpha v\bigr |_\mathcal{W}^2\bigr)\,dl\biggr)
\\
&&{}\times\int
_0^t\bigl(C(R)+\bigl |\Lambda^\alpha
v_n\bigr |_\mathcal{W}^2+\bigl |\Lambda^\alpha
v\bigr |_\mathcal {W}^2\bigr)\,dl,
\end{eqnarray*}
where $\Theta_n=\sup_{[0,T]}|\Lambda^{s-\alpha+1+\sigma
_1}(z-z_n)|$. We
conclude $\theta_n\rightarrow\theta$ in $C([0,T];\mathcal{W})$. Now,
since $\tau>T$, if
$S=\sup_{t\in[0,T]}|\Lambda^s\theta(t)|^2$, then $S<R$ and we find
$\delta>0$ (depending only on $R$ and $S$) and $n_0\in\mathbb{N}$ such
that $\Theta^2_n<\delta$ and
$|v_n-v|^2_\mathcal{W}<\delta$ for all $n\geq n_0$, and so
\[
\bigl |\theta_n(t)\bigr |_\mathcal{W}\leq\bigl |v_n(t)-v(t)\bigr |_\mathcal{W}+
\Theta _n+\bigl |\theta (t)\bigr |_\mathcal{W}\leq2\sqrt{\delta}+\sqrt{S}
\leq\sqrt{R-\delta}.
\]
Then $\tau_n>T$ for all $n\geq n_0$.
\end{pf}
%
%$\hfill\Box$

\begin{Theorem}\label{the6.2.3}
Fix $\alpha>\frac{2}{3}$. Suppose Hypothesis~\ref{hypE.3}
holds and for $x\in\mathcal{W}$ let $P_x$ be the distribution of the
solution of (\ref{eq3.1}) with initial value $\theta(0)=x$. Then for every
$x\in\mathcal{W}$ and every $T>0$, the image measure of $P_x$ at time
$T$ is fully supported on $\mathcal{W}$.
\end{Theorem}

\begin{pf} Fix $x\in\mathcal{W}$ and $T>0$. We need to show that for every
$y\in\mathcal{W}$ and $\varepsilon>0$, $P_x[|\theta_T-y|_\mathcal
{W}<\varepsilon]>0$.
Let $\bar{y}\in\mathcal{W}\cap D(A_\alpha)$ such that $A_\alpha
\bar
{y}\in\mathcal{W}_1$ and $|y-\bar{y}|_\mathcal{W}<\frac
{\varepsilon
}{2}$. Choose $R>0$ such that $3|x|_\mathcal{W}^2<R$ and
$3|y|_\mathcal{W}^2<R$.
Then by Theorem~\ref{the6.1.3},
\begin{eqnarray*}
P_x\bigl[|\theta_T-y|_\mathcal{W}<
\varepsilon\bigr]&\geq&P_x\biggl[|\theta _T-
\bar{y}|_\mathcal{W}<\frac{\varepsilon}{2}\biggr]\geq P_x
\biggl[|\theta _T-\bar {y}|_\mathcal{W}<\frac{\varepsilon}{2},
\tau_R>T\biggr]
\\
&=&P_x^{(R)}\biggl[|\theta_T-
\bar{y}|_\mathcal{W}<\frac{\varepsilon
}{2},\tau _R>T\biggr].
\end{eqnarray*}
By Lemma~\ref{lem6.2.1}, there is a control $\bar{\omega}\in
W^{l,p}([0,T];D(A_\alpha^{\alpha_1}))$, with $l,p$ and $\alpha_1$
chosen as in Lemma~\ref{lem6.2.2}, such that the solution $\bar
{\theta}$ to the
control problem (\ref{eq6.11}) corresponding to $\bar{\omega}$
satisfies $\bar
{\theta}(0)=x, \bar{\theta}(T)=\bar{y}$ and $|\bar{\theta
}(t)|_\mathcal
{W}^2\leq\frac{2}{3}R$. By Lemma~\ref{lem6.2.2}, there exists
$\delta>0$ such
that for all $\omega\in W^{l,p}([0,T];D(A_\alpha^{\alpha_1}))$ with
$|\omega-\bar{\omega}|_{W^{l,p}([0,T];D(A_\alpha^{\alpha
_1}))}<\delta$,
we have
\[
\bigl |\theta(T,\omega)-\bar{y}\bigr |_\mathcal{W}<\frac{\varepsilon}{2} \quad
\mbox{and}\quad\sup_{t\in[0,T]}\bigl |\theta(t,\omega )\bigr |^2_\mathcal{W}<R,
\]
where $\theta(\cdot,\omega)$ is the solution to the control problem
(\ref{eq6.11}) corresponding to $\omega$ and starting at $x$. Hence,
\[
P_x^{(R)}\biggl[|\theta_T-
\bar{y}|_\mathcal{W}<\frac{\varepsilon
}{2},\tau _R>T\biggr]\geq
P_x^{(R)}\bigl[|\eta-\bar{\omega}|_{W^{l,p}([0,T];D(A_\alpha
^{\alpha_1}))}<\delta\bigr],
\]
where $\eta_t=\theta_t-x+\int_0^t(A_\alpha\theta_s+\chi_R(|\theta
_s|_\mathcal{W}^2)u\cdot\nabla\theta_s)\,ds$, hence $\theta_T=\theta
(T,\eta)$, and the right-hand side of the inequality above is strictly
positive since by Hypothesis~\ref{hypE.3}
$\eta$ is a Gaussian process in $D(A_
\alpha^{\alpha_1})$.
\end{pf}

\begin{Theorem}\label{the6.2.4}
Let $\alpha>\frac{2}{3}$ and suppose Hypothesis~\ref{hypE.3}
holds. Then there exists a unique invariant measure $\nu$ on $\mathcal
{W}$ for the transition semigroup $(P_t)_{t\geq0}$. Moreover:
\begin{enumerate}[(iii)]
\item[(i)]The invariant measure $\nu$ is ergodic.

\item[(ii)]The transition semigroup $(P_t)_{t\geq0}$ is $\mathcal{W}$-strong
Feller, irreducible and, therefore, strongly mixing. Furthermore,
$P_t(x,dy),t>0,x\in\mathcal{W}$, are mutually equivalent.

\item[(iii)]There exist $0<\delta_1<\frac{s+\alpha-1}{2\alpha}$
and $0<\gamma
_0<1$ such that
\[
\int\bigl |A_\alpha^{\delta_1} x\bigr |^{2\gamma_0}_\mathcal{W}\,d
\nu<\infty.
\]
\end{enumerate}
\end{Theorem}

\begin{pf}
By similar methods as the proof of Theorem~\ref{the5.3.2}, we obtain the
existence of the invariant measures. In fact, under Hypothesis~\ref{hypE.3},
we could
choose the following approximation:
\[
d\theta_n(t)+A_\alpha\theta_n(t)
\,dt+u_n(t)\cdot\nabla\theta _n(t)\,dt=
k_{\delta_n}*G\,dW(t),
\]
with initial data $\theta_n(0)=x\in H^s, u_n$ satisfying (\ref
{eq1.3}) with
$\theta$ replaced by $\theta_n$ and $k_{\delta_n}$ is the periodic
Poisson kernel as in the proof of Theorem~\ref{the3.3}. By the same
arguments as
Theorems~\ref{theA.1} and~\ref{theA.2},
we obtain that there exist a unique solution to
the above equation with $\theta_n\in C([0,\infty),H^s)\cap
L^2_{\mathrm{loc}}([0,\infty),H^{s+\alpha})$ $P$-a.s.
Then do the same
calculations for $\theta_n$ as in Lemma~\ref{lem5.3.1}, we obtain that
there exists $0<\gamma_0<1, 0<\tilde{\delta}_1<s+\alpha-1$ such that
\[
E\biggl[\int_0^t\bigl |
\Lambda^{\tilde{\delta}_1+s} \theta_n\bigr |^{2\gamma_0} \,dr\biggr]\leq
C(1+t) \bigl(\|x\|_{L^q}^q+1\bigr). %\end{array}
\]
Choose $x_0\in H^1$ and define
\[
\mu_t=\frac{1}{t}\int_0^tP_r^*
\delta_{x_0} \,dr.
\]
Since by similar arguments as in the proof of Theorem~\ref{theA.1}, we have
$P$-a.s. $\theta_n\rightarrow\theta$ in $L^2([0,T],H)$ and for
$2\alpha
\delta_1\leq\tilde{\delta}_1$, $0<\gamma_0<1$
\[
\int\bigl |A_\alpha^{\delta_1} x\bigr |_{H^s}^{2\gamma_0}
\mu_t(dx)=\frac
{1}{t}E_{x_0}\biggl[\int
_0^t\bigl |A_\alpha^{\delta_1}
\theta\bigr |_{H^s}^{2\gamma
_0}\, dr\biggr],
\]
by the above estimates we have for $t>1$
\[
\int\bigl |A_\alpha^{\delta_1} x\bigr |^{2\gamma_0}_{H^s}
\mu_t(dx)\leq C.
\]
This implies that $\mu_t$ is tight on $ H^s$.
Hence, any limit point of $\mu_t$ is an invariant measure for
$(P_t)_{t\geq0}$.
Therefore, by Doob's theorem, the strongly mixing property is a
consequence of Theorem~\ref{the6.1.2} and Theorem~\ref{the6.2.3}.
\end{pf}

\begin{remark}[(Mildly degenerate noise)]\label{rem6.2.5}
We can also consider the
ergodicity of the equation driven by a mildly degenerate noise as in
\cite{EH01}.
For this, we have to use an extension of the Bismut--Elworthy--Li
formula. We have the same problem as explained in Remark~\ref
{rem6.1.1}. So, we
can just get the result for $\alpha>2/3$.
\end{remark}

%s6.3 #&#
\subsection{Exponential convergence for \texorpdfstring{$\alpha>\frac{2}{3}$}{$alpha>\frac{2}{3}$}}

In this subsection, we assume that $\alpha>\frac{2}{3}$ and
$s>3-2\alpha
$. Then under Hypothesis~\ref{hypE.3} the associated O--U process
$z\in C([0,\infty),
H^{2+\delta_0})$ for some $0<\delta_0<s+2\alpha-3$.

\begin{Lemma}\label{lem6.3.1}
Fix $\alpha>2/3$. Let $\theta$ denote the solution of
(\ref{eq3.1}) and take $p>\frac{2}{3\alpha-2}$, then for every
$R_0\geq1$,
there exist values $T_1=T_1(R_0)$ and $\tilde{C}_1=\tilde{C}_1(R_0)$
such that if $\sup_{t\in[0,T_1]}\|\theta(t)\|_{L^p}^p\leq R_0$, and
$\sup_{t\in[0,T_1]}|\Lambda^{s+2\alpha-1-\varepsilon}z(t)|^2\leq R_0$
for some $0<\varepsilon<3\alpha-2-\frac{2}{p}$, then $|\Lambda
^{s+\delta}\theta(T_1)|^2\leq\tilde{C}_1$ for some $\delta>0$.
\end{Lemma}

\begin{pf}
For $v=\theta-z$, we have the following estimate:
\begin{eqnarray*}
\frac{1}{2}\frac{d}{dt}|v|^2+\kappa\bigl |
\Lambda^\alpha v\bigr |^2&=&\bigl\langle-u \cdot\nabla(v+z),v\bigr
\rangle=\langle-u \cdot\nabla z,v\rangle
\\
&\leq& C\|\nabla z\|_{L^\infty}\bigl[|v|^2+|v|\cdot |z|\bigr],
\end{eqnarray*}
which implies that there exist $\tilde{C}_0=\tilde{C}_0(R_0)>0$ and
for $P$-a.s. $\omega$, $\exists0<t_0(\omega)<1$ such that
\[
\bigl |\Lambda^\alpha\theta(t_0)\bigr |^2\leq
\tilde{C}_0.
\]

For any $\tilde{r}>0$ with $\tilde{r}-\alpha+1+\sigma_3<s+2\alpha
-1-\varepsilon$ for $\sigma_3=\frac{2}{p}$, we have the following a
priori estimate for $v$, $r=\frac{\alpha}{\alpha-{1}/{2}-
{1}/{p}}$:
%
%e6.14 #&#
%
\begin{eqnarray}
\label{eq6.14} %\begin{array}{cc}
&&\frac{d}{dt}\bigl |\Lambda^{\tilde{r}}
v\bigr |^2+2\kappa \bigl |\Lambda^{\tilde{r}+\alpha}v\bigr |^2\nonumber
\\
&&\qquad \leq 2\bigl |\bigl
\langle\Lambda^{\tilde
{r}-\alpha
}\nabla\cdot(u\theta), \Lambda^{\tilde{r}+\alpha} v\bigr
\rangle \bigr |
\\
&&\qquad \leq C\bigl |\Lambda^{\tilde{r}+\alpha} v\bigr |\cdot\bigl |\Lambda^{\tilde
{r}-\alpha
+1+\sigma_3}\theta\bigr |\cdot\|
\theta\|_{L^p}
\nonumber
\\[-1pt]
&& \qquad\leq \frac{\kappa}{4}\bigl |\Lambda^{\tilde{r}+\alpha}v\bigr |^2 +C\|\theta
\|_{L^p}^r\bigl |\Lambda^{\tilde{r}}v\bigr |^2 +C\bigl |
\Lambda^{\tilde{r}-\alpha+1+\sigma_3}z\bigr |^2\cdot\|\theta\| _{L^p}^2,
\nonumber
\end{eqnarray}
where we used $\operatorname{div}u=0$ in the first inequality and
Lemmas~\ref{lem2.1}, \ref{lem2.2} in
the second inequality and the interpolation inequality (\ref{eq2.3}) and
Young's inequality in the last inequality.
We choose the approximation $v_n$ as in the proof of Theorem~\ref
{theA.1} with
initial time $t=0$ replaced by initial time $t=t_0(\omega)$. Then by a
similar argument as in the proof of Theorem~\ref{theA.1} we have the following
$L^p$-norm estimate of $v_n$,\vspace*{-1pt}
\[
\frac{d}{dt}\|v_n\|_{L^p}^p
\leq Cp\|\nabla z\|_\infty\bigl(\|v_n\| _{L^p}^p+
\|z\|_{L^p}\|v_n\|_{L^p}^{p-1}\bigr).
\]
Thus, we have\vspace*{-2pt}
\[
\frac{d}{dt}\|v_n\|_{L^p}\leq C\|\nabla
z\|_\infty\bigl(\|v_n\| _{L^p}+\|z\|_{L^p}\bigr).
\]
Then by Gronwall's lemma and $s> 3-2\alpha$, we obtain the uniform
$L^p$-norm estimates as (\ref{eqA.6}) for $v_n$. Moreover, by (\ref
{eq6.14}) and
Gronwall's lemma, we obtain the uniform $H^{\alpha}$-norm estimates as
(\ref{eqA.7}) for $v_n$. By a similar argument as in the proof of
Theorem~\ref{theA.1},
we have $v_n$ converges to some process $\tilde{v}$ in $L^2([t_0,T],H)$
such that $\tilde{v}+z$ is the solution of (\ref{eq3.1}) in
$[t_0,T]$. Then by
the uniqueness proof in Theorem~\ref{the4.2}, we have $\tilde{v}=v$, which
implies for $P$-a.s. $\omega$, $v\in L_{\mathrm{loc}}^\infty
([t_0,\infty
),H^\alpha)\cap L^2_{\mathrm{loc}}([t_0,\infty),H^{2\alpha})$. Therefore,
(\ref{eq6.14}) also holds for $v$ with $\tilde{r}=\alpha$, which
implies that\vspace*{-1pt}
\[
\bigl |\Lambda^\alpha v(t)\bigr |^2+\kappa\int_{t_0}^t\bigl |
\Lambda^{2\alpha} v(l)\bigr |^2\,dl\leq\bigl(\bigl |\Lambda^\alpha
v(t_0)\bigr |^2+C(R_0)\bigr) \bigl(\exp
\bigl[C(R_0)t\bigr]+1\bigr),
\]
which implies that there exist $\tilde{C}_1=\tilde{C}_1(R_0)>0$ and
$\tilde{T}_0(R_0)$ such\vspace*{1pt} that\break $|\Lambda^\alpha v(\tilde{T}_0)|\leq
\tilde
{C}_1(R_0)$.
Moreover, there exists $t_1=t_1(\omega)>t_0(\omega)$ such that
$|\Lambda^{2\alpha}v(t_1)|\leq\tilde{C}_1$.
Using (\ref{eq6.14}) for $\tilde{r}=2\alpha$ and by similar
arguments as above,
we obtain that there exists $T_0=T_0(R_0)$ independent of $\omega$ such
that\break $|\Lambda^{2\alpha}v(T_0)|\leq\tilde{C}_1$.
Then we proceed analogously and obtain that there exists
$T_1=T_1(R_0)>T_0(R_0)$ such\vspace*{2pt} that $|\Lambda^{s+\delta}v(T_1)|\leq
\tilde
{C}_1$ for some $0<\delta<3\alpha-2-\sigma_3-\varepsilon$.
\end{pf}
%
%$\hfill\Box$

\begin{Lemma}\label{lem6.3.2}
Let $\alpha>2/3$. Suppose Hypothesis~\ref{hypE.3} holds with
$s>3-2\alpha$. Then for each $R\geq1$ there exist $T_1>0$ and a compact
subset $K\subset\mathcal{W}$ such that\vspace*{-1pt}
\[
\inf_{\|x\|_{L^p}\leq R}P_{T_1}(x, K)>0
\]
for $p$ in Lemma~\ref{lem6.3.1}.
\end{Lemma}

\begin{pf}
Define $K:=\{x\dvtx |\Lambda^{s+\delta} x|^2\leq\tilde{C}_1(R_0)\}
$, where $\tilde{C}_1(R_0),\delta$ comes from the previous lemma. By
Lemma~\ref{lem6.3.1}, for $R\leq R_0$, we have\vspace*{-1pt}
\begin{eqnarray*}
\inf_{\|x\|_{L^p}\leq R}P_{T_1}(x, K)&\geq&\inf
_{\|x\|_{L^p}\leq
R}\Bigl(1-P_x\Bigl[\sup_{t\in[0,T_1]}\bigl |
\Lambda^{s+2\alpha-1-\varepsilon}z(t)\bigr |^2> R_0\Bigr]
\\[-1pt]
&&\phantom{\hspace*{85pt}} {}-P_x\Bigl[\sup_{t\in[0,T_1]}
\bigl \|\theta(t)\bigr \| _{L^p}^p>R_0\Bigr]\Bigr).
\end{eqnarray*}
Under Hypothesis~\ref{hypE.3}, since $z$ is a Gaussian process, one deduces
that there exist $\eta, C>0$ such that
\[
P_x\Bigl[\sup_{t\in[0,T_1]}\bigl |\Lambda^{s+2\alpha-1-\varepsilon}z(t)\bigr |^2>
R_0\Bigr]\leq Ce^{-\eta({R_0^2}/{T_1})}
\]
(see, e.g., \cite{FR07}, Proposition~15).
Also by Theorem~\ref{the3.3}, we obtain
\[
\sup_{\|x\|_{L^p}\leq R}P_x\Bigl[\sup_{t\in[0,T_1]}
\bigl \|\theta(t)\bigr \| _{L^p}^p>R_0\Bigr]\leq\sup
_{\|x\|_{L^p}\leq R}\frac{E_x[\sup_{t\in
[0,T_1]}\|\theta(t)\|_{L^p}^p]}{R_0}\leq\frac{C(R)}{R_0}.
\]
Choosing $R_0$ big enough, we prove the assertion.
\end{pf}
%
%$\hfill\Box$

The exponential convergence now follows from Lemma~\ref{lem6.3.2} and an
abstract result of \cite{GM05}, Theorem~3.1. For $p>\frac{2}{3\alpha
-2}$ let
$V\dvtx L^p\rightarrow\mathbb{R}$ be a measurable function and define
$\|\phi
\|_{\mathrm{V}}:=\sup_{x\in L^p}\frac{|\phi(x)|}{V(x)}$ and $\|\nu\|_{\mathrm{V}}:=\sup_{\|
\phi\|_{\mathrm{V}}\leq1}\langle\nu,\phi\rangle$ for a signed measure $\nu$.

\begin{Theorem}\label{the6.3.3}
Let $\alpha>2/3$. Assume that Hypothesis~\ref{hypE.3} holds
with $s>3-2\alpha$ and let $V(x):=1+\|x\|_{L^p}^p$ for $p>\frac
{2}{3\alpha-2}$. Then there exist $C_{\exp}>0$ and $a>0$ such that
\[
\bigl \|P_t^*\delta_{x_0}-\mu\bigr \|_{\mathrm{var}}\leq
\bigl \|P_t^*\delta_{x_0}-\mu \bigr \| _{\mathrm{V}}\leq
C_{\exp}\bigl(1+\|x_0\|_{L^p}^p
\bigr)e^{-at}
\]
for all $t>0$ and $x_0\in L^p$, where $\|\cdot\|_{\mathrm{var}}$ is the
total variation distance on measures.
\end{Theorem}

\begin{pf}
By \cite{GM05}, Theorem~3.1, we need to verify the following four conditions:
\begin{enumerate}[3.]
\item[1.]the measures $(P_t(x,\cdot))_{t>0,x\in L^p}$ are equivalent,
\item[2.]$x\rightarrow P_t(x,\Gamma)$ is continuous in $\mathcal{W}$
for all
$t>0$ and all Borel sets $\Gamma\subset H$,
\item[3.]for each $R\geq1$ there exist $T_1>0$ and a compact subset
$K\subset
\mathcal{W}$ such that
\[
\inf_{\|x\|_{L^p}\leq R}P_{T_1}(x, K)>0,
\]
\item[4.]there exist $k,b,c>0$ such that for all $t\geq0$,
\[
E^{P_x}\bigl[\bigl \|\theta(t)\bigr \|_{L^p}^p\bigr]\leq k
\|x\|_{L^p}^pe^{-bt}+c.
\]
\end{enumerate}

Condition 1 can be verified by \cite{GM05}, Lemma~3.2, and
$P_t(x,\mathcal
{W})=1$ for $x\in L^p$ since for fixed $t>0$ the solution $\theta$ will
go into $H^s$ space if the initial value $x\in L^p$.
Other conditions can be verified by Theorem~\ref{the6.2.4}, Lemma~\ref
{lem6.3.2} and
Proposition~\ref{pro5.1.5}.
\end{pf}
%
%$\hfill\Box$

\begin{Remark}\label{rem6.3.4}
For $\alpha>\frac{3}{4}$, we could get a better result
following a similar argument as in \cite{R08}. Namely, there exist
$C_{\exp
}>0$ and $a>0$ such that
\[
\bigl \|P_t^*\delta_{x_0}-\mu\bigr \|_{\mathrm{TV}}\leq
\bigl \|P_t^*\delta_{x_0}-\mu\bigr \| _{\mathrm{V}}\leq
C_{\exp}\bigl(1+|x_0|^2\bigr)e^{-at}
\]
for all $t>0$ and $x_0\in H$. Here, $P_t$ could be every Markov
selection obtained in Theorem~\ref{theC.5} associated to the solution of
equation (\ref{eq3.1}). The reason why $\alpha>\frac{3}{4}$ is
needed is as follows.

As in Theorem~\ref{the6.1.2}, we can prove $P_t$ is $H^s$-strong
Feller with
$s>3-3\alpha$. And for a solution $\theta$ of equation (\ref{eq3.1}) starting
from $x\in H$, we can only prove that it will enter $H^\alpha$ under
Hypothesis~\ref{hypE.3}. If the process $\theta$ enters $H^s$, we can
prove that it
satisfies the above four conditions. Hence, to obtain exponential
convergence for every $x\in H$, we need the process starting from $x\in
H$ to enter $H^s$. Hence, we need $3-3\alpha<s\leq\alpha$, that is,
$\alpha>\frac{3}{4}$.
\end{Remark}

\begin{appendix}
\section*{Appendix A}\label{appA}

In this appendix, we construct a measurable map associated with the
stochastic quasi-geostrophic equation, which will be used in the proof
of Section~\ref{sec6}.
This proof is similar as done in \cite{ZZ12}, Section~3. Here, we give it
for the reader's convenience.

Assume that for any $m< 2+\sigma$,
$z\in C((0,\infty), H^m)$ with $\sigma$ in Hypothesis~\ref{hypE.1}. Then
consider the following equation:
\setcounter{equation}{0}
\begin{equation}
\label{eqA.1} \frac{dv}{dt}+A_\alpha v+(u_v+u_z)
\cdot\nabla(v+z)=0.
\end{equation}
For (\ref{eqA.1}), we obtain the following existence and uniqueness
result if
the initial value starts from $H^1$.

\begin{theorem}\label{theA.1}
Fix $\alpha>1/2$. Suppose that for any $m< 2+\sigma$,
$z\in C((0,\infty),\allowbreak   H^m)$. For any $v_0\in H^1$, there exists a unique
solution $v\in L^\infty_{\mathrm{loc}}([0,\infty);\break H^1)
\cap L^2_{\mathrm{loc}}([0,\infty); H^{1+\alpha})$ of equation
(\ref{eqA.1}) with $v(0)=v_0$, that is, for any $\varphi\in
C^1(\mathbb{T}^2)$
\begin{eqnarray*}
&&\bigl\langle v(t),\varphi\bigr\rangle-\langle v_0,\varphi\rangle
\\
&&\qquad{}+
\int_{0}^t\bigl\langle A_\alpha^{1/2}v(r),A_\alpha^{1/2}
\varphi\bigr\rangle\, dr-\int_{0}^t \bigl
\langle(u_v+u_z) (r)\cdot\nabla\varphi,(v+z) (r)\bigr
\rangle\,dr=0,
\end{eqnarray*}
where $u_v,u_z$ satisfy (\ref{eq1.3}) with $\theta$ replaced by
$v,z$, respectively.
\end{theorem}

\begin{pf}
We construct an approximation of (\ref{eqA.1}) by a similar
construction as in the proof of Theorem~\ref{the3.3}.

We pick a smooth $\phi\geq0$, with $\operatorname{supp}\phi\subset[1,2], \int_0^\infty
\phi=1$, and for $\delta>0$ let
\[
U_\delta[\theta](t):=\int_0^\infty\phi(
\tau) \bigl(k_\delta*R^\bot \theta \bigr) (t-\delta\tau)\,d
\tau,
\]
where $k_\delta$ is the periodic Poisson kernel in $\mathbb{T}^2$
given by
$\widehat{k_\delta}(\zeta)=e^{-\delta|\zeta|},\zeta\in\mathbb{Z}^2$,
and we set $\theta(t)=0, t<0$. We take a zero sequence $\delta_n,
n\in
\mathbb{N}$, and consider the equation
\begin{equation}
\label{eqA.2} d v_n(t)+A_\alpha v_n(t)
\,dt+u_n(t)\cdot\nabla\bigl(v_n(t)+z\bigr)\,dt= 0,
\end{equation}
with initial data $v_n(0)=v_0$ and $u_n=U_{\delta_n}[v_n+z]$.
For a fixed $n$, this is a linear equation in $v_n$ on each subinterval
$[t_k,t_{k+1}]$ with $t_k=k\delta_n$, since $u_n$ is determined by the
values of $v_n$ on the two previous subintervals. By a similar argument
as in the proof of Theorem~\ref{the3.3}, we obtain the existence and uniqueness
of a solution $v_n\in C([0,T],H^{1})\cap L^2([0,T],H^{1+\alpha})$ to
(\ref{eqA.2}) (for more details, we refer to \cite{ZZ12}). Now we
take any $p$
satisfying $\frac{2}{2\alpha-1}<p<\infty$. From now on, we fix such $p$
and we have $H^1\subset L^p $ by Lemma~\ref{lem2.2}. Since the
periodic Riesz
transform is bounded on $L^p$, we have for $t>0$
\begin{equation}
\label{eqA.3} \sup_{[0,t]}\bigl \|U_\delta[\theta]
\bigr \|_{L^p}\leq C\sup_{[0,t]}\|\theta\| _{L^p},
\end{equation}
and also
\begin{equation}
\label{eqA.4} \int_{0}^t\bigl \|U_\delta[
\theta]\bigr \|_{L^p}^p\,d\tau\leq C\int_{0}^t
\| \theta\| _{L^p}^p\,d\tau.
\end{equation}

By Lemma~\ref{lem5.1.4}, we obtain for $v_n$ the following inequality
by taking
inner product with $|v_n|^{p-2}v_n$ in $L^2$:
\begin{eqnarray}
\label{eqA.5} \frac{d}{dt}\|v_n\|_{L^p}^p+2
\lambda_1\|v_n\|_{L^p}^p&\leq& p\bigl |
\bigl\langle u_n\cdot\nabla(v_n+z), |v_n|^{p-2}v_n
\bigr\rangle\bigr |
\nonumber
\\[-8pt]
\\[-8pt]
&\leq&p\| \nabla z\|_\infty\|u_n\|_{L^p}
\|v_n\|_{L^p}^{p-1},
\nonumber
\end{eqnarray}
where we used $\operatorname{div}u_n=0$ and H\"{o}lder's inequality in the
last inequality.
Therefore,
\begin{eqnarray*}
&&\bigl \|v_n(t)\bigr \|_{L^p}^p-\bigl \|v_n(0)
\bigr \|_{L^p}^p+\int_{0}^t2
\lambda_1 \bigl \| v_n(\tau)\bigr \|_{L^p}^p\,d
\tau
\\
&&\qquad\leq\varepsilon\int_{0}^t\bigl(
\|u_n\|_{L^p}^p+\|v_n
\|_{L^p}^p\bigr)\, d\tau +pC(\varepsilon)\int
_{0}^t\|\nabla z\|_\infty^{{p}/{(p-1)}}
\|v_n\| _{L^p}^{p}\,d\tau
\\
&&\qquad\leq\varepsilon\int_{0}^t
\|v_n\|_{L^p}^p\,d\tau +pC(\varepsilon)\int
_{0}^t\|\nabla z\|_\infty^{{p}/{(p-1)}}
\|v_n\|_{L^p}^{p}\,d\tau +C\int
_0^t\|z\|_{L^p}^p\,d\tau,
\end{eqnarray*}
where we used (\ref{eqA.4}) in the last inequality.
Then Gronwall's lemma and $H^1\subset L^p $ yield that for any $T\geq0$
\begin{equation}
\label{eqA.6} \sup_{t\in[0,T]}\bigl \|v_n(t)
\bigr \|_{L^p}\leq C,
\end{equation}
where $C$ is a constant independent of $n$.

Moreover, we get the following estimate by taking the inner product in
$L^2$ with $\Lambda e_k$ for (\ref{eqA.2}), multiplying both sides by
$\langle
v,\Lambda e_k\rangle$ and summing up over $k$:
\begin{eqnarray*}
&&\frac{1}{2}\frac{d}{dt}|\Lambda v_n|^2+
\kappa\bigl |\Lambda^{1+\alpha
} v_n\bigr |^2
\\
&&\qquad \leq \bigl |
\Lambda^{1-\alpha}\bigl(u_n\cdot\nabla(v_n+z)
\bigr)\bigr |\bigl |\Lambda ^{1+\alpha}v_n\bigr |
\\
&&\qquad \leq C\bigl |\Lambda^{1+\alpha} v_n\bigr |\bigl[\bigl |\Lambda^{2-\alpha
+\sigma_1}(v_n+z)\bigr |
\|u_n\|_{L^p}+\bigl |\Lambda^{2-\alpha+\sigma_1}u_n\bigr |\|
v_n+z\|_{L^p}\bigr],
\end{eqnarray*}
where $\sigma_1=2/p<(2\alpha-1)$ and we used Lemma~\ref{lem2.1} in
the last inequality.
Hence, we obtain that for $r=\frac{2\alpha}{2\alpha-1-\sigma_1}$,
\begin{eqnarray*}
&&\frac{1}{2}\bigl(\bigl |\Lambda v_n(t)\bigr |^2-\bigl |\Lambda
v_n(0)\bigr |^2\bigr)+\kappa\int_{0}^t\bigl |
\Lambda^{1+\alpha} v_n\bigr |^2\,d\tau
\\
&&\qquad\leq C\int_{0}^t\bigl |\Lambda
^{1+\alpha} v_n\bigr |\bigl[\bigl |\Lambda^{2-\alpha+\sigma_1}(v_n+z)\bigr |
\|u_n\| _{L^p}+\bigl |\Lambda^{2-\alpha+\sigma_1}u_n\bigr |
\|v_n+z\|_{L^p}\bigr]\,d\tau
\\
&&\qquad\leq\frac{\kappa}{2}\int_{0}^t\bigl |
\Lambda^{1+\alpha} v_n\bigr |^2d\tau
\\
&&\qquad\quad{}+C\Bigl[\sup
_{t\in[0,T]}\bigl(\bigl \|v_n(t)+z(t)\bigr \|_{L^p}^r+
\bigl \|v_n(t)+z(t)\bigr \|^2_{L^p}\bigr)+1\Bigr]
\\
&&\qquad\quad\quad{}\times\int
_{0}^t|\Lambda v_n|^2+
\bigl |\Lambda^{2-\alpha+\sigma_1}z\bigr |^2\,d\tau,
\end{eqnarray*}
where we used (\ref{eqA.3}), (\ref{eqA.4}), the interpolation
inequality (\ref{eq2.3}) and
Young's inequality in the last inequality.
By Gronwall's lemma and (\ref{eqA.6}), we get that for $v_0\in H^1$
\begin{equation}
\label{eqA.7} \sup_{0\leq t\leq T}\bigl |\Lambda v_n(t)\bigr |^2+
\kappa\int_{0}^T\bigl |\Lambda ^{1+\alpha}
v_n\bigr |^2\,d\tau\leq C,
\end{equation}
where $C$ is a constant independent of $n$. Now decompose $v_n$ as
\[
v_n(t)=v_0-\int_0^tA_\alpha
v_n(s)\,ds-\int_0^t
\bigl(u_n(s)\cdot\nabla \bigl(v_n(s)+z(s)\bigr)\bigr)
\,ds.
\]
By (\ref{eqA.7}), we obtain
\[
\biggl\|\int_0^\cdot A_\alpha
v_n(s)\,ds\biggr\|_{W^{1,2}(0,T,H^{-\alpha})}\leq C
\]
and
\[
\biggl\|\int_0^\cdot\bigl(u_n(s)\cdot
\nabla\bigl(v_n(s)+z(s)\bigr)\bigr)\,ds\biggr\| _{W^{1,2}(0,T,H^{-3})}\leq C.
\]
So, we have proved
\[
\|v_n\|_{W^{1, 2}([0,T], H^{-3})}\leq C,
\]
where $C$ is a constant independent of $n$.
By the compactness embedding $W^{1,2}([0,T],H^{-3})\cap
L^2([0,T],H^{1+\alpha})\subset L^2([0,T],H^1)$ we have that there
exists a subsequence of $v_n$ converging in $L^2([0,T],H^1)$ to a
solution $v\in\break L^\infty_{\mathrm{loc}}([0,\infty);H^1)
\cap L^2_{\mathrm{loc}}([0,\infty); H^{1+\alpha})$ of equation (\ref
{eqA.1}). Thus, (\ref{eqA.7}) is
also satisfied for $v$. Uniqueness can be deduced from a similar
argument as in the proof of Theorem~\ref{the4.2}.
\end{pf}
%
%$\hfill\Box$

\begin{theorem}\label{theA.2}
Fix $\alpha>1/2$. Suppose that for any $m< 2+\sigma$,
$z\in C([0,\infty),\allowbreak   H^m)$. The solution $v$ obtained in Theorem~\ref
{theA.1} is
in $ C([0,\infty);H^1 )$.
\end{theorem}

\begin{pf} It is sufficient to show that
\[
\Lambda\frac{dv}{dt}\in L^2_{\mathrm{loc}}\bigl([0,\infty
); H^{-\alpha}\bigr).
\]
For $\varphi$ smooth enough, we have
\begin{eqnarray*}
\biggl|\biggl\langle\frac{dv}{dt},\Lambda\varphi\biggr\rangle\biggr|&=&|
\kappa\bigl\langle -\Lambda^{\alpha}v,\Lambda^{1+\alpha}\varphi\bigr
\rangle-\bigl\langle \bigl(u\cdot \nabla(\Lambda\varphi)\bigr),v+z\bigr\rangle
\\
&\leq&\bigl[\kappa\bigl |\Lambda^{1+\alpha
}v\bigr |+C\bigl |\Lambda^{2-\alpha}\bigl(u
\cdot(v+z)\bigr)\bigr |\bigr]\bigl |\Lambda^\alpha\varphi\bigr |
\\
&\leq&C\bigl[\bigl |\Lambda^{1+\alpha}v\bigr |+\bigl |\Lambda^{2-\alpha+\sigma_1} (v+z)\bigr |\|v+z
\|_{L^p}\bigr]\bigl |\Lambda^\alpha\varphi\bigr |, %\end{array}
\end{eqnarray*}
where $0<\sigma_1<2\alpha-1, p=\frac{2}{\sigma_1}$ and we used
Lemma~\ref{lem2.1} in the last inequality.
Then
\[
\biggl\|\Lambda\frac{dv}{dt}\biggr\|_{H^{-\alpha}}\leq C\bigl(\|v+z\|
_{L^p}+1\bigr)\bigl |\Lambda ^{1+\alpha}v\bigr |+C\|v+z\|_{L^p}\bigl |
\Lambda^{2-\alpha+\sigma_1}z\bigr |.
\]
By (\ref{eqA.6}) and (\ref{eqA.7}), we obtain for $0<T<\infty$
\[
\int_{0}^T\biggl\|\Lambda\frac{dv}{dt}(\tau)
\biggr\|_{H^{-\alpha}}^2\,d\tau <\infty,
\]
which implies that $v\in C([0,\infty);H^1 )$.
\end{pf}
%
%$\hfill\Box$

\begin{theorem}\label{theA.3}
Fix $\alpha>1/2$. Suppose that for any $m< 2+\sigma$,
$z\in C([0,\infty),\allowbreak   H^m)$. For any fixed $t>0$, the map $v_0\mapsto
v(t, v_0)$ is a continuous map from $H^1$ into itself, where $v(t,v_0)$
is the solution of equation (\ref{eqA.1}) with $v(0)=v_0$.
\end{theorem}

\begin{pf}
Let $v_1,v_2$ be two solutions of (\ref{eqA.1}) and $\zeta=v_1-v_2,
\theta_1=v_1+z,\theta_2=v_2+z$. Then $\zeta$ satisfies the following
equation:
\[
\biggl(\frac{d}{dt}\zeta,\varphi\biggr)+\kappa\bigl(\Lambda^\alpha
\zeta,\Lambda ^\alpha \varphi\bigr)=-(u_1\cdot\nabla\zeta,
\varphi)-(u_\zeta\cdot\nabla \theta _2,\varphi),
\]
where $u_1,u_\zeta$ satisfy (\ref{eq1.3}) with $\theta$ replaced by
$\theta
_1,\zeta$, respectively.

Taking $\varphi=\Lambda e_k$, multiplying both sides by $\langle\zeta
,\Lambda e_k\rangle$ and summing up over $k$ we have the following
estimate since $v_i\in C([0,\infty);H^1 )\cap L^2_{\mathrm
{loc}}([0,\infty
); H^{1+\alpha})$, $i=1,2$, by Theorems~\ref{theA.1} and~\ref{theA.2}:
\begin{eqnarray*}
&&\frac{1}{2}\frac{d}{dt}|\Lambda\zeta|^2+\kappa\bigl |\Lambda
^{1+\alpha} \zeta\bigr |^2
\\
&&\qquad  = -\bigl\langle\Lambda(u_1\cdot
\nabla\zeta ),\Lambda \zeta\bigr\rangle-\bigl\langle u_\zeta\cdot\nabla
\theta_2,\Lambda^{2}\zeta \bigr\rangle
\\
 &&\qquad\leq C\bigl |\Lambda^{1+\alpha} \zeta\bigr |\bigl[\bigl |\Lambda^{2-\alpha}(u_\zeta
\theta _2)\bigr |+\bigl |\Lambda^{2-\alpha}(u_1\zeta)\bigr |\bigr]
\\
&&\qquad \leq C\bigl |\Lambda^{1+\alpha} \zeta \bigr |\bigl[\bigl |\Lambda^{2-\alpha+\sigma_1}\zeta\bigr |\bigl |
\Lambda^{\sigma_2}\theta_2\bigr | +\bigl |\Lambda^{2-\alpha+\sigma_1}
\theta_2\bigr |\bigl |\Lambda^{\sigma_2}\zeta\bigr |
\\
&&\quad\phantom{\hspace*{70pt}}{}+\bigl |\Lambda^{2-\alpha+\sigma_1}\theta_1\bigr |\bigl |\Lambda^{\sigma
_2}\zeta\bigr |
+\bigl |\Lambda^{2-\alpha+\sigma_1}\zeta\bigr |\bigl |\Lambda^{\sigma_2}\theta_1\bigr |
\bigr]
\\
&& \qquad\leq \frac{\kappa}{2}\bigl |\Lambda^{1+\alpha} \zeta\bigr |^2
 \\
&&\quad\qquad {}+C\bigl[|
\Lambda \theta _2|^r+|\Lambda\theta_1|^r
%  {}
  +\bigl |\Lambda^{1+\alpha}v_2\bigr |^2+\bigl |\Lambda
^{2-\alpha+\sigma_1}z\bigr |^2 +\bigl |\Lambda^{s+\alpha}v_1\bigr |^2
\bigr]|\Lambda\zeta|^2,
\end{eqnarray*}
where $r=\frac{2\alpha}{2\alpha-1-\sigma_1}, \sigma_2=1-\sigma_1$ for
some $0<\sigma_1<(2\alpha-1)$ and we used Lemma~\ref{lem2.1} in the second
inequality and Lemma~\ref{lem2.2}, the interpolation inequality (\ref{eq2.3}),
$H^{1}\subset H^{\sigma_2} $ and Young's inequality in the last inequality.
Then Gronwall's lemma yields that\vspace*{-1pt}
\begin{eqnarray*}
\hspace*{-5pt} &&|\Lambda\zeta|^2\leq C\bigl |\Lambda\zeta(0)\bigr |^2\exp\biggl\{
\int_{0}^T\bigl |\Lambda \theta_2(
\tau)\bigr |^r+\bigl |\Lambda\theta_1(\tau)\bigr |^r
\\
\hspace*{-5pt}&&\phantom{\hspace*{125pt}}{}+\bigl |
\Lambda^{2-\alpha
+\sigma
}z\bigr |^2+\bigl |\Lambda^{1+\alpha}v_1(
\tau)\bigr |^2+ \bigl |\Lambda^{1+\alpha}v_2(\tau)\bigr |^2
\,d\tau\biggr\}.
\end{eqnarray*}
Thus, the result follows.\vadjust{\goodbreak}
\end{pf}
%
%$\hfill\Box$

Now for $v_0\in H^1, \overline{W}\in C(\mathbb{R}^+,H^{-1-\varepsilon
_0}) $ we define\vspace*{-1pt}
\[
v(t,\overline{W},v_0):=\cases{ %
v
\bigl(t,v_0,z(\overline{W})\bigr),& \quad$\mbox{if } z(\overline {W})
\in C\bigl(\mathbb{R}^+,H^{m}\bigr) \mbox{ for } m<2+\sigma$,
\cr
0,&
\quad$\mbox{otherwise}$, %\end{array}
}
\]
where $v(t,v_0,z(\overline{W}))$ is the solution to (\ref{eqA.1}) we obtained
in Theorem~\ref{theA.1}.

Combining Theorems~\ref{theA.1}--\ref{theA.3} we obtain the following results.

\begin{theorem}\label{theA.4}
Fix $\alpha>1/2$.
$v\dvtx\mathbb{R}^+\times C(\mathbb{R}^+,H^{-1-\varepsilon
_0})\times
H^1\mapsto H^1$, $(t,\overline{W},\allowbreak  v_0)\mapsto v(t,\overline{W},v_0)$ is
a measurable map.
\end{theorem}

\begin{pf}
By Theorems~\ref{theA.1}--\ref{theA.3}
$t\mapsto v(t,\overline{W},v_0)$ and
$v_0\mapsto v(t,\overline{W},v_0)$ is continuous. Then it is sufficient
to prove that if $z_n\rightarrow z$ in $C(\mathbb{R}^+,H^{m}),
m<2+\sigma$, $v_n\rightarrow v$ in $C([0,T],H^1)$, where $v_n=v(\cdot
,v_0,z_n), v=v(\cdot,v_0,z)$. By the same arguments as in the proof of
Theorem~\ref{theA.1}, we have the following estimate:\vspace*{-1pt}
\[
\sup_{[0,T]}|\Lambda v_n|^2\leq C(T),
\qquad\sup_{[0,T]}|\Lambda v|^2\leq C(T),
\]
and\vspace*{-1pt}
\[
\int_0^{T}\bigl |\Lambda^{1+\alpha}v_n(l)\bigr |^2
\,dl\leq C(T),\qquad\int_0^{T}\bigl |\Lambda
^{1+\alpha}v(l)\bigr |^2\,dl\leq C(T).
\]
Since $v,v_n\in C([0,+\infty),H^1)\cap L^2_{\mathrm{loc}}((0,+\infty
),H^{1+\alpha})$, we obtain\vspace*{-1pt}
\begin{eqnarray*}
&&\frac{d}{dt}\bigl |\Lambda(v-v_n)\bigr |^2+2
\kappa\bigl |\Lambda^{1+\alpha} (v_n-v)\bigr |^2
\\[-1pt]
&&\qquad=\bigl\langle(u_{v_n}+u_{z_n})\cdot
\nabla(v_n+z_n), \Lambda ^{2}(v-v_n)
\bigr\rangle
\\[-1pt]
&&\qquad\quad{}-\bigl\langle(u_v+u_z)\cdot\nabla(v+z),
\Lambda ^{2}(v-v_n)\bigr\rangle
\\
&&\qquad=\bigl[\bigl\langle(u_{v_n}-u_v)\cdot
\nabla(v_n+z_n),\Lambda ^{2}(v-v_n)
\bigr\rangle
\\
&&\qquad\quad{}+\bigl\langle(u_v+u_z)\cdot
\nabla(v_n-v),\Lambda ^{2}(v-v_n)\bigr\rangle
\\
&&\qquad\quad{}+\bigl\langle(u_{z_n}-u_z)\cdot
\nabla (v_n+z_n),\Lambda ^{2}(v-v_n)
\bigr\rangle
\\
&&\qquad\quad\phantom{\hspace*{9pt}}{}+\bigl\langle(u_v+u_z)\cdot
\nabla(z_n-z),\Lambda ^{2}(v-v_n)\bigr\rangle
\bigr], %\end{array}
\end{eqnarray*}
where $u_{v_n},u_{z_n}$ satisfy (\ref{eq1.3}) with $\theta$ replaced by
$v_n,z_n$ respectively.
For the first term on the right-hand side, we have
\begin{eqnarray*}
&&\bigl |\bigl\langle(u_{v_n}-u_v)\cdot
\nabla(v_n+z_n),\Lambda ^{2}(v-v_n)
\bigr\rangle\bigr |
\\
&&\qquad = \bigl |\bigl\langle\Lambda^{1-\alpha}\nabla\cdot
\bigl((u_{v_n}-u_v) (v_n+z_n)
\bigr),\Lambda^{1+\alpha}(v-v_n)\bigr\rangle\bigr |
\\
&&\qquad \leq  C\bigl |\Lambda^{1+\alpha}(v-v_n)\bigr |\bigl |\Lambda^{{2-\alpha}+\sigma
_1}(v-v_n)\bigr |\bigl |
\Lambda^{\sigma_2}(v_n+z_n)\bigr |
\\
&&\qquad\quad  {}+C\bigl |\Lambda^{1+\alpha}(v-v_n)\bigr |\bigl |\Lambda^{{2-\alpha}+\sigma
_1}(v_n+z_n)\bigr |\bigl |
\Lambda^{\sigma_2}(v-v_n)\bigr |
\\
&&\qquad \leq  \frac{\kappa}{4}|\Lambda ^{1+\alpha}(v-v_n)|^2+C
\bigl(C(T)+\bigl |\Lambda^{1+\alpha}v_n\bigr |^2\bigr)\bigl |\Lambda
(v-v_n)\bigr |^2
\\
&&\qquad\quad  {}+c\bigl |\Lambda^{{2-\alpha}+\sigma_1}z_n\bigr |^2\bigl |
\Lambda(v-v_n)\bigr |^2. %\end{array}
\end{eqnarray*}
Here, $\sigma_1,\sigma_2$ are as (\ref{eq6.6}) and we used
$\operatorname{div}
(u_{v_n}-u_v)=0$ in the first equality and Lemmas~\ref{lem2.1}
and~\ref{lem2.2} in the
first inequality and the interpolation inequality (\ref{eq2.3}) and Young's
inequality in the last step. The other term can be estimated similarly.
Then we obtain
\begin{eqnarray*}
&&\frac{d}{dt}\bigl |\Lambda(v-v_n)\bigr |^2+2
\kappa\bigl |\Lambda^{1+\alpha} (v_n-v)\bigr |^2
\\
&&\qquad\leq\kappa\bigl |\Lambda^{1+\alpha} (v_n-v)\bigr |^2+C
\bigl(C(T)+\bigl |\Lambda ^{1+\alpha} v_n\bigr |^2+\bigl |
\Lambda^{1+\alpha} v\bigr |^2\bigr)
\\
&&\phantom{\hspace*{123pt}}{}\times \bigl(\bigl |\Lambda (v-v_n)\bigr |^2+\bigl |
\Lambda ^{2-\alpha+\sigma_1}(z-z_n)\bigr |^2\bigr). %\end{array}
\end{eqnarray*}
Gronwall's lemma yields that
\begin{eqnarray*}
\bigl |\Lambda(v-v_n) (t)\bigr |^2&\leq&\Theta_n\exp
\biggl(C\int_0^t\bigl(C(T)+\bigl |\Lambda
^{1+\alpha
} v_n\bigr |^2+\bigl |\Lambda^{1+\alpha}
v\bigr |^2\bigr)\,dl\biggr)
\\
&&{}\times\int_0^t
\bigl(C(T)+\bigl |\Lambda ^{1+\alpha} v_n\bigr |^2+\bigl |
\Lambda^{1+\alpha} v\bigr |^2\bigr)\,dl,
\end{eqnarray*}
where $\Theta_n=\sup_{[0,T]}|\Lambda^{2-\alpha+\sigma_1}(z-z_n)|$. Then
the results follow.
\end{pf}
%
%$\hfill\Box$

\section*{Appendix B}\label{appB}

In this appendix, we prove the following lemma
to complete the proof of Theorem~\ref{the3.3}.

\renewcommand{\thelemma}{B.\arabic{lemma}}
\begin{lemma}\label{lemB.1}
For any $x_0\in B_0$ defined in the proof of Theorem~\ref{the3.3},
there exists $Q_{x_0}\in\mathcal{Q}_{x_0}$ such that the map
$x_0\mapsto Q_{x_0}$ from $B_0$ to $\mathcal{P}(\Omega_0^{t_1^n})$ is
measurable with respect to $\mathcal{B}_{t_1^n}$.
\end{lemma}

\begin{pf}
Let $\mathcal{B}_{t_1^n}^1$ be the Borel $\sigma$-algebra on
$\tilde{B}_0:=\{x(\cdot)1_{[0,t_1^n]}(\cdot
)+x(t_1^n)\* 1_{[t_1^n,\infty
)}(\cdot)\dvtx x\in B_0\}$ with  the topology induced by $\sup_{0\leq
t\leq
t_1^n}\|x(t)\|_{H^3}$. Since $\{\sup_{0\leq t\leq t_1^n}\|x(t)\|
_{H^3}<a\}\in\mathcal{B}_{t_1^n}$, we know $\mathcal
{B}_{t_1^n}^1\subset\mathcal{B}_{t_1^n}$.
It suffices to prove that if for $\{x_m,m\in\mathbb{N}\cup\{0\}\}
\subset\tilde{B}_0$, $\sup_{0\leq t\leq t_1^n}\|x_m(t)-x_0(t)\|
_{H^3}\rightarrow0$ and $Q_m\in\mathcal{Q}_{x_m}$, then for some
subsequence $m_k$, $Q_{m_k}$ weakly converges to some $Q\in\mathcal
{Q}_{x_0}$, because then \cite{SV79}, Lemma~12.1.8, Theorem~12.1.10, implies
the existence of a $Q_{x_0}\in\mathcal{Q}_{x_0(\cdot
)1_{[0,t_1^n]}(\cdot
)+x_0(t_1^n)1_{[t_1^n,\infty)}(\cdot)}$ such that the map $x_0\mapsto
x_0(\cdot)1_{[0,t_1^n]}(\cdot)+x(t_1^n)1_{[t_1^n,\infty)}(\cdot
)\mapsto
Q_{x_0}$ from $B_0$ to $\tilde{B}_0$ to$\mathcal{P}(\Omega_0^{t_1^n})$
is measurable with respect to $\mathcal{B}_{t_1^n}$. Moreover, by
$Q_{x_0}\in\mathcal{Q}_{x_0}$, the result follows.

\textit{Step} 1: We prove that $(Q_m)_{m\in\mathbb{N}}$ is tight in
$\mathbb
{S}:=C([t_1^n,+\infty),H^1)\cap L^q_{\mathrm{loc}}([t_1^n,\allowbreak  +\infty),H^3)$
for some $q\in\mathbb{N}$. Define for each $m\in\mathbb{N}$,
\[
M^m(t,x):=\sum_{i=1}^\infty
M_i^m(t,x)e_i,
\]
where $M_i^m$ is given in the proof of Theorem~\ref{the3.3} (Step 2)
with $x_0$
replaced by $x_m$. Then $(M^m(t,x))_{t\geq t_1^n}$ is a continuous
$H^3$-valued $\mathcal{B}_t$-martingale with respect to $Q_m$ and the
following equality holds in $H^1$:
\setcounter{equation}{0}
\renewcommand{\theequation}{B.\arabic{equation}}
\begin{eqnarray}
\label{eqB.1} x(t)&=&x_m\bigl(t_1^n\bigr)-
\int_{t_1^n}^{t\wedge t_2^n} \bigl(A_\alpha
x(s)+U_{\delta
_n}[x_m](s)\cdot\nabla x(s)\bigr)
\,ds
\nonumber
\\[-8pt]
\\[-8pt]
&&{}+M^m(t),\qquad Q_m\mbox{-a.s.}
\nonumber
\end{eqnarray}
By H\"{o}lder's inequality and (M3), (M1) for $Q_m$, we have
\begin{eqnarray}
\label{eqB.2} &&E^{Q_m}\biggl[\sup_{s\neq t\in[t_1^n,t_2^n]}\biggl(\biggl\|
\int_s^t A_\alpha x(r)+U_{\delta_n}[x_m](r)
\cdot\nabla x(r)\,dr\biggr\|_{H^1}^\gamma \Big/|t-s|^{\gamma-1}
\biggr)\biggr]
\nonumber
\\
&&\qquad\leq CE^{Q_m}\biggl[\int_{t_1^n}^{t_2^n}
\bigl \|A_\alpha x(r)+U_{\delta
_n}[x_m](r)\cdot\nabla x(r)
\bigr \|_{H^1}^\gamma\,dr\biggr]
\nonumber
\\[-8pt]
\\[-8pt]
&&\qquad\leq CE^{Q_m}\Bigl[\sup_{t_1^n\leq r\leq t_2^n}\bigl \| x(r)
\bigr \|_{H^3}^\gamma\Bigl(1+\sup_{0\leq r\leq
t_1^n}
\bigl \|x_m(r)\bigr \|_{H^1}^\gamma\Bigr) \Bigr]
\nonumber
\\
&&\qquad\leq C\bigl(\bigl \| x_m\bigl(t_1^n
\bigr)\bigr \|_{H^3}^\gamma +1\bigr) \Bigl(1+\sup
_{0\leq r\leq t_1^n}\bigl \|x_m(r)\bigr \|_{H^1}^\gamma
\Bigr),
\nonumber
\end{eqnarray}
where $C$ is independent of $m$.
For $t_1^n\leq s<t\leq t_2^n$ and $q\in\mathbb{N}$, we have
\begin{eqnarray*}
E^{Q_m}\bigl \|M^m(t,x)-M^m(s,x)
\bigr \|_{H^1}^{2q}&\leq& C_q E^{Q_m}\biggl(
\int_s^t\bigl \|\Lambda\bigl(k_{\delta_n}*G
\bigl(x(r)\bigr)\bigr)\bigr \|_{L_2(U;H)}^2\,dr\biggr)^q
\\
&\leq& C_q |t-s|^{q-1} \int_s^t
E^{Q_m}\bigl \|G\bigl(x(r)\bigr)\bigr \|_{L_2(U;H)}^{2q}\,dr
\\
&\leq&C_q |t-s|^{q-1}\int_s^t
E^{Q_m}\bigl(\bigl |x(r)\bigr |^{2q}+1\bigr)\,dr
\\
&\leq&C_q |t-s|^{q} \bigl(\bigl |\Lambda^3x_m
\bigl(t_1^n\bigr)\bigr |^{2q}+1\bigr), %\end{array}
\end{eqnarray*}
where we used Hypothesis~\ref{hypG.1} in the third inequality and (M3)
in the last
inequality. By Kolmogorov's criterion for any $\beta\in(0,\frac
{q-1}{2q})$, we get
\begin{equation}
\label{eqB.3} \qquad E^{Q_m}\biggl(\sup_{s\neq t\in[t_1^n,t_2^n]}
\frac{\|M^m(t,x)-M^m(s,x)\|
_{H^1}^{2q}}{|t-s|^{q\beta}}\biggr)\leq C\bigl(\bigl |\Lambda^3x_m
\bigl(t_1^n\bigr)\bigr |^{2q}+1\bigr).
\end{equation}
Combining (\ref{eqB.1})--(\ref{eqB.3}) and $Q_m(\{x\dvtx
x(s)=x(t_2^n),s\in[t_2^n,+\infty)\}
)=1$, we obtain for $\beta_1=1-\frac{1}{\gamma}$ and any $T>0$
\[
\sup_{m\in\mathbb{N}}E^{Q_m}\biggl(\sup_{s\neq t\in[t_1^n,T]}
\frac{\|
x(t)-x(s)\|_{H^1}}{|t-s|^{\beta_1}}\biggr)<\infty.
\]
Thus, by (M3) for $Q_m$ and \cite{GRZ09}, Lemma~4.3, $(Q_m)_{m\in
\mathbb{N}}$
is tight in $\mathbb{S}$.

Without loss of generality, we assume that $Q_m$ weakly converges to
some probability measure $Q$ in $\mathbb{S}$. We need to prove $Q\in
\mathcal{Q}_{x_0}$.

\textit{Step} 2: By Skorohod's representation theorem, there exist a
probability space $(\tilde{\Omega},\tilde{\mathcal{B}},\tilde{P})$ and
$\mathbb{S}$-valued random variable $\tilde{x}_m$ and $\tilde{x}$
such that:
\begin{enumerate}[(ii)]
\item[(i)] $\tilde{x}_m$ has the law $Q_m$ for each $m\in\mathbb{N}$;
\item[(ii)]
$\tilde{x}_m\rightarrow\tilde{x}$ in $\mathbb{S}$, $\tilde
{P}$-a.e., and $\tilde{x}$ has the law $Q$.
\end{enumerate}

First, we easily deduce that
\begin{eqnarray*}
Q\bigl(x\bigl(t_1^n\bigr)=x_0
\bigl(t_1^n\bigr)\bigr)&=&\tilde{P}\bigl(\tilde{x}
\bigl(t_1^n\bigr)=x_0\bigl(t_1^n
\bigr)\bigr)
\\
&=&\lim_{m\rightarrow\infty} Q_m\bigl(x
\bigl(t_1^n\bigr)=x_m\bigl(t_1^n
\bigr)\bigr)=1,
\\
Q\bigl(x(t)=x\bigl(t_2^n\bigr),t\geq
t_2^n\bigr)&=&\tilde{P}\bigl(\tilde{x}(t)=\tilde {x}
\bigl(t_2^n\bigr),t\geq t_2^n
\bigr)
\\
&=&\lim_{m\rightarrow\infty} Q_m\bigl(x(t)=x
\bigl(t_2^n\bigr),t\geq t_2^n
\bigr)=1.
\end{eqnarray*}
For $q\in\mathbb{N}$, set
\[
\xi_q(x):=\sup_{r\in[t_1^n,t_2^n]}\bigl \|x(r)\bigr \|_{H^3}^{2q}+
\int_{t_1^n}^{t_2^n}\bigl \|x(r)\bigr \|_{H^3}^{2(q-1)}
\bigl \|x(r)\bigr \|_{H^{3+\alpha}}^2\,dr.
\]
Then
\begin{eqnarray*}
E^Q\bigl(\xi_q(x)\bigr)&=&E^{\tilde{P}}\bigl(
\xi_q(\tilde{x})\bigr)\leq\liminf_{m\rightarrow
\infty}
E^{Q_m}\bigl(\xi_q(x)\bigr)\leq\liminf_{m\rightarrow\infty}C
\bigl(\bigl \| x_m\bigl(t_1^n\bigr)
\bigr \|_{H^3}^{2q}+1\bigr)
\\
&\leq& C\bigl(\bigl \|x_0
\bigl(t_1^n\bigr)\bigr \|_{H^3}^{2q}+1
\bigr).
\end{eqnarray*}
Thus, (M1) and (M3) follow.

Now we want to show that $(M_i(t,x))_{ t\geq t_1^n}$ in the proof of
Theorem~\ref{the3.2} (Step 2) is a continuous $\mathcal
{B}_t$-martingale with
respect to $Q$, whose square variation process is given by
\[
\langle M_i\rangle(t,x)=\int_{t_1^n}^{t\wedge t_2^n}
\bigl \|(k_{\delta
_n}*G)^*\bigl(x(s)\bigr) (e_i)\bigr \|_U^2
\,ds.
\]
Since $\sup_{0\leq t\leq t_1^n}\|x_m(t)-x_0(t)\|_{H^3}\rightarrow0$ and
$\tilde{x}_m\rightarrow\tilde{x}$ in $\mathbb{S}$, we have
\begin{eqnarray*}
&&\lim_{m\rightarrow\infty}E^{\tilde{P}}\int
^{t^n_2}_{t_1^n}\bigl |\bigl\langle U_{\delta_n}[x_m](s)
\cdot\nabla\tilde {x}_m(s)+A_\alpha\tilde{x}_m(s)
\\
&&\phantom{\hspace*{64pt}}{}-U_{\delta_n}[x_0](s)
\cdot\nabla \tilde {x}(s)-A_\alpha\tilde{x}(s),e_i \bigr
\rangle\bigr | \,ds
\\
&&\qquad\leq\lim_{m\rightarrow\infty}E^{\tilde{P}}\int
^{t^n_2}_{t_1^n}\bigl |\bigl\langle\bigl(U_{\delta_n}[x_m](s)-U_{\delta
_n}[x_0](s)
\bigr)\cdot\nabla\tilde{x}_m(s)
\\
&&\phantom{\hspace*{98pt}}{}+U_{\delta_n}[x_0](s)
\cdot \nabla\bigl(\tilde{x}_m(s)-\tilde{x}(s)\bigr)
\\
&&\phantom{\hspace*{151pt}}{}+A_\alpha\bigl(\tilde{x}_m(s)-\tilde {x}(s)
\bigr),e_i \bigr\rangle\bigr | \,ds
\\
&&\qquad=0, %\end{array}
\end{eqnarray*}
which implies that for $t\geq t_1^n$
\begin{equation}
\label{eqB.4} \lim_{m\rightarrow\infty} E^{\tilde
{P}}\bigl |M_i^m(t,
\tilde{x}_m)-M_i(t,\tilde{x})\bigr |=0. %
\end{equation}
Then we obtain for $t_1^n\leq s<t$,
\[
E^Q\bigl(M_i(t,x)|\mathcal{B}_s
\bigr)=M_i(s,x). %
\]
On the other hand, by the B--D--G inequality, we have
\[
\sup_m E^{\tilde{P}}\bigl |M_i^m(t,
\tilde{x}_m)\bigr |^{2q}\leq C\sup_m\int
_{t_1^n}^{t_2^n}E^{\tilde{P}}\bigl(
\bigl \|(k_{\delta_n}*G)^*\bigl(\tilde {x}_m(s)\bigr)
(e_i)\bigr \|_{U}^{2q}\bigr)\,ds<+\infty.
\]
By (\ref{eqB.4}), we have
\[
\lim_{m\rightarrow\infty} E^{\tilde{P}}\bigl |M_i(t,
\tilde {x}_m)-M_i(t,\tilde {x})\bigr |^2=0.
\]
Then we obtain
\begin{eqnarray*}
&&E^Q\biggl(M_i^2(t,x)-\int
_{t_1^n}^t\bigl \|(k_{\delta_n}*G)^*\bigl(x(r)\bigr)
(e_i)\bigr \| _{U}^2\,dr|\mathcal{B}_s
\biggr)
\\
&&\qquad=M_i^2(s,x)-\int_{t_1^n}^s
\bigl \|(k_{\delta_n}*G)^*\bigl(x(r)\bigr) (e_i)\bigr \|_{U}^2
\,dr. %
\end{eqnarray*}
Now the results follow.
\end{pf}
%
%$\hfill\Box$

\setcounter{section}{2}
\section{Markov selections in the general case}\label{appC}

In this appendix, we will use \cite{GRZ09}, Theorem~4.7, to get an
almost sure
Markov family $(P_{x})_{x\in L^2}$ for equation (\ref{eq3.1}). Here,
we will
use the same notation as in \cite{GRZ09}. Below we choose
\[
H=\mathbb{Y}=L^2\bigl(\mathbb{T}^2\bigr)
\]
and
\[
\mathbb{X}=\bigl(H^{2+2\alpha}\bigr)^*,\qquad
\mathbb{X}^*=H^{2+2\alpha}. %
\]
Then $\mathbb{X}$ is a Hilbert space and $\mathbb{X}^*\subset\mathbb
{Y}$ compactly. Let $\mathcal{E}=\{e_i,i\in\mathbb{N}\}$ be the
orthonormal basis of $H$ introduced in Section~\ref{sec2}.
We define the operator $\mathcal{A}$ as follows: for $\theta\in
C^\infty
(\mathbb{T}^2)$
\[
\mathcal{A}(\theta):=-\kappa(-\Delta)^\alpha\theta-u\cdot
\nabla \theta, %
\]
where $u$ satisfies (\ref{eq1.3}). Then by Lemma~\ref{lemC.3}
below, $\mathcal{A}$ can
be extended to an operator $\mathcal{A}\dvtx H\rightarrow\mathbb
{X}$. For
$\theta$ not in $H$ define $\mathcal{A}(\theta):=\infty$.

Set
\[
\Omega:=C\bigl([0,\infty);\mathbb{X}\bigr), %
\]
and let $\mathcal{B}$ denote the $\sigma$-field of Borel sets of
$\Omega
$ and let $\mathcal{P}(\Omega)$ denote the set of all probability
measures on $(\Omega,\mathcal{B})$. Define the canonical process
$x\dvtx\Omega\rightarrow\mathbb{X}$ as
\[
x_t(\omega)=\omega(t). %
\]
For each $t$, $\mathcal{B} _t=\sigma(x_s\dvtx0\leq s\leq t)$. Given
$P\in
\mathcal{P}(\Omega)$ and $t>0$, let $P(\cdot|\mathcal{B}_t )(\omega
)$ denote
a regular conditional probability distribution of $P$ given $\mathcal
{B}_t $. In particular, $P(\cdot|\mathcal{B}_t )(\omega)\in\mathcal
{P}(\Omega)$ for every $\omega\in\Omega$ and for any bounded
$\mathcal
{B}$-measurable function $f$ on $\Omega$
\[
E^P[f|\mathcal{B}_t]=\int
_\Omega f(y)P(dy|\mathcal{B}_t),\qquad P
\mbox{-a.s.}, %
\]
and there exists a $P$-null set $N\in\mathcal{B}_t$ such that for every
$\omega$ not in $N$
\[
P(\cdot|\mathcal{B}_t) (\omega)_{|\mathcal{B}_t}=
\delta_\omega \qquad(=\mbox{Dirac measure at }\omega), %
\]
hence,
\[
P\bigl(\bigl\{y\dvtx y(s)=\omega(s),s\in[0,t]\bigr\}|
\mathcal{B}_t\bigr) (\omega)=1. %
\]
In particular, we can consider $P(\cdot|\mathcal{B}_t )(\omega)$ as a
measure on $(\Omega^t,\mathcal{B}^t)$, that is,
\[
P(\cdot|\mathcal{B}_t ) (\omega)\in\mathcal{P}\bigl(
\Omega^t\bigr), %
\]
where $\Omega^t:=C([t,\infty);\mathbb{X})$ and $\mathcal
{B}^t:=\sigma
(x_s\dvtx s\geq t)$.

We say $P\in\mathcal{P}(\Omega)$ is concentrated on the paths with
values in
$H$, if there exists $A\in\mathcal{B} $ with $P(A)=1$ such that
$A\subset\{\omega\in\Omega\dvtx x_t(\omega)\in H,\forall t\geq0\}$.
The set of such measures is denoted by $\mathcal{P}_H(\Omega)$.
The shift operator $\Phi_t\dvtx\Omega\rightarrow\Omega^t$ is
defined by
\[
\Phi_t(\omega) (s)=\omega(s-t), \qquad s\geq t.
\]

Following \cite{GRZ09}, Definitions 2.5, we introduce the following notions.

\begin{definition}\label{defC.1}
A family $(P_x)_{x\in H}$ of probability measures
in $\mathcal{P}_H(\Omega)$, is called an \emph{almost sure Markov
family} if for any $A\in\mathcal{B} $, $x\mapsto P_x(A)$ is $\mathcal
{B}(H)/\mathcal{B}([0,1])$-measurable, and for each $x\in H$ there
exists a Lebesgue null set $T_{P_x}\subset(0,\infty)$ such that for
all $t$ not in $T_{P_x}$ and $P_x$-almost all $\omega\in\Omega$
\[
P_x(\cdot|\mathcal{B}_t ) (
\omega)=P_{\omega(t)}\circ\Phi^{-1}_t. %
\]
\end{definition}

We now introduce the following notion of a martingale solution to
equation (\ref{eq3.1}) and write $x(t)$ instead of $x_t$.

\begin{definition}\label{defC.2}
Let $x_0\in H$. A probability measure $P\in\mathcal
{P}(\Omega)$ is called a martingale solution of equation
(\ref{eq3.1}) with initial
value $x_0$, if:
\begin{enumerate}[(M3)]
\item[(M1)]$P(x(0)=x_0)=1$ and for any $n\in\mathbb{N}$
\[
P\biggl\{x\in\Omega\dvtx \int_0^n
\bigl \|\mathcal{A}\bigl(x(s)\bigr)\bigr \|_{\mathbb{X}}\, ds+\int_0^n
\bigl \| G\bigl(x(s)\bigr)\bigr \|^2_{L_2(U;H)}\,ds<+\infty\biggr\}=1;
\]
\item[(M2)]for every $l\in\mathcal{E}$, the process
\[
M_l(t,x):=_\mathbb{X}\bigl\langle x(t),l\bigr
\rangle_{\mathbb{X}^*}-\int^t_0{ }
_\mathbb{X}\bigl\langle\mathcal{A}\bigl(x(s)\bigr),l\bigr
\rangle_{\mathbb{X}^*}\,ds %
\]
is a continuous square-integrable $\mathcal{B}_t$-martingale under $P$,
whose quadratic variation process is
given by
\[
\langle M_l\rangle(t,x):=\int_0^t
\bigl \|G^*\bigl(x(s)\bigr) (l)\bigr \|_U^2\,ds, %
\]
where the asterisk denotes the adjoint operator of $G(x(s))$;
\item[(M3)]for any $p\in\mathbb{N}$, there exist a continuous positive
real function $t\mapsto C_{t,p}$ (only depending on $p$ and $\mathcal
{A}, G$), a lower semicontinuous positive real functional $\mathcal
{N}_p\dvtx\mathbb{Y}\rightarrow[0,\infty]$, and a Lebesgue null set
$T_P\subset(0,\infty)$ such that for all $0\leq s\in[0,\infty
)\backslash T_P$ and for all $t\geq s$
\[
E^P\biggl[\sup_{r\in[s,t]}\bigl |x(r)\bigr |^{2p}+
\int_s^t\mathcal{N}_p\bigl(x(r)
\bigr)\, dr\Big|\mathcal {B}_s \biggr]\leq C_{t-s}
\bigl(\bigl |x(s)\bigr |^{2p}+1\bigr). %
\]
\end{enumerate}
\end{definition}

First, we prove the following lemma.

\renewcommand{\thelemma}{C.3}
\begin{lemma}\label{lemC.3}
For any $\theta_1,\theta_2\in C^\infty(\mathbb{T}^2)$,
\begin{eqnarray*}
\bigl \|(-\Delta)^\alpha\theta_1-(-
\Delta)^\alpha\theta_2\bigr \|_{\mathbb
{X}}&\leq& C_1|
\theta_1-\theta_2|, %
\\
\|u_1\cdot\nabla\theta_1-u_2
\cdot\nabla\theta_2\|_{\mathbb
{X}}&\leq& C_2\bigl(|
\theta_1|+|\theta_2|\bigr)|\theta_1-
\theta_2| %
\end{eqnarray*}
for constants $C_1,C_2$. In particular, the operator $\mathcal{A}\dvtx
C^\infty(\mathbb{T}^2)\rightarrow\mathbb{X}$ extends to an operator
$\mathcal{A}\dvtx H\rightarrow\mathbb{X}$ by continuity.
\end{lemma}

\begin{pf}
We only prove the second assertion, the first can be proved
analogously. By the Sobolev embedding theorem, we have
\begin{eqnarray*}
&&\|u_1\cdot\nabla\theta_1-u_2\cdot\nabla
\theta_2\|_{\mathbb
{X}}
\\
&&\qquad=\sup_{w\in C^\infty(\mathbb{T}^2):\|w\|_{H^{2+2\alpha
}}\leq
1}\bigl |\langle u_1\cdot\nabla
\theta_1-u_2\cdot\nabla\theta_2,w\rangle\bigr |
\\
&&\qquad=\sup_{w\in C^\infty(\mathbb{T}^2):\|w\|_{H^{2+2\alpha
}}\leq1}\bigl |\langle u_1\cdot\nabla w,
\theta_1\rangle-\langle u_2\cdot\nabla w,\theta
_2\rangle\bigr |
\\
&&\qquad=\sup_{w\in C^\infty(\mathbb{T}^2):\|w\|_{H^{2+2\alpha
}}\leq1}\bigl |\bigl\langle(u_1-u_2)
\cdot\nabla w,\theta_1\bigr\rangle+\langle u_2\cdot \nabla
w,\theta_1-\theta_2\rangle\bigr |
\\
&&\qquad\leq C\Bigl[\sup_{w\in C^\infty(\mathbb
{T}^2):\|w\|_{H^{2+2\alpha}}\leq1}\|\nabla w\|_{C(\mathbb
{T}^2)}
\Bigr]\bigl(|u_1-u_2|\cdot|\theta_1|+|
\theta_1-\theta_2|\cdot|u_2|\bigr)
\\
&&\qquad\leq C\bigl(|\theta_1|+|\theta_2|\bigr)|
\theta_1-\theta_2|.
\end{eqnarray*}
In the last inequality, we use (\ref{eq2.1}) and the constant $C$
changes from
line to line.\end{pf}
%
%$\hfill\Box$

In order to use \cite{GRZ09}, Theorem~4.7, we define the functional
$\mathcal
{N}_1$ on $\mathbb{Y}$ as follows:
\[
\mathcal{N}_1(\theta):=\cases{ %
\bigl |\Lambda^\alpha\theta\bigr |^2,&\quad$\mbox{if }
\theta\in H^{\alpha}$,
\cr
+\infty,&\quad$\mbox{otherwise}$. %\end{array}
} %
\]
It is obvious that $\mathcal{N}_1\in\mathfrak{U}^2$, defined in
\cite{GRZ09},
Section~4. We recall that a lower semicontinuous function $\mathcal
{N}\dvtx\mathbb{Y}\rightarrow[0,\infty]$ belongs to $\mathfrak
{U}^2$ if
$\mathcal{N}(x)=0$ implies $x=0$, $\mathcal{N}(cy)\leq c^2\mathcal
{N}(y), \forall c\geq0,y\in\mathbb{Y}$ and $\{y\in\mathbb{Y}\dvtx
\mathcal
{N}(y)\leq1\}$ is relatively compact in~$\mathbb{Y}$.

\renewcommand{\thetheorem}{C.4}
\begin{theorem}\label{theC.4}
Let $\alpha\in(0,1)$ and assume $G$ satisfies Hypothesis~\ref{hypG.1}
with $\rho_1=0$. Then for each $x_0\in H$, there exists a martingale
solution $P\in\mathcal{P}(\Omega)$ starting from $x_0$ to equation
(\ref{eq3.1}) in the sense of Definition~\ref{defC.2}.
\end{theorem}

\begin{pf}
We only need to check (C1)--(C3) in \cite{GRZ09}, Section~4, for the
above $\mathcal{A}$ and $G$.

The demi-continuity condition (C1) holds since Lemma~\ref{lemC.3}
and Hypothesis~\ref{hypG.1}
imply demi-continuity of $\mathcal{A}$ and $G$.

The coercivity condition (C2) follows, because noting that for $\theta
\in\mathbb{X}^*$
\[
\langle u\cdot\nabla\theta, \theta\rangle=0, %
\]
we have
\[
\bigl\langle\mathcal{A}(\theta),\theta\bigr\rangle=-
\mathcal{N}_1(\theta). %
\]

Also the growth condition (C3) is clear since by Lemma~\ref{lemC.3}
\[
\bigl \|\mathcal{A}(\theta)\bigr \|_{\mathbb{X}}\leq C|\theta|^2
\]
and
\[
\bigl \|G(\theta)\bigr \|_{L_2(K;H)}\leq C\bigl(|\theta|+1\bigr). %
\]
\upqed
\end{pf}
%
%$\hfill\Box$

The set of all such martingale solutions with initial value $x_0$ is
denoted by $\mathcal{C}(x_0)$. Using \cite{GRZ09}, Theorem~4.7, we now
obtain the following.

\renewcommand{\thetheorem}{C.5}
\begin{theorem}\label{theC.5}
Let $\alpha\in(0,1)$. Assume $G$ satisfies Hypothesis~\ref{hypG.1} with
\mbox{$\rho_1=0$}. Then there exists an almost sure Markov family
$(P_{x_0})_{x_0\in H}$ for equation (\ref{eq3.1}) and $P_{x_0}\in
\mathcal{C}(x_0)$
for each $x_0\in H$.
\end{theorem}
\end{appendix}


%suskaldyti doi

% imsref loaded by audrone.aklyte, 2014-02-27 09:02:28
%
% imsref loaded by audrone.aklyte, 2014-03-03 10:16:29

\printaddresses

\end{document}